\documentclass[10pt]{article}
\usepackage{amssymb}           
\usepackage{amsmath}           
\usepackage{amsthm}           
\usepackage{amsbsy}           
\usepackage[dvips]{graphics}           
\usepackage{amscd}           

\theoremstyle{plain}           
\newtheorem{theorem}{Theorem}[section]

\newtheorem{lemma}[theorem]{Lemma}           
\newtheorem{corollary}[theorem]{Corollary}           
\newtheorem{proposition}[theorem]{Proposition}

\theoremstyle{definition}           
\newtheorem{definition}{Definition}[section]           
\theoremstyle{remark}      
        
\newtheorem{remark}[theorem]{Remark}

\newcommand{\Z}{{\mathbb{Z}}}           
\newcommand{\C}{{\mathbb{C}}}

\newcommand{\s}{{\Sigma}_{0,\infty}}

\newcommand{\M}{{{\cal M}}}           
           
\newcommand{\B}{{{\cal B}}}

\newcommand{\D}{{{\cal D}_0^s}}

\def\longlongrightarrow{\relbar\joinrel\relbar\joinrel\relbar\joinrel\rightarrow}      

\def\C{\mathbb{C}}           
\def\Z{\mathbb{Z}}           
\def\Q{\mathbb{Q}}

\DeclareMathOperator{\Mor}{Mor}      
\DeclareMathOperator{\Gr}{Gr}          
           
\begin{document}

\title{On a universal mapping class group of genus zero\\
{\small \em  GAFA 2004, to appear}}          
 \author{           
\begin{tabular}{cc}           
 Louis Funar &  Christophe Kapoudjian\\           
\small \em Institut Fourier BP 74, UMR 5582            
&\small \em Laboratoire Emile Picard, UMR 5580\\           
\small \em University of Grenoble I &\small \em University of Toulouse           
III\\           
\small \em 38402 Saint-Martin-d'H\`eres cedex, France           
&\small \em 31062 Toulouse cedex 4, France\\           
\small \em e-mail: {\tt funar@fourier.ujf-grenoble.fr}           
& \small \em e-mail: {\tt ckapoudj@picard.ups-tlse.fr} \\           
\end{tabular}           
}

\date{}
\maketitle           
       
\begin{abstract}      
The aim of this paper is to introduce a group containing        
the mapping class groups of all genus zero surfaces.       
Roughly speaking, such a group is intended       
to be a discrete analogue of the diffeomorphism group of the circle. One      
defines indeed a {\it universal mapping class group of genus      
zero}, denoted $\B$. The latter is a nontrivial extension of the Thompson      
group $V$ (acting on the Cantor set) by an inductive limit of pure mapping class groups      
of all genus zero surfaces. We prove that $\B$ is a finitely presented group, and give an explicit presentation of it.      
      
\vspace{0.2cm}      
      
      
      
\end{abstract}

\section{Introduction}           
The problem of considering all mapping class groups together has been           
considered first by Moore and Seiberg (\cite{MS}), who worked with           
the somewhat imprecise {\em duality groupoid},  an essential  ingredient           
in their definition of conformal field theories in dimension two.              
Rigorous proofs of their results were obtained first by K.Walker           
(\cite{W}), who worked out the axioms of a topological quantum           
field theory with corners in dimension three, and were further improved and given a definitive treatment in \cite{BK,BK2,FG,HLS}          
from quite different perspectives. The answer provided in these papers          
is a groupoid containing all mapping class groups and having a           
finite (explicit) presentation, in which generators come from           
surfaces with $\chi=-1$ and relations from surfaces of $\chi=-2$,           
illustrating the so-called Grothendieck principle. The extra          
structure one considers in the tower of mapping class groups (of          
surfaces with boundary) is the exterior multiplication law           
(inducing a monoidal category structure) coming from gluing together           
surfaces along boundary components. In some sense the duality           
groupoid is the smallest groupoid in which the tower of mapping class          
groups and the exterior multiplication law fit together.           
       
\vspace{0.1cm}       
\noindent        
Nevertheless this answer is not completely satisfactory. We would like          
to obtain a universal mapping class group, whose category of          
representations corresponds to the (groupoid) representations of the       
duality          
groupoid. This group would be a discrete analogue of the           
diffeomorphism group of the circle. This analogy           
suggests us to look for a connection between the Thompson group           
and the mapping class groups. We fulfill this program       
for the genus zero situation in the present paper, and will give a partial treatment for           
the full tower (of arbitrary genus) surfaces in a forthcoming article.           
          
\vspace{0.1cm}       
\noindent            
We introduce first the {\em universal mapping class group}            
of genus zero --- denoted $\B$ --- by a geometric construction (see Section 2): the elements of      
$\B$ are mapping classes of a certain surface $\s$ of genus zero,    
which is homeomorphic      
to a sphere minus a Cantor set. The mapping classes are assumed to preserve      
asymptotically       
some extra structure on $\s$, called the {\em rigid structure} of $\s$. The defining      
properties of the resulting group $\B$ are: first it contains            
{\em uniformly} all mapping class groups of holed           
spheres and second, it surjects onto the Thompson group $V$.

\begin{definition}           
Denote by $K^*(g,n)={K^*}(\Sigma_{g,n})$  the pure mapping   
class group of the     
$n$-holed orientable surface of genus $g$, which            
consists of classes of orientation preserving    
homeomorphisms of $\Sigma_{g,n}$, modulo isotopies which are pointwise fixing   
the boundary. Denote by $\M(g,n)={\M}(\Sigma_{g,n})$ the full   
mapping class group, consisting of mapping classes of homeomorphisms   
which respect a fixed parametrization of the boundary circles,  
allowing them to be permuted among themselves. This group is   
related to $K^*(g,n)$ by the short exact sequence  
\[1\to K^*(g,n)\to \M(g,n) \to {\cal S}_n\to 1\]   
where ${\cal S}_n$ stands for the permutation group on $n$ elements.   
\end{definition}

\noindent When $g=0$ and  $n$ goes to infinity, this  short   
exact sequence        
stabilizes to give rise to the exact sequence           
\[1\to  K^*_{\infty}\to {\mathcal B} \to V\to 1 \]           
where $V$ is the Thompson group acting on the Cantor set (see \cite{CFP}),      
$K^*_{\infty}$ is an inductive limit $\bigcup_{n} K^*(0,3\cdot 2^n)$, and ${\mathcal B}$           
is the universal mapping class group of genus zero.             
The inductive limit takes into account a suitable            
natural injection $K^*(0,m)\hookrightarrow  K^*(0,2m)$ which will            
become obvious in the sequel.

\vspace{0.1cm}              
\noindent The main result of this paper is (see Section 3, Theorem \ref{pres0}          
for a more precise statement and explicit relations) as follows:

\begin{theorem}\label{main}          
The group $\B$ is finitely presented.           
\end{theorem}          
          
\noindent This is somewhat unexpected since $\B$ contains all the genus zero mapping class          
groups, and the kernel of the surjection to          
$V$ is an infinitely generated group. We will discuss  first the relationship between the group $\B$ and the         
the duality groupoid considered in \cite{FG} (see Section 4), then        
present the proof of the main result in  Section 5. The method used to obtain the presentation is greatly inspired from           
Hatcher-Thurston's approach (\cite{HT}) to the mapping class groups of    
compact surfaces, insofar as it exploits the action of $\B$ on the simply    
connected Hatcher-Thurston complex ${\cal HT}(\s)$ of the    
(non-compact) surface $\s$. However, the finiteness of the presentation is    
somewhat miraculous, and relies on a deep connection between the Hatcher-Thurston      
complex ${\cal HT}(\s)$ and the Ptolemy-Thompson group $T$, the subgroup      
of $V$ acting on the circle. Indeed, the    
presentation of $T$ enables us to find a simply connected subcomplex of ${\cal    
  HT}(\s)$, which has only a finite number of orbits of 2-cells under the action of    
$\B$.   
  
\vspace{0.1cm}              
\noindent   
It is likely that, by strengthening the methods of the present paper,  
one is able to prove that the group $\B$ is    
${\rm FP}_{\infty}$.

\vspace{0.1cm}              
\noindent   
In the last section (Section 6), we investigate the relationship between the    
group $\B$ and a group recently discovered by M. Brin (see   
\cite{bri}), of which we became aware upon the completion of our   
work. This is the {\it braided Thompson group} $BV$, an    
extension of $V$ by an inductive limit of Artin pure braid groups. The   
relation between $\B$ and $BV$ is the same as that between the   
mapping class group of the holed sphere and the usual braid group, and   
we prove therefore that $\B$ contains $BV$.     
   
\vspace{0.1cm}              
\noindent   
The original motivation in \cite{bri}  for considering $BV$ was    
the intimate relationship between   coherence questions in categories   
with multiplication  and Thompson's groups. Specifically, let us    
consider a category    
${\mathcal C}$ with functorial multiplication $\otimes$, identity element, natural     
isomorphisms $\alpha: (A\otimes B)\otimes C\to A\otimes(B\otimes C)$,    
$c:A\otimes B\to B\otimes A$ and $t:A\to A$. M. Brin (\cite{bri})   
constructed groups  and epimorphisms:    
$\xi_2:G_2({\mathcal C}, \otimes, \alpha, c)\to V$,    
$\xi_3:G_3({\mathcal C}, \otimes, \alpha, c)\to BV$ such that       
$({\mathcal C},\otimes, \alpha, c)$ is a {\em symmetric, monoidal   
category} if and only if $\xi_2$ is an isomorphism, and that    
$({\mathcal C},\otimes, \alpha, c)$ is a {\em braided, tensor category} if   
and only if $\xi_3$ is an isomorphism. Along these lines one can   
construct (but this is beyond the scope of this paper) a group and an epimorphism   
$\xi_4:G_4({\mathcal C}, \otimes, \alpha, c, t)\to \B$, with the   
property that $({\mathcal C},\otimes, \alpha, c, t)$ is a {\em ribbon   
category}  if and only if $\xi_4$ is an isomorphism.

\vspace{0.1cm}              
\noindent   
It should    
be mentioned that the link between (some of) Thompson's groups and the    
braid groups has been revealed for the first time in the work of    
P. Greenberg and V. Sergiescu (\cite{gr-se}), where    
they construct an acyclic extension ${\cal A}$ of $F'$ --- the derived subgroup    
of Thompson's group $F\subset T$ --- by the stable braid group $B_{\infty}$. However, the new approach to the group ${\cal    
  A}$ given in    
\cite{ka-se} clarifies the differences between this group and    
the groups $\B$ and    
$BV$:  ${\cal  A}$ has a description as a mapping class group braiding    
a countable family of punctures on a certain non-compact surface (obtained    
from $\s$ by adding some tubes with punctures, see    
\cite{ka-se} for the details), while $\B$ and    
$BV$ rather braid the ends at infinity of the surface $\s$.

\vspace{0.1cm}       
\noindent Another place where the duality groupoid and the tower of mapping          
class groups enter as a key ingredient is in Grothendieck's approach           
to ${\rm Gal}(\overline{\Q}/\Q)$, as  further developed by Ihara,    
Drinfeld and  explained and explored in a series of papers    
by Lochak and Schneps (see          
e.g. \cite{HLS,LS1,LS}). One would like to understand the           
possible equalities between the following groups:           
\[{\rm Gal}(\overline{\Q}/\Q)\subset \widehat{\rm GT}\subset {\rm Out}({\cal  W}), \]          
where $\widehat{\rm GT}$           
is the group of Grothendieck-Teichm\"uller          
group (as introduced by Drinfeld in \cite{D}), and           
${\rm Out}({\cal W})$ is the outer automorphism group of the           
 (various) towers ${\cal W}$ of profinite completions    
of mapping class groups.     
One version of ${\cal W}$ consists of the genus zero surfaces,          
$\{\M(0,n)\}_{n\geq 3}$ and the gluing homomorphisms,          
while another version consists of all surfaces.           
          
\vspace{0.1cm}   
\noindent It is known (see \cite{LS1}) that $\widehat{\rm GT}$           
acts naturally and faithfully on          
the profinite groups $\widehat{K^*(0,n)}$, respecting the    
(gluing surfaces)           
homomorphisms between these groups. In \cite{LNS, NS},           
the authors extended these results to higher genus.           
In the same way $\widehat{\rm GT}$ acts on a suitable   
completion $\widehat{\B}$ of the {\em group} $\B$ (cf. \cite{ka2}). The group  $\widehat{\B}$ is a relative profinite completion of $\B$ with   
respect to the morphism $\B\to V$. In fact, the extended   
Grothendieck-Teichm\"uller group $\widehat{\rm GT}_e$, which is an extension   
of $\widehat{\rm GT}$ by $\widehat{\Z}$, embeds into ${\rm   
  Out}(\widehat{\B})$. It is likely that the cokernel of the embedding   
$\widehat{\rm GT}_e \to {\rm Out}(\widehat{\B})$ is ${\rm  Out}(V)$. The computation of ${\rm Out}(V)$ (which one might    
reasonably conjecture to be $\Z/2\Z$) would permit   
to replace the tower of mapping class groups    
by a single group.

\vspace{0.3cm}

\noindent {\bf Acknowledgements.} We are indebted to Joan         
Birman, Matthew Brin, Thomas Fiedler and Bill Harvey for 
their comments and useful discussions. We are grateful to the referee for  
his/her careful reading of the paper and for his/her numerous  
comments which led to a considerable improvement of the   
accuracy and of the quality of the exposition.     
Significant impetus for writing this paper    
was provided by the stimulating discussions we had with    
Vlad Sergiescu, who introduced both of us to Thompson's groups.

\section{The construction}           
The main step in obtaining the universal mapping class group $\B$ is           
to shift from the compact surfaces to an infinite surface, and to           
consider those homeomorphisms having a nice behaviour at infinity.\\

\subsection{The genus zero infinite surface $\s$}           
\begin{definition}\label{s}           
Let $\s$ be the infinite surface $\s$ of genus zero, built up as           
an inductive limit of finite subsurfaces $S_m$, $m\geq 0$: $S_0$ is a 3-holed sphere, and $S_{m+1}$ is obtained from $S_m$ by gluing a copy of a 3-holed sphere along each boundary           
component of $S_m$. The surface $\s$ is oriented and all           
homeomorphisms considered in the sequel will be           
orientation preserving, unless the opposite is explicitly stated.           
\end{definition}

\begin{figure}      
\begin{center}           
\includegraphics{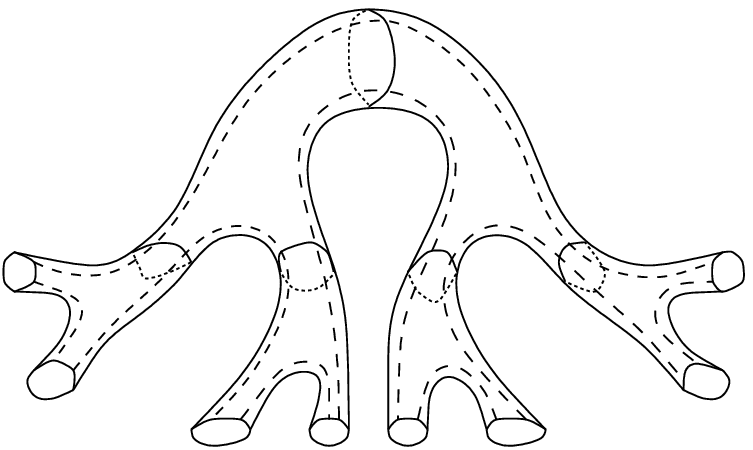}           
\caption{Canonical pants decomposition and rigid structure on $\s$}      
\end{center}           
\end{figure}           
           
\begin{definition}        
\begin{enumerate}           
\item A {\em pants decomposition} of the surface $\s$ is a maximal collection            
of distinct nontrivial simple closed curves on $\s$ which are pairwise            
disjoint and non-isotopic. The complementary regions (which are 3-holed           
spheres) are called {\em pairs of pants}.            
\item A {\em rigid structure} (see \cite{FG}) on $\s$ consists of two   
pieces of data:
\begin{itemize}   
\item a pants           
  decomposition, and      
\item  a {\em prerigid structure}, i.e. a countable collection of disjoint    
line segments embedded into $\s$, such that the complement of their   
union in $\s$ has two connected components.   
\end{itemize}    
These pieces must be {\em compatible} in the following sense:    
first, the traces of the prerigid    
structure on each pair of pants (i.e. the intersections with the pairs   
of pants) are made up of three connected          
components, called {\em seams}. Second,    
for each pair of boundary circles of a given      
  pair of pants, there is exactly one seam joining the two circles.   
    
One says then that the pants decomposition and the prerigid structure are {\em subordinate} to the rigid           
  structure.               
\item By construction, $\s$ is naturally equipped with a pants           
decomposition, which will be referred to below as the            
{\em canonical (pants) decomposition}. One fixes a prerigid structure      
(called the {\em canonical prerigid structure})    
compatible with the canonical decomposition (cf. Figure 1).    
The resulting rigid structure is the           
{\em canonical rigid structure}. Note that it is not canonically   
defined.       
\item  The complement in $\s$ of the union of lines      
 of the canonical prerigid structure has two components: we distinguish one of           
them as the {\em visible side} of $\s$.           
\item A pants decomposition (resp. (pre)rigid structure) is {\em asymptotically           
    trivial} if outside a compact subsurface of $\s$, it coincides with the           
  canonical pants decomposition (resp. canonical (pre)rigid structure).           
\end{enumerate}           
\end{definition}

\subsection{The universal mapping class group $\B$ of genus zero}           
      
\begin{definition}\label{groupB}           
\begin{enumerate}           
\item A compact subsurface $\Sigma_{0,n}\subset \s$ is {\em admissible} if           
its boundary is contained in the canonical decomposition. The {\em level} of a           
(not necessarily admissible) compact subsurface $\Sigma_{0,n}\subset \s$ is           
the number $n$ of its boundary components.              
\item Let $\varphi$ be a homeomorphism of   
  $\s$. One says that            
$\varphi$ is {\em asymptotically rigid} if there           
exists an admissible subsurface $\Sigma_{0,n}\subset \s$ such           
  that: $\varphi(\Sigma_{0,n})$ is also admissible, and  the    
restriction of $\varphi: \s -\Sigma_{0,n}\to \s-\varphi(\Sigma_{0,n})$           
 is {\em rigid}, meaning that it respects the traces of the    
canonical rigid structure, mapping the pants           
 decomposition into the pants decomposition, the seams into the           
 seams, and the visible side into the visible side. Such a surface   
 $\Sigma_{0,n}$ is called a {\em support} for $\varphi$.  Note that we   
 are not using the word ``support" in the usual sense, as the map outside the   
 support defined above might well not being the identity.    
   
\end{enumerate}           
\noindent It is easy to see that the isotopy classes of asymptotically    
rigid homeomorphisms form a group, which one    
denotes by $\B$,    
and which will be called the {\em universal mapping class group in           
genus zero}.            
\end{definition}

\begin{remark}\label{PSL} The subgroup of $\B$ consisting of those mapping classes represented by globally rigid homeomorphisms is isomorphic to the    
  group of automorphisms of the planar tree ${\cal T}$ of the surface (see    
  Definition \ref{arbre}) which respect the local orientation of the edges around    
  each vertex: this is the group $PSL(2,\Z)$. In the notation of    
  Section 3, this subgroup is freely generated by the elements $\alpha^2$ (of order 2) and    
  $\beta$ (of order 3).    
\end{remark}              
    
\begin{remark}\label{comp}           
Subsets of $\s$ have a visible side, which is the intersection with   
the visible side of $\s$, and a hidden side which is the   
complement of the former. In particular    
each boundary circle of an admissible subsurface           
$\Sigma_{0,n}\subset \s$ has a visible and a hidden side, which are    
both half-circles. The full mapping class group   
${\M}(\Sigma_{0,n})$ may be equivalently     
defined as the group of isotopy classes of orientation preserving    
homeomorphisms of $\Sigma_{0,n}$    
which permute the boundary components but preserve their visible sides. The    
isotopies stabilize the half-circles. There is an obvious injective morphism           
$i_*:{\M}(\Sigma_{0,n})\hookrightarrow \B$, obtained by extending rigidly a           
homeomorphism representing a mapping class of ${\M}(\Sigma_{0,n})$.   
   
\vspace{0.1cm}   
\noindent           
Furthermore, if $\Sigma$ and $\Sigma'$ are two admissible subsurfaces,           
we denote by ${\M}(\Sigma)\cap {\M}(\Sigma')$ the intersection in           
$\B$ of the natural images of ${\M}(\Sigma)$ and   
${\M}(\Sigma')$, i.e. the set of those mapping classes of           
${\M}(\Sigma\cap\Sigma' )$ which extend rigidly to both           
$\Sigma$ and $\Sigma'$. The compatibility of the embeddings of the           
various mapping class groups into $\B$ is summarized in Figure   
\ref{dia}.

\begin{figure}
\begin{center}           
\begin{picture}(0,0)%
\includegraphics{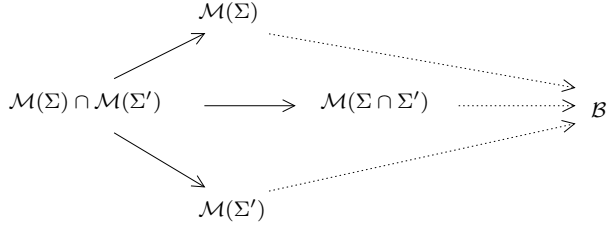}
\end{picture}%
\setlength{\unitlength}{2486sp}%
\begingroup\makeatletter\ifx\SetFigFont\undefined%
\gdef\SetFigFont#1#2#3#4#5{%
  \reset@font\fontsize{#1}{#2pt}%
  \fontfamily{#3}\fontseries{#4}\fontshape{#5}%
  \selectfont}%
\fi\endgroup%
\begin{picture}(5805,2194)(226,-1919)   
\put(6031,-871){\makebox(0,0)[lb]{\smash{\SetFigFont{8}{9.6}{\rmdefault}{\mddefault}{\updefault}${\cal B}$}}}   
\put(226,-781){\makebox(0,0)[lb]{\smash{\SetFigFont{8}{9.6}{\rmdefault}{\mddefault}{\updefault}${\cal M}(\Sigma)\cap {\cal M}(\Sigma')$}}}   
\put(2116,-1861){\makebox(0,0)[lb]{\smash{\SetFigFont{8}{9.6}{\rmdefault}{\mddefault}{\updefault}${\cal M}(\Sigma')$}}}   
\put(2116,119){\makebox(0,0)[lb]{\smash{\SetFigFont{8}{9.6}{\rmdefault}{\mddefault}{\updefault}${\cal M}(\Sigma)$}}}   
\put(3331,-781){\makebox(0,0)[lb]{\smash{\SetFigFont{8}{9.6}{\rmdefault}{\mddefault}{\updefault}${\cal M}(\Sigma\cap \Sigma')$}}}   
\end{picture}

\end{center}           
\caption{Diagram of embeddings into $\B$ of the various ${\M}(\Sigma)$}\label{dia}           
\end{figure}

\end{remark}           
           
\begin{definition}           
Let $\Sigma$ be an admissible subsurface, $K^*(\Sigma)$ its the pure mapping class group. Each inclusion           
$\Sigma\subset\Sigma'$ induces an injective embedding $j_{\Sigma,\Sigma'}:           
K^*(\Sigma)\rightarrow K^*(\Sigma')$ (though not always a morphism           
$\M(\Sigma)\rightarrow \M(\Sigma')$). The collection           
$\{K^*(\Sigma), j_{\Sigma,\Sigma'}\}$ is a direct system, and one denotes by           
$K^*_{\infty}$ its direct limit.     
\end{definition}

\subsection{The group $\B$ as an extension of Thompson's group $V$}           
\begin{definition}[Tree for $\s$]\label{arbre}           
\begin{enumerate}           
\item Let ${\cal T}$ be the infinite binary tree. There is a natural projection            
$q:\s\to {\cal T}$, such that the pullback of the set of edge midpoints           
is the set of circles of the canonical pants           
decomposition. The projection $q$ admits a continuous cross-section,           
that is, one may embed ${\cal T}$ in the visible side of $\s$, with   
one vertex on each pair of pants of          
the canonical pants decomposition, and one edge transverse to each circle.    
Since the  visible side of $\s$ is a planar surface, ${\cal T}$ will be viewed as          
a planar tree.           
\item A {\em finite binary tree} $X$ is a finite subtree of ${\cal T}$ whose           
  internal vertices are all 3-valent. Its terminal vertices (or 1-valent           
  vertices) are called {\em leaves}. One denotes by ${\cal L}(X)$ the set of           
  leaves of $X$, and calls the number of leaves the {\em level} of $X$.           
\end{enumerate}           
\end{definition}

\begin{figure}      
\begin{center}           
\includegraphics{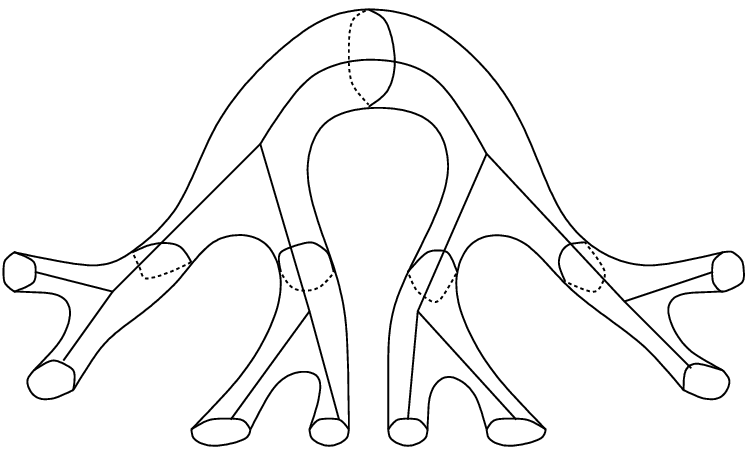}        
\caption{The tree ${\cal T}$ for  $\s$}      
\end{center}      
\end{figure}           
               
\begin{definition}[Thompson's group $V$]           
\begin{enumerate}           
\item A symbol $(T_1,T_0,\sigma)$ is a triple consisting of two finite           
  binary trees $T_0$ and $T_1$ of the same level, together with a bijection $\sigma:{ \cal L}(T_0)\rightarrow           
  { \cal L}(T_1)$.           
\item If $X$ is a finite binary subtree of ${\cal T}$ and $v$ is a leaf of           
  $X$, one defines the finite binary subtree $\partial_v X$ as the union of $X$           
  with the two edges which are the descendants of $v$. Viewing $\partial_v X$        
  as a subtree of the planar tree ${\cal T}$, one may distinguish the left descendant from the right descendant           
  of $v$. Accordingly, one denotes by $v_l$ and $v_r$ the leaves of the two new           
  edges of $\partial_v X$.            
\item  Let ${\cal R}$ be the equivalence relation on the set of symbols           
  generated by the following relations:           
$$(T_1,T_0,\sigma)\sim_v (\partial_{\sigma(v)}T_1,\partial_v T_0,\partial_v           
\sigma)$$           
where $v$ is any leaf of $T_0$, and $\partial_v \sigma$ is the natural           
extension of $\sigma$ to a bijection    
${\cal L}(\partial_v T_0)\rightarrow {\cal L}(\partial_{\sigma(v)}   
T_1)$ which maps $v_l$ and  $v_r$ to $\sigma(v)_l$ and           
$\sigma(v)_r$, respectively. One denotes by $[T_1,T_0,\sigma]$ the class of a           
symbol $(T_1,T_0,\sigma)$, and by $V$ the set of equivalence classes for           
the relation ${\cal R}$. Given two elements of $V$, one may represent them by           
two symbols of the form $(T_1,T_0,\sigma)$ and $(T_2,T_1,\tau)$ respectively, and define           
the product           
$$[T_2,T_1,\tau]\cdot[T_1,T_0,\sigma]=[T_2,T_0,\tau\circ\sigma]$$           
This product endows $V$ with a group structure,    
with neutral element $[T,T,id_{{\cal           
    L}(T)}]$, where $T$ is any finite binary subtree. The present   
group $V$ is Thompson's group $V$ (cf. \cite{CFP}).              
\end{enumerate}           
\end{definition}           
      
\begin{remark}      
We warn the reader that our definition of the group $V$ is slightly different      
from the standard one (as given in \cite{CFP}), since the binary trees occurring      
in the symbols are {\it unrooted} trees. Nevertheless, the present group $V$ is isomorphic to      
the group denoted by the same letter in \cite{CFP}.       
\end{remark}       
      
\noindent We introduce Thompson's group $T$, the subgroup of $V$ acting on the      
circle (see \cite{GS}), which will play a key role in the proof of Theorem \ref{pres0}.          
          
\begin{definition}[Ptolemy-Thompson's group $T$]\label{tho}          
\begin{enumerate}          
\item Let $T_{S_0}$ be the smallest finite binary subtree          
  of ${\cal T}$ containing ${q}(S_0)$. Choose a cyclic (counterclockwise          
  oriented, with respect to the orientation of $\s$)          
  labeling of its leaves by $1,2,3$. Extend inductively this cyclic labeling          
  to a cyclic labeling by $\{1,...,n\}$ of the leaves of any finite    
  binary subtree of $\,{\cal T}$          
  containing $T_{S_0}$, in the following way: if $T_{S_0}\subset X\subset \partial_v X$, where $v$          
  is a leaf of a cyclically labeled tree $X$, then there is a unique cyclic          
  labeling of the leaves of $\partial_v X$ such that:          
\begin{itemize}          
\item if $v$ is not the leaf 1 of $X$, then the leaf 1 of $\partial_v X$    
  coincides with the leaf 1 of $X$;          
\item if $v$ is the leaf 1 of $X$, then the leaf 1 of $\partial_v X$ is the          
  left descendant of $v$.          
\end{itemize}          
\item Thompson's group $T$ (also called Ptolemy-Thompson's group) is   
  the subgroup of $\,V$ consisting of elements          
  represented by symbols $(T_1,T_0,\sigma)$, where $T_1$ and $T_0$ contain          
  ${q}(S_0)$, and $\sigma: { \cal L}(T_0)\rightarrow  { \cal L}(T_1)$ is a cyclic permutation. The cyclicity of        
  $\sigma$ means that there exists some integer $i_0$, $1\leq i_0\leq n$ (if        
  $n$ is the level of $T_0$ and $T_1$), such that $\sigma$ maps the $i^{th}$        
  leaf of $T_0$ onto the $(i+i_0)^{th}$ (mod $n$) leaf of $T_1$, for $i=1,\ldots, n$.          
\end{enumerate}          
\end{definition}            
           
\begin{proposition}\label{sequence}           
We have the following exact sequence:           
\[ 1 \to K^*_{\infty} \to \B \to V\to 1\]            
Moreover this extension splits over the Ptolemy-Thompson group   
$T\subset V$ i.e. there exists a section  $T\hookrightarrow \B$.     
\end{proposition}           
           
\begin{proof}           
Let us define the projection $\B \to V$.            
Consider $\varphi\in \B$ and let $\Sigma$ be a support for           
$\varphi$. We introduce the           
  symbol           
  $(T_{\varphi(\Sigma)},T_{\Sigma},\sigma(\varphi))$, where    
$T_{\Sigma}$ (resp. $T_{\varphi(\Sigma)}$) denotes the minimal finite binary subtree of ${\cal T}$ which contains ${q}(\Sigma)$ (resp. ${q}(\varphi(\Sigma))$), and $\sigma(\varphi)$ is the bijection induced by $\varphi$ between the set of leaves of both trees. The image of $\varphi$ in $V$ is the            
class of this triple, and it is easy to check that this correspondence induces           
a well-defined and surjective morphism $\B\rightarrow V$. The kernel is the subgroup of           
isotopy classes of homeomorphisms inducing the identity outside a   
support,   
and hence is the direct limit of the pure mapping class groups.   
   
\vspace{0.1cm}    
\noindent          
Denote by ${\bf T}$ the subgroup of $\B$ consisting of mapping classes           
represented by asymptotically rigid homeomorphisms preserving the whole visible           
side  of $\s$. The image of ${\bf T}$ in $V$ is the subgroup of elements           
represented by symbols $(T_1,T_0,\sigma)$, where $\sigma$ is a bijection preserving the cyclic order of the labeling of the leaves of the   
trees. Thus, the image of ${\bf T}$ is Ptolemy-Thompson's group    
$T\subset V$. Finally, the kernel of the epimorphism ${\bf T}\to T$ is           
trivial. In the following, we shall identify $T$ with ${\bf T}$.                
\end{proof}

\subsection{Universality of the group $\B$}          
We will show below that $\B$ is, in some sense, the smallest           
group containing uniformly all genus zero mapping class groups.   
   
\vspace{0.1cm}           
\noindent The {\em tower}  $\{\M(\Sigma),\; \Sigma\subset \s\; \mbox{admissible}\}$ is the collection   
of mapping class groups of admissible subsurfaces $\Sigma\subset \s$,   
endowed with the following  additional structure: for any     
subsurfaces $\Sigma$ and $\Sigma'$, there is defined the    
subgroup $\M(\Sigma)\cap \M(\Sigma')$ of $\M(\Sigma)$, $\M(\Sigma')$,   
and $\M(\Sigma\cap \Sigma')$.     
         
\begin{definition}          
The tower $\{\M(\Sigma),\; \Sigma\subset \s\; \mbox{admissible}\}$  maps (respectively embeds)   
in the group          
$\Gamma$ if there are given homomorphisms (respectively embeddings)          
$i_{\Sigma}:\M(\Sigma)\to \Gamma$ satisfying the property           
expressed in the diagram of Figure \ref{dia}.           
\end{definition}          
          
\noindent For $x\in \M(\Sigma)$, $\Sigma\subset \s$ admissible,           
let $\widehat{x}:\s\to \s$ be its rigid extension outside $\Sigma$ (see Remark      
\ref{comp}). Let $\widehat{\Sigma}\subset \s$ be any admissible    
surface, not necessarily  invariant by $\widehat{x}$, which contains    
$\Sigma$.          
There exists then (at least) one element $\lambda\in T$ such that           
$\lambda(\widehat{x}(\widehat{\Sigma}))=\widehat{\Sigma}$. In fact, $T$          
acts transitively on the set of admissible surfaces of fixed topological type.             
In particular, $\lambda\circ\widehat{x}\mid_{\widehat{\Sigma}}$           
keeps $\widehat{\Sigma}$ invariant.           
          
\begin{definition}\label{tower}          
The tower  $\{\M(\Sigma),\; \Sigma\subset \s\; \mbox{admissible}\}$            
is left $T$-equivariantly mapped to $\Gamma$ if           
\[i_{\widehat{\Sigma}}(\lambda\circ\widehat{x}\mid_{\widehat{\Sigma}})=          
\rho_{\widehat{x}(\widehat{\Sigma}),\widehat{\Sigma}}(\lambda)          
i_{\Sigma}(x) \]          
for all $x, \Sigma,\widehat{\Sigma}$ and $\lambda$ as above. Here one          
supposes that $\rho_{\Sigma,\Sigma'}(\lambda)\in \Gamma$ (defined for           
$\lambda\in T$ such that $\lambda(\Sigma)=\Sigma'$) is a groupoid          
representation of $T$, i.e.           
\[ \rho_{\Sigma',\Sigma''}(\lambda')\rho_{\Sigma,\Sigma'}(\lambda)=          
\rho_{\Sigma,\Sigma''}(\lambda'\lambda) \]          
whenever this makes sense.           
          
\vspace{0.1cm}   
\noindent The right $T$-equivariance is defined in the same way, but using           
the right action of $T$.           
The tower is $T$-equivariantly mapped if it is           
right and left $T$-equivariantly mapped to $\Gamma$ and           
the left and right groupoid representations of $T$ agree.           
\end{definition}

\begin{proposition}          
$\B$ is universal among $T$-equivariant mappings of the tower          
$\{\M(\Sigma),\; \Sigma\subset \s\;$ $ \mbox{admissible}\}$. More precisely, any           
$T$-equivariant map into $\Gamma$ is induced from a uniquely defined       
homomorphism  $j:\B\to \Gamma$.           
\end{proposition}          
\begin{proof}          
Let $x\in \B$. Assume that $x$ has an invariant support          
$\Sigma\subset \s$. We must set then           
\[ j(x) = i_{\Sigma}(x|_{\Sigma})\in \Gamma\]          
\begin{lemma}          
$j(x)$ does not depend upon the choice of the invariant support           
$\Sigma\subset \s$.          
\end{lemma}          
\begin{proof}          
Let $\Sigma'$ be another invariant support. Then $\Sigma\cap \Sigma'$           
is admissible and also invariant by $x$. Let us show that the          
intersection is nonvoid whenever $x$ is not the identity. Assume the contrary. Then the connected           
component of $\s-\Sigma$ containing $\Sigma'$ must be fixed, since           
$\Sigma'$ is invariant and $x$ can only permute the components    
of $\s-\Sigma$. Then the restriction of $x$ at this component should be the          
identity. In particular, $x|_{\Sigma'}$ is identity. Since $\Sigma'$          
is an invariant  support for $x$ one obtains that $x$ is the identity.          
   
\vspace{0.1cm}   
\noindent    
Now $\Sigma\cap \Sigma'$ is also a support for $x$. In particular      
$x|_{\Sigma\cap \Sigma'}\in \M(\Sigma)\cap \M(\Sigma')$ and thus   
from  definition \ref{tower} we derive that           
$i_{\Sigma}(x|_{\Sigma})=i_{\Sigma\cap\Sigma'}(x|_{\Sigma\cap\Sigma'})=          
i_{\Sigma'}(x|_{\Sigma'})$.           
\end{proof}          
          
\noindent By the results of the next section, $\B$ is generated by four mapping    
classes with invariant supports. Notice however           
that not all elements have invariant supports. So each $x\in\B$ may be    
written as $x=x_1x_2 \cdots x_p\in \B$ where $x_1,x_2,\ldots,x_p$ have invariant    
supports, and one defines           
\[ j(x)=j(x_1) \cdots j(x_p)\]          
One needs to show that this definition is coherent with the previous          
one, and thus yields indeed a group homomorphism.           
For this purpose it suffices  to show that, using the first definition,          
$j(xy)=j(x)j(y)$ holds, when $x$, $y$ and $xy$ have invariant supports      
$\Sigma_x, \Sigma_y$ and $\Sigma$, respectively. Choose $\lambda\in T$ such that $\lambda(y(\Sigma))=\Sigma$.           
Then           
\[ i_{\Sigma}(xy|_{\Sigma})=          
i_{\Sigma}(x\lambda^{-1}|_{\Sigma})i_{\Sigma}(\lambda y|_{\Sigma})=          
i_{\Sigma_x}(x|_{\Sigma_x})\rho_{\Sigma, y(\Sigma)}(\lambda^{-1})          
\rho_{y(\Sigma), \Sigma}(\lambda)i_{\Sigma_y}(y|_{\Sigma_y})=\]          
\[=i_{\Sigma_x}(x|_{\Sigma_x})i_{\Sigma_y}(y|_{\Sigma_y})\]   
where we used the fact that $y(\Sigma)=x^{-1}(\Sigma)$.           
This proves the claim.           
\end{proof}

\section{A presentation for $\B$}           
Let us define now the elements of $\B$ described in the Figures 4           
to 7. Specifically:            
\begin{itemize}           
\item $t$ is a right Dehn twist around a circle C parallel to the boundary      
  component of $S_0$ labeled 3. Given an outward orientation of the surface,   
  this means that $t$ maps an arc crossing C transversely to an arc   
  which turns right as it approaches C. The effect of the twist on the           
  seams is shown on the picture. An invariant support is $S_0$.              
\item $\pi$ is the braiding, acting as a braid in           
  $\M(S_0)$. Assume that $S_0$  is identified with   
  the complex domain $\{|z|\leq 7, |z-3|\geq 1, |z+3|\geq 1\}\subset \C$.    
  A specific homeomorphism in the mapping class of $\pi$    
  is the composition of the counterclockwise   
 rotation of $180$ degrees around the origin --- which exchanges   
the small boundary circles labeled 1 and 2 in the figure ---  
with a map which rotates of $180$ degrees in the clockwise direction   
each boundary circle. The latter can be constructed as follows.     
  
\vspace{0.1cm}  
\noindent   
Let $A$ be an annulus in the plane, which we suppose for simplicity   
to be $A=\{1\leq |z|\leq 2\}$.    
The homeomorphism $D_{A, C}$ acts as the counterclockwise  
rotation of $180$ degrees  
on the boundary circle $C$ and   
keeps the other  boundary   
component pointwise   
fixed:   
\[D_{A, C}(z)=   
\left\{ \begin{array}{ll}  
z \exp(\pi\sqrt{-1}(2-|z|)), & \mbox{ if }  C=\{|z|=1\}\\  
z \exp(\pi\sqrt{-1}(|z|-1)), & \mbox{ otherwise}  
\end{array}\right.\]

\vspace{0.1cm}  
\noindent   
The map we wanted is    
$D^{-1}_{A_0, C_0}D^{-1}_{A_1, C_1}D^{-1}_{A_2, C_2}$, where   
$A_0=\{6\leq |z|\leq 7\}$, $C_0=\{|z|=7\}$,   
$A_1=\{1\leq |z-3|\leq 2\}$, $C_1=\{|z-3|=1\}$,   
$A_2=\{1\leq |z+3|\leq 2\}$, and  $C_2=\{|z+3|=1\}$.    
   
\noindent A  support for $\pi$ is $S_0$.    
One has pictured also the images of the rigid structure.           
   
\item $\beta$ is a rotation of order 3. It is the unique globally rigid   
  mapping class which permutes counterclockwise and cyclically the three   
  boundary circles of $S_0$. An invariant support for $\beta$ is $S_0$.   
   
\item $\alpha$ is a rotation of order 4. Let $\Sigma_{0,4}$ be the 4-holed sphere    
consisting of the union of  $S_0$           
(labeled P in the picture) with the pair of pants above $S_0$ (labeled   
Q). The element $\alpha$ is the unique mapping class which   
  preserves globally the prerigid structure and  permutes counterclockwise and   
  cyclically the four boundary circles of $\Sigma_{0,4}$. An invariant support   
  for $\alpha$ is $\Sigma_{0,4}$.   
   
\end{itemize}

\begin{figure}        
\begin{center}         
\includegraphics{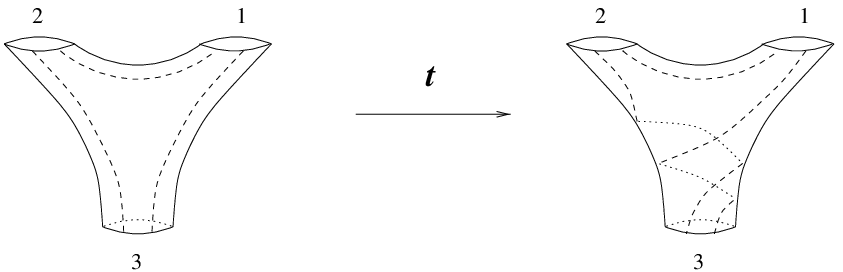}           
\caption{The twist $t$}\label{t}      
\end{center}           
\end{figure}           
           
\begin{figure}        
\begin{center}         
\includegraphics{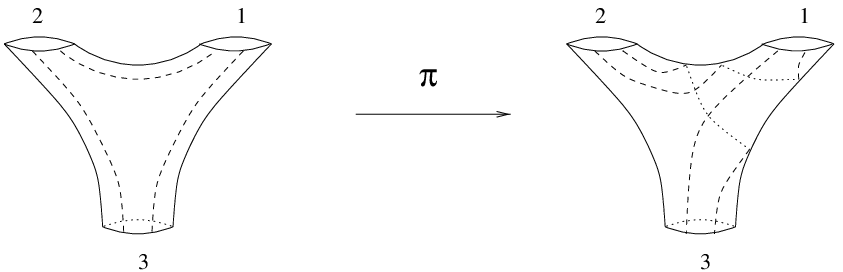}           
\caption{The braiding $\pi$}\label{pi}       
\end{center}          
\end{figure}

\begin{figure}       
\begin{center}          
\includegraphics{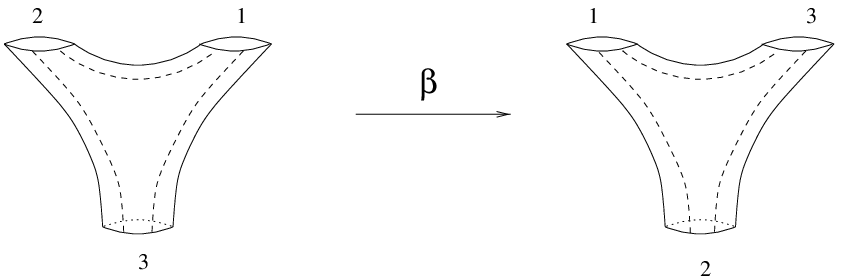}           
\caption{The rotation $\beta$}\label{beta}      
\end{center}           
\end{figure}

\begin{figure}           
\begin{center}      
\includegraphics{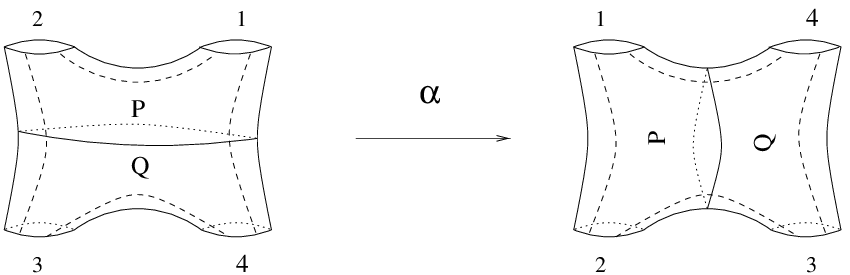}           
\caption{The rotation $\alpha$}\label{F}           
\end{center}      
\end{figure}

\vspace{0.1cm}    
\noindent           
In the following, the notation $g^h$ for two elements $g$ and $h$ of a group stands for $g^{-1}hg$.           
\begin{theorem}\label{pres0}           
The group $\B$ has the following presentation:            
\begin{itemize}           
\item Generators: $t$, $\pi$, $\beta$, and $\alpha$.              
           
\item Relations:            
\begin{enumerate}           
\item Relations at the level of the pair of pants.          
           
\begin{enumerate}           
\item $[t,t_i]=[t,\pi]=1$, $i=1,2$, where $t_1=\beta t           
  \beta^{-1}$, $t_2=\beta^{-1}t\beta$   
\item $t_1 \pi = \pi t_2$, $t_2 \pi = \pi t_1$            
\item $\pi^2=tt_1^{-1}t_2^{-1}$           
\item  $\beta^3=1$          
\item $\beta=t\pi^{\beta}\pi$           
\end{enumerate}           
           
\item Relations coming from the triangle singularities.             
\begin{enumerate}           
\item $(\beta\alpha)^5=1$          
\item $(\alpha\pi)^3=t_2^{-1}$            
\end{enumerate}           
           
\item Relations coming from permutations.            
\begin{enumerate}           
\item $\alpha^4=1$   
\end{enumerate}           
\item Relations coming from commutativity (one sets $t_3=\alpha^2 t_1\alpha^2$   
  and $t_4=\alpha^2 t_2\alpha^2$).           
           
\begin{enumerate}           
\item $[\pi, \alpha^2\pi\alpha^2]=1$            
\item $[t_3,\pi]=[t_4,\pi]=1$     
\item $\alpha t_1 \alpha^{-1}= t_2$          
\item $[t_1,t_3]=1$, $[t_3, \pi^{\beta}t_3\left(\pi^{\beta}\right)^{-1}]=1$, $[t_3, \pi^{\beta}\pi^{\alpha^{2}}\left(\pi^{\beta}\right)^{-1}]=1$\item $[\pi^{\alpha^{2}},           
  \pi^{\beta}\pi^{\alpha^{2}}\left(\pi^{\beta}\right)^{-1}]=1$       
\item $[\beta\alpha\pi^{\beta}, \pi^{\alpha^{2}}]=1$       
\item  $\beta\alpha t_2 (\beta\alpha)^{-1}= \pi^{\beta}t_3 (\pi^{\beta})^{-1}$,      
  $\beta\alpha t_3 (\beta\alpha)^{-1}=\pi^{\beta}t_4 (\pi^{\beta})^{-1}$, $\beta\alpha t_4 (\beta\alpha)^{-1}=t_2$

\end{enumerate}           
\item Consistency relations.            
\begin{enumerate}           
\item $t=\alpha^2t\alpha^2$           
\end{enumerate}           
           
\item Lifts of relations in $T$            
\begin{enumerate}           
\item           
  $(\alpha^2\pi^{\beta}\alpha^3\beta^2)^2=(\pi^{\beta}\alpha^3\beta^2\alpha^2)^{2}$           
\item $\alpha^2\beta\alpha^2\pi           
  \pi^{\beta}\alpha^3\beta^2\alpha^2\beta^2\alpha^2=           
(\alpha^3\beta^2\alpha^2\beta\alpha^2\pi^{\beta})^2$            
\end{enumerate}           
\end{enumerate}           
\end{itemize}    
The terminology used for the classification of the relations is borrowed from \cite{FG}.            
\end{theorem}           
           
\noindent       
As a corollary, one obtains a new presentation of Thompson's group $V$, with 3           
generators and 11 relations (the presentation in \cite{CFP} has 4           
generators and 14 relations).            
\begin{corollary}\label{Vpres}           
Thompson's group $V$ has  the following presentation:            
\begin{itemize}           
\item Generators:  $\pi$, $\beta$, and $\alpha$.              
           
\item Relations:        
           
\begin{enumerate}           
\item $\pi^2=1$           
\item $\beta^3=1$           
\item $\alpha^4=1$           
\item $\pi^{\beta}\pi=\beta$           
\item $(\beta\alpha)^5=1$           
\item $(\alpha\pi)^3=1$            
\item $[\pi, \alpha^2\pi\alpha^2]=1$            
\item $[\beta\alpha\pi^{\beta}, \pi^{\alpha^{2}}]=1$           
\item $[\pi^{\alpha^{2}},           
  \pi^{\beta}\pi^{\alpha^{2}}\pi^{\beta}]=1$           
\item           
  $(\alpha^2\pi^{\beta}\alpha^3\beta^2)^2=(\pi^{\beta}\alpha^3\beta^2\alpha^2)^{2}$           
\item $\alpha^2\beta\alpha^2 \beta^2   \alpha^3\beta^2\alpha^2\beta^2\alpha^2=           
(\alpha^3\beta^2\alpha^2\beta\alpha^2\pi^{\beta})^2$            
\end{enumerate}           
\end{itemize}           
\end{corollary}           
           
\begin{proof}           
The pure mapping class group $K^*_{\infty}$ is the normal subgroup of $\B$           
generated by $t$, and $V\cong \B/K^*_{\infty}$. Notice that relation           
{\it 11} comes from {\it 6(b)}, after replacing $\pi \pi^{\beta}$ by $\beta^2$.           
\end{proof}

\section{$\B$ versus the stable duality groupoid of genus zero $\D$}       
      
\subsection{The stable duality groupoid of genus zero $\D$}      
\begin{figure}           
\includegraphics{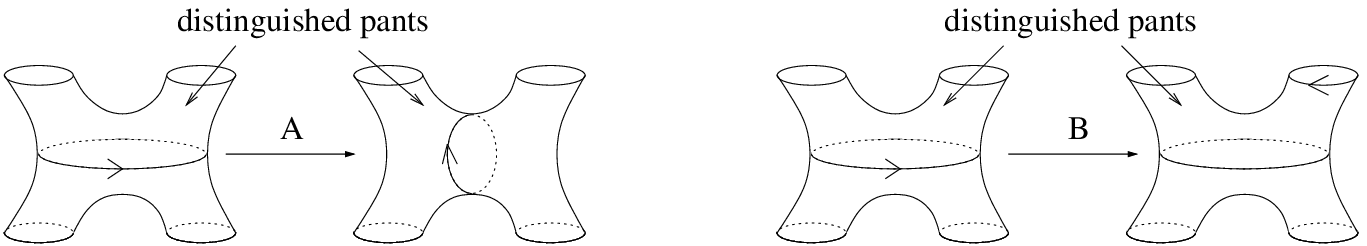}    
\caption{Moves in the Ptolemy groupoid}\label{fi1}            
\end{figure}

\begin{definition}\label{stable}           
The {\em stable duality groupoid of genus zero}, $\D$, is the category   
defined as follows:\\           
           
\noindent $\bullet$ The objects are the isotopy classes of asymptotically trivial rigid           
  structures of $\s$, with  a distinguished          
  oriented circle (among the circles of the rigid structure). The pair of pants bounded by the distinguished circle,          
  which induces on the latter the opposite orientation, is called the {\em          
    distinguished pair of pants} (see Figure \ref{fi1}).\\   
           
\noindent $\bullet$ The  morphisms are words in the moves    
  $T,B,\Pi$ and $A$, defined as follows:\\           
    
-- The moves $T$ and $\Pi$ change a rigid structure as $t$ and           
  $\pi$ do on Figures \ref{t} and \ref{pi}, respectively, where the represented          
  supports must be viewed as the distinguished pairs of pants, the circles        
  labeled 3 being the distinguished circles. Thus $T$ and           
  $\Pi$ leave unchanged the pants decomposition subordinate to the rigid           
  structure and the distinguished circle.\\   
     
-- The moves $B$ and $A$ change a rigid structure as $\beta$ and           
  $\alpha$ do on Figures \ref{beta} and \ref{F}, respectively: on Figure        
  \ref{beta}, the represented supports must be viewed as the distinguished        
  pairs of pants, the circles labeled 3 being the distinguished circles; on Figure \ref{F}, the           
  circles which separate the adjacent pairs of pants are the distinguished           
  circles, and the distinguished pants are labeled by P (see also Figure \ref{fi1}). Thus, the moves $B$           
  and $A$ leave unchanged the prerigid structure subordinate to the rigid   
  structure.\\

\noindent $\bullet$ The relations between the moves are encoded in the assumption that for           
  any rigid structure $r$, the group of morphisms ${\Mor}(r,r)$ must  
  be trivial.  
  
\vspace{0.1cm}  
\noindent            
One denotes by ${\cal R}ig(\s)$ the set of objects of $\D$. If $r$ and $r'$ belong to      
$ {\cal  R}ig(\s)$ and $W\in {\Mor}(r,r')$ (or its inverse) is a move, one      
denotes its target $r'$ by $W\cdot r$.              
\end{definition}           
              
\begin{remark}   
The moves are only changing the rigid structure   
locally. In particular, the moves $B$ and $A$ do not rotate the whole picture. Rather in the first   
case, it just changes which circle is distinguished, and similarly   
for $A$.     
\end{remark}   
   
\begin{definition}           
Let $M$ be the free group on $\{T,\Pi,A,B\}$. The group $M$ acts on ${\cal R}ig(\s)$ via           
$$M\times {\cal R}ig(\s)\rightarrow {\cal R}ig(\s)$$           
$$\left( W=W_1\cdots           
W_n,r \right)\mapsto W_1 \cdots W_{n-1} W_n\cdot r$$    
where $W_i^{\pm 1}\in\{T,\Pi,A,B\}$. Let           
$K$ be the kernel of the resulting morphism $M\rightarrow Aut({\cal R}ig(\s))$. The           
quotient $M/K$ is called the {\em group of moves} of the groupoid $\D$.           
\end{definition}

\begin{proposition}\label{move}           
The group $\B$ and the group of moves $M/K$ are naturally anti-isomorphic. The           
anti-isomorphism is induced by    
$W=W_n\cdots W_1\in M\mapsto w=w_1\cdots w_n\in\B$,           
where $w_i^{\pm 1}=\alpha$ if $W_i^{\pm 1}=A$, $w_i^{\pm 1}=\beta$ if           
$W_i^{\pm 1}=B$, $w_i^{\pm 1}=\pi$ if $W_i^{\pm 1}=\Pi$, and $w_i^{\pm 1}=t$ if           
$W_i^{\pm 1}=T$. For all $W=W_n\cdots W_1\in M/K$, $W(r_*)=w(r_*)$, where           
$r_*$ is the canonical rigid structure.    
\end{proposition}           
           
\begin{proof}           
We first prove the last assertion. By definition of the generators $\alpha,           
\beta,\pi,t$, we have indeed $W_1(r_*)=w_1(r_*)$. But $W_2$ acts on $W_1(r_*)$           
as the conjugate of $w_2$ by $w_1$:           
$W_2(W_1(r_*))=(w_1w_2w_1^{-1})\cdot w_1(r_*)=w_1(w_2(r_*))$. Inductively we check           
that $W(r_*)=w(r_*)$.   
   
\vspace{0.1cm}           
\noindent Since $M$ is free, the anti-morphism $W\in M\mapsto w \in \B$ is           
well-defined. It is surjective, since the image contains the generators of the           
group $\B$. Let $W=W_1\cdots W_n \in M$ be in the kernel of the           
anti-isomorphism i.e. $w=w_n\cdots w_1=1$ in $\B$. We prove that $W$           
belongs to $K$.    
   
\vspace{0.1cm}           
\noindent   
For any rigid structure  $r$, there exists some $w'\in           
\B$ such that $r=w'(r_*)$. Choose $W'\in M$ in the preimage of $w'$. Then           
$W(r)=WW'(r_*)=w'w(r_*)=w'(r_*)=r$, so $W\in K$. Clearly, $K$ lies in the           
kernel of the anti-morphism, since if $w\in \B$ fixes pointwise    
a rigid structure, then it must be trivial. Hence $M/K$ is anti-isomorphic to $\B$.           
\end{proof}             
           
\begin{remark}\label{K}           
We note from the proof above that $K$ coincides with the stabilizer of any           
rigid structure. Thus, for all $r,r'$ in ${\cal R}ig(\s)$, there exists a unique      
$W\in M/K$ such that $r'=W\cdot r.$            
\end{remark}       
              
\subsection{$\D$ is a stabilization of ${\mathcal D}_0$, the duality groupoid   
  in  genus zero}   
   
The purpose of this subsection is to relate $\D$ with the the duality   
groupoid evoked in the Introduction. Its content will not be used in the sequel.         
           
\begin{definition}\label{dual0}           
The duality groupoid ${\mathcal D}$ considered in \cite{FG,MS}   
consists of the transformations of           
rigid structures. The duality groupoid of genus zero, ${\mathcal D}_0$, is the   
subgroupoid of ${\mathcal D}$ consisting of the transformations involving only   
genus zero surfaces:\\   
   
\noindent $\bullet$ Its objects are pairs $(\Sigma,r_{\Sigma})$,      
where $\Sigma$ is a compact surface of genus zero, and $r_{\Sigma}$ is a rigid      
structure on $\Sigma$, which, in the sense of \cite{FG}, consists of a   
decomposition of $\Sigma$ into pairs of pants, a      
collection of seams on the pairs of pants, a numbering of the pairs of pants,      
and for each pair of pants, a labeling by 1,2,3 of the boundary circles. They are defined up to isotopy.\\   
   
\noindent $\bullet$ Its morphisms are   
changes of rigid structure.

\vspace{0.1cm}  
\noindent                                           
There is a tensor structure $\otimes$ on ${\mathcal D}_0$, which corresponds   
to the connected sum along boundary components.   
\end{definition}

\noindent    
A finite presentation  for ${\mathcal D}$ has been obtained   
in \cite{FG}. It is easy to check that the subset of generators            
$T_1,B_{23},R,F$ and $P$ together with the relations provided in    
\cite{FG}, which are involving only these generators,    
form a presentation for ${\mathcal D}_0$.            
   
\begin{proposition}\label{dual}           
The duality groupoid of genus zero, ${\mathcal D}_0$, has the following presentation:            
\begin{enumerate}           
\item[] Generators: $T_1, R, B_{23}, F, P$ and their inverses.            
\item[]  Relations (Moore-Seiberg equations): \\           
1. at the level of a pair of pants:\\           
\ a) $T_1B_{23}=B_{23}T_1$,\ $T_2B_{23}=B_{23}T_3$, \ $T_3B_{23}=B_{23}T_2$,           
where $T_2=R^{-1}T_1R$ and $T_3=RT_1R^{-1}$,\\           
\ b) $B_{23}^2=T_1T_2^{-1}T_3^{-1}$\\           
\ c) $R^3=1$\\      
\ d) $R=B_{23}RB_{23}R^{-1}T_1$\\           
2. relations defining  inverses:\\           
\ a) $P^{(12)}F^2=1$\\           
\ b) $T_3^{-1}B_{23}^{-1}S^2=1$\\           
3. relations coming from ``triangle singularities'':\\           
\ a) $P^{(13)}{R^{(2)}}^2F^{(12)}{R^{(2)}}^2F^{(23)}{R^{(2)}}^2           
F^{(12)}{R^{(2)}}^2F^{(23)}{R^{(2)}}^2F^{(12)}=1$\\           
\ b) $T_3^{(1)}FB_{23}^{(1)}FB_{23}^{(1)}FB_{23}^{(1)}=1$\\           
4. relations coming from the symmetric groups:\\           
\ a) $P^2=1$ \\           
\ b) $P^{(23)}P^{(12)}P^{(23)}=P^{(12)}P^{(23)}P^{(12)}$\\           
\end{enumerate}           
\end{proposition}   
           
\begin{remark} We used the convention that superscripts tell us           
on which factors of the tensor product the move acts. Here the  tensor           
structure is implicit. The pentagon relation 3.a) is corrected here           
since it was erroneously stated in \cite{FG} p.608 (see p.635) while      
the relation 1.d) was omitted.            
\end{remark}

\vspace{0.1cm}  
\noindent   
 The following proposition clarifies in which sense the stable groupoid in   
 genus zero, $\D$, is indeed a stabilization of the duality      
groupoid ${\mathcal D}_0$.  
  
\vspace{0.1cm}  
\noindent  
  View the group of moves $M/K$ of the groupoid $\D$ as the tautological      
groupoid (with a unique object, and $M/K$ as the set of morphisms).   
       
\begin{proposition}\label{quotient}      
There exists a      
surjective morphism of groupoids  $s: {\mathcal D}_0\rightarrow M/K$.      
\end{proposition}      
      
\begin{proof}      
First, $s$ maps the objects of ${\mathcal D}_0$ on the unique object of $M/K$.      
Let $(\Sigma,r_{\Sigma})$ be an object of ${\mathcal D}_0$. The source and the      
target of any morphism of ${\mathcal D}_0$ may be represented by two      
rigid structures $r_{\Sigma}$ and $r_{\Sigma}'$ on the same support      
$\Sigma$. So, let $\varphi:(\Sigma,r_{\Sigma})\rightarrow      
(\Sigma,r_{\Sigma}')$ be such a morphism. One may identify $\Sigma$ with an      
admissible subsurface of $\s$ (in a non-unique way). Let $r$ (resp. $r'$) be      
the object of $\D$ with the following properties:   
\begin{itemize}      
\item  it coincides ith the canonical rigid structure of $\s$   
  outside $\Sigma$;     
\item  it is induced by $r_{\Sigma}$ (resp. $r_{\Sigma}'$) on $\Sigma$;      
\item  its oriented circle is the labeled 1 circle of the labeled 1 pair of      
pants of $r_{\Sigma}$ (resp. $r_{\Sigma}'$), viewed as the distinguished pair  
of pants of $r$ (resp. $r'$), and oriented as explained in      
Definition \ref{stable}.   
\end{itemize}      
There is a unique $W\in M/K$ such that $r'=W\cdot r$ (see Remark      
\ref{K}). Moreover, $W$ does not depend on the choice of the representative      
$(\Sigma\subset \s,r_{\Sigma})$, but only on $\varphi$. One sets $s(\varphi)=W$ and this      
completes the definition of $s$. It is easy to check that $s(\varphi \psi)=      
s(\varphi)s(\psi)$.   
   
\vspace{0.1cm}      
\noindent If $\varphi$ is a morphism of type $T_1^{()},R^{()}, B_{23}^{()}$ or $F^{()}$ (where      
the superscript tells us on which factors of the tensor product the move acts),      
then $s(\varphi)$ is conjugate (by some element $W\in M/K$ depending on the superscript) to a morphism of type $T, B, \Pi$ or $A^{-1}$      
respectively. This easily      
implies the surjectivity of $s$. Notice that the image by $s$ of a move $P$      
which transposes the numberings of two pairs of pants is a word in $A^2$ and $B$.\end{proof}

\subsection{The stable duality groupoid $\D$ generalizes the universal Ptole\-my groupoid}      
      
The universal Ptolemy groupoid ${\cal P}t$ appears in   
\cite{pe0,pe,im,LS}. We translate its definition to a language related to our          
framework:           
           
\begin{definition}           
The universal Ptolemy groupoid ${\cal P}t$ is the full subgroupoid of $\D$           
whose objects are the isotopy classes of asymptotically trivial pants decompositions    
(with distinguished  oriented circles) which are compatible    
with the canonical prerigid structure. Hence           
the morphisms are the composites of the $A$ and $B$ moves (which suffice,   
since $A^{-1}=A^3$ and $B^{-1}=B^2$).   
\end{definition}           
             
\noindent From Proposition \ref{move}, we may deduce that the group of moves of ${\cal P}t$           
is anti-isomorphic to the subgroup of $\B$ generated by $\alpha$ and $\beta$:           
this is precisely the Ptolemy-Thompson group $T\subset V$, cf. Proposition           
\ref{sequence}. From the presentation of the group $T$, a presentation of the           
groupoid ${\cal P}t$ has been obtained in \cite{LS}:  ${\cal P}t$ is            
presented by the generators $A$, $B$, and  relations:            
$$A^4=1,\;\; B^3=1,\;\; (AB)^5=1$$      
$$ [BAB,A^2BAB A^2]=1,\;\; [BAB, A^2 B^2A^2            
BABA^2BA^2]=1$$            
This means that if a sequence of $A$ and $B$ moves starts and finishes at            
the same object, then the sequence is a product of conjugates of the sequences            
of the presentation. Equivalently, Ptolemy-Thompson's group $T$ is presented          
by the generators $\alpha$, $\beta$, and the relations:          
$${\alpha}^4=1, {\beta}^3=1, (\beta\alpha)^5=1,          
[\beta\alpha\beta,\alpha^2\beta\alpha\beta\alpha^2]=1, [\beta\alpha\beta,\alpha^2\beta^2\alpha^2\beta\alpha\beta\alpha^2\beta\alpha^2]=1$$                
            
\begin{definition}\label{PT}           
Let $\Gr({\cal P}t)$ be the graph whose vertices are the objects of   
${\cal P}t$, with edges corresponding to $A$ or $B$ moves.    
\end{definition}          
          
\noindent Thus, $\Gr({\cal P}t)$ is          
  a subcomplex of the classifying space of the category ${\cal P}t$, and          
  is exactly the Cayley graph of Ptolemy-Thompson's group $T$, generated by          
  $\alpha$ and $\beta$.  This assertion relies on the fact that in the   
  Cayley graph of a group (with a chosen set of generators),    
 there is an edge between two elements $g_1$ and $g_2$ if and only if    
 $g_1^{-1}g_2$       is a generator of the group.

\begin{remark}\label{PSL2}    
The subgroup of $T$ generated by $\alpha^2$ and $\beta$ is isomorphic to    
$PSL(2,\Z)$, viewed as the group of orientation preserving automorphisms of    
the planar tree ${\cal T}$ (see also Remark \ref{PSL}). Using the duality    
between ${\cal T}$ and the canonical pants decomposition, one deduces   
that, given    
two objects $r_1$ and $r_2$ of ${\cal P}t$ differing only by the   
position and the orientation of the    
distinguished circles, there exists a unique sequence of $B$ moves and squares    
of $A$ moves connecting $r_1$ to $r_2$.    
\end{remark}

\section{The group $\B$ is finitely presented}           
The method we develop here follows the approach of Hatcher-Thurston in their      
proof that the mapping class groups of           
compact surfaces (\cite{HT}) are finitely presented. However, the           
Hatcher-Thurston complex of the (non-compact) surface $\s$, even           
  restricted to the asymptotically canonical pants decompositions, is too large for providing a finite presentation of ${\cal B}$. Instead, the ${\cal B}$-complex we construct (\S5.3, Definition \ref{k3}) will           
  use Hatcher-Thurston complexes of compact holed-spheres with bounded           
  levels only, together with simplicial complexes closely related to those   
  used by K. Brown in the study of finiteness properties of the Thompson group           
  $V$ (Brown-Stein complexes of bases, cf. \cite{br2}). By adding some cells to link both types of complexes and to kill some combinatorial loops, we shall obtain a           
  simply connected ${\cal B}$-complex, with a finite number of cells modulo           
  ${\cal B}$. By a standard theorem (\cite{br}), it will follow that           
  ${\cal B}$ is finitely presented.  
  
\vspace{0.1cm}  
\noindent  
As several complexes will appear throughout this section, we give below a road   
map for all of them. An arrow $``X\longrightarrow Y"$ means that the complex $X$ is   
introduced and studied for proving that the complex $Y$ is connected and   
simply connected.

\begin{center}   
\begin{picture}(0,0)%
\includegraphics{road.pstex}%
\end{picture}%
\setlength{\unitlength}{3315sp}%
\begingroup\makeatletter\ifx\SetFigFont\undefined%
\gdef\SetFigFont#1#2#3#4#5{%
  \reset@font\fontsize{#1}{#2pt}%
  \fontfamily{#3}\fontseries{#4}\fontshape{#5}%
  \selectfont}%
\fi\endgroup%
\begin{picture}(6177,3268)(586,-2774)   
\put(2161,-421){\makebox(0,0)[lb]{\smash{\SetFigFont{10}{12.0}{\rmdefault}{\mddefault}{\updefault}${\cal HT}_{red}(\s)$}}}   
\put(4231,-1051){\makebox(0,0)[lb]{\smash{\SetFigFont{10}{12.0}{\rmdefault}{\mddefault}{\updefault}$\widetilde{\cal DP}(\s)$}}}   
\put(6391,-1051){\makebox(0,0)[lb]{\smash{\SetFigFont{10}{12.0}{\rmdefault}{\mddefault}{\updefault}${\cal DP}(\s)$}}}   
\put(721,299){\makebox(0,0)[lb]{\smash{\SetFigFont{10}{12.0}{\rmdefault}{\mddefault}{\updefault}$\Gr({\cal P}t)$}}}   
\put(6391,-2491){\makebox(0,0)[lb]{\smash{\SetFigFont{10}{12.0}{\rmdefault}{\mddefault}{\updefault}${\cal DP}_5^+(\s)$}}}   
\put(2431,-1861){\makebox(0,0)[lb]{\smash{\SetFigFont{10}{12.0}{\rmdefault}{\mddefault}{\updefault}${\cal BS}_{3,7}(p)$}}}   
\put(991,-2716){\makebox(0,0)[lb]{\smash{\SetFigFont{10}{12.0}{\rmdefault}{\mddefault}{\updefault}${\cal BS}(p)$}}}   
\put(586,-1186){\makebox(0,0)[lb]{\smash{\SetFigFont{10}{12.0}{\rmdefault}{\mddefault}{\updefault}${\cal HT}(\s)$}}}   
\end{picture}

\end{center}          
    
\noindent {\it Nota bene:} From now on, the topological objects associated to $\s$,   
namely, the circles, the pants decompositions and the (pre)rigid structures,   
will be considered up to isotopy. This way, the group $\B$ acts on them.

\subsection{Hatcher-Thurston complexes of the infinite surface}           
           
\begin{definition} Let ${\cal HT}(\s)$ denote the Hatcher-Thurston           
  complex of $\s$:          
\begin{enumerate}          
\item The vertices are the asymptotically trivial            
  pants decompositions of $\s$.           
\item The edges correspond to pairs of pants decompositions $(p,p')$ which    
differ by a local A move, i.e. $p'$ is obtained from $p$ by replacing   
one curve $\gamma$ in $p$ by a curve $\gamma'$ that intersects   
$\gamma$ twice.    
   
\item The 2-cells are introduced to fill in the cycles of moves of            
  the following types: triangular cycles, square cycles (of disjointly   
  supported A moves) and pentagonal cycles (cf. \cite{HT, HLS}).          
\end{enumerate}            
\end{definition}

\begin{remark}         
A pants decomposition is codified by a Morse function on the           
surface. Then the A move  is the elementary non-trivial change induced by           
a small isotopy among smooth functions which crosses once transversally the           
discriminant locus made of functions in Cerf's stratification,  and it           
is therefore a local change in this respect, too.         
\end{remark}        
\begin{remark}          
We should mention that an $A$ move involves unoriented circles, unlike           
the $A$ move in the Ptolemy groupoid.         
\end{remark}   
   
\begin{figure}    
\begin{center}          
\includegraphics{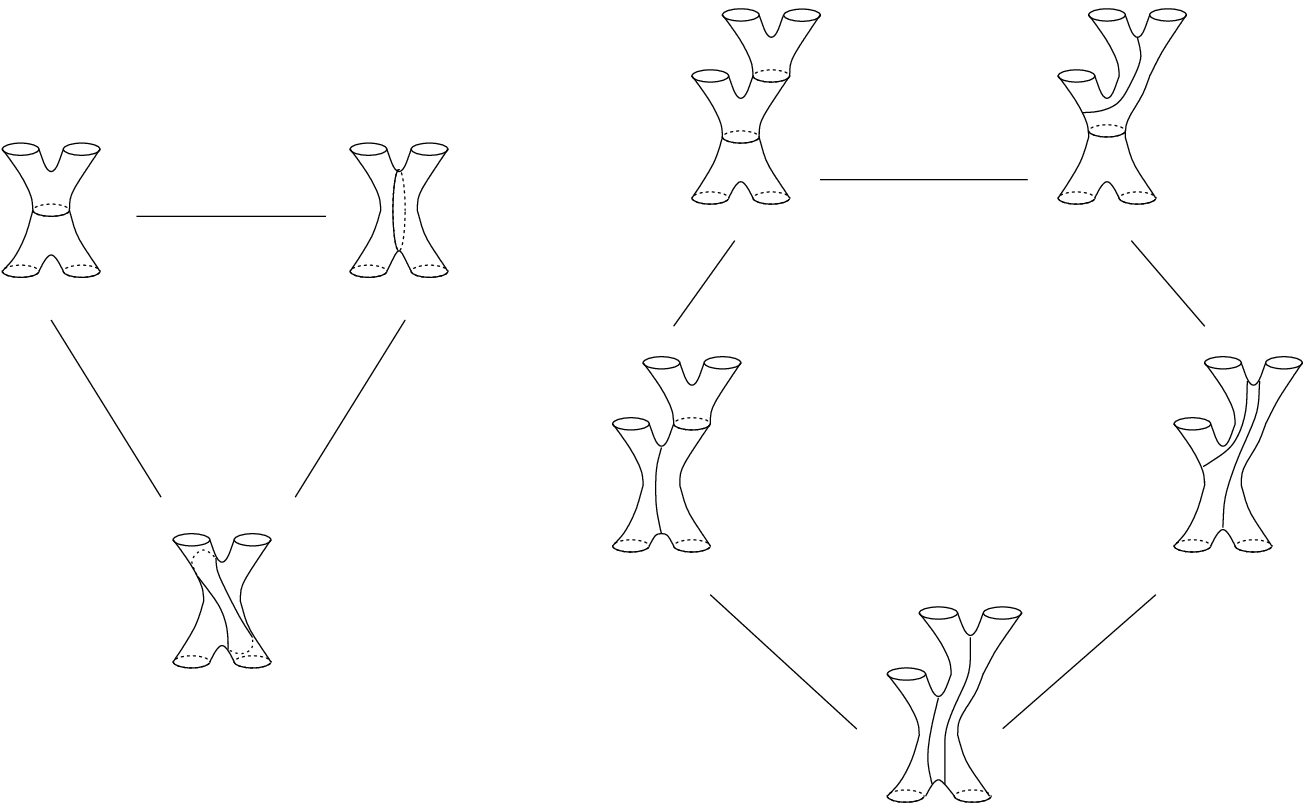}    
\caption{Triangular cycle and pentagonal cycle in ${\cal HT}(\s)$}\label{ficell}\end{center}            
\end{figure}

\begin{remark}\label{detail}   
A triangular (resp. pentagonal) 2-cell $\sigma$ is determined by a connected subsurface   
$\Sigma\subset \s$ of level 4 (resp. 5), an asymptotically trivial pants decomposition $p_{\sigma}$ on   
$\s\setminus \Sigma$, and a cycle of three (resp. five) A   
moves inside $\Sigma$ (cf. Figure \ref{ficell}).  
  
\vspace{0.1cm}   
\noindent  A square 2-cell, generically denoted $(DC)$, is   
determined by two disjoint or adjacent connected subsurfaces $\Sigma$ and $\Sigma'$ of   
level 4, an asymptotically trivial pants decomposition $p_{(DC)}$ on $\s\setminus (\Sigma\cup \Sigma')$, and two A moves,   
supported in $\Sigma$ and $\Sigma'$, respectively. One defines an integer   
associated to $(DC)$, called the {\it distance} $d_{(DC)}$: this is the minimal integer $r\geq 0$ such that there exists   
pairs of pants $P_1,\ldots,P_r$ belonging to $p_{(DC)}$, with $P_1$ adjacent to $\Sigma$, $P_r$ adjacent to $\Sigma'$, and $P_i$ adjacent to   
$P_{i+1}$, for all $i=1,\ldots, r-1$.   
\end{remark}         
   
\begin{proposition}\label{HT}             
The complex ${\cal HT}(\s)$ is connected and simply connected. The group $\B$   
acts cellularly on it, with one orbit of 0-cells, one orbit of 1-cells, one   
orbit of triangular 2-cells, one orbit of pentagonal 2-cells, but  countably   
many orbits of square 2-cells. Two square 2-cells $(DC)$ and $(DC')$ are in  
the same orbit if, and only if, $d_{(DC)}=d_{(DC')}$.             
\end{proposition}   
   
\begin{proof} The first assertion results from \cite{HT,HLS}, by describing the complex as an            
inductive limit of the Hatcher-Thurston complexes ${\cal HT}(S)$ of the            
admissible  supports $S\subset{\s} $, with respect to the various   
inclusions $S\subset S'$. Equivalently, if one wants to connect $p$ to   
$p'$, it is enough to do so in the pants complex   
for some surface ``containing"   
both, and similarly for the simple connectedness.  
  
\vspace{0.1cm}   
\noindent   
That $\B$ acts transitively on the set of 0-cells and on the set of 1-cells is   
obvious. Consider next two triangular 2-cells $\sigma_1$ and $\sigma_2$. Let   
$\Sigma_i$ ($i=1,2$) be the subsurface of level 4 of $\s$ which supports the   
three A moves involved in $\sigma_i$ ($i=1,2$) and $p_{\sigma_i}$ ($i=1,2$) be the   
pants decomposition outside $\Sigma_i$ (cf. Remark \ref{detail}). Denote by $\overline{\sigma_i}$ the   
triangular 2-cell corresponding to $\sigma_i$, viewed in ${\cal H   
  T}(\Sigma_i)$. Since ${\cal M}(0,4)$ acts   
transitively on the 2-cells of ${\cal H T}(\Sigma_{0,4})$, one deduces that there exists a   
homeomorphism from $\Sigma_1$ to $\Sigma_2$ mapping $\overline{\sigma_1}$ onto   
$\overline{\sigma_2}$. The pants decomposition $p_{\sigma_1}$ and $p_{\sigma_1}$   
being asymptotically trivial, one can extend $h$ to an asymptotically rigid   
homeomorphism of $\s$ which maps $p_{\sigma_1}$ onto $p_{\sigma_2}$. The   
resulting mapping class belongs to $\B$, and maps the cell $\sigma_1$ onto   
$\sigma_2$.  
  
\vspace{0.1cm}   
\noindent    
Using the transitivity of the action of ${\cal M}(0,5)$ on the pentagonal   
2-cells of $\Sigma_{0,5}$, one proves similarly that there exists a   
unique orbit of pentagonal 2-cells in ${\cal HT}(\s)$.  
  
\vspace{0.1cm}   
\noindent    
Finally, one proves by similar arguments that two square 2-cells $(DC)$ and   
$(DC')$ are equivalent modulo $\B$ if and only if $d_{(DC)}=d_{(DC')}$.   
\end{proof}

\begin{figure}  
\begin{center}            
\begin{picture}(0,0)%
\includegraphics{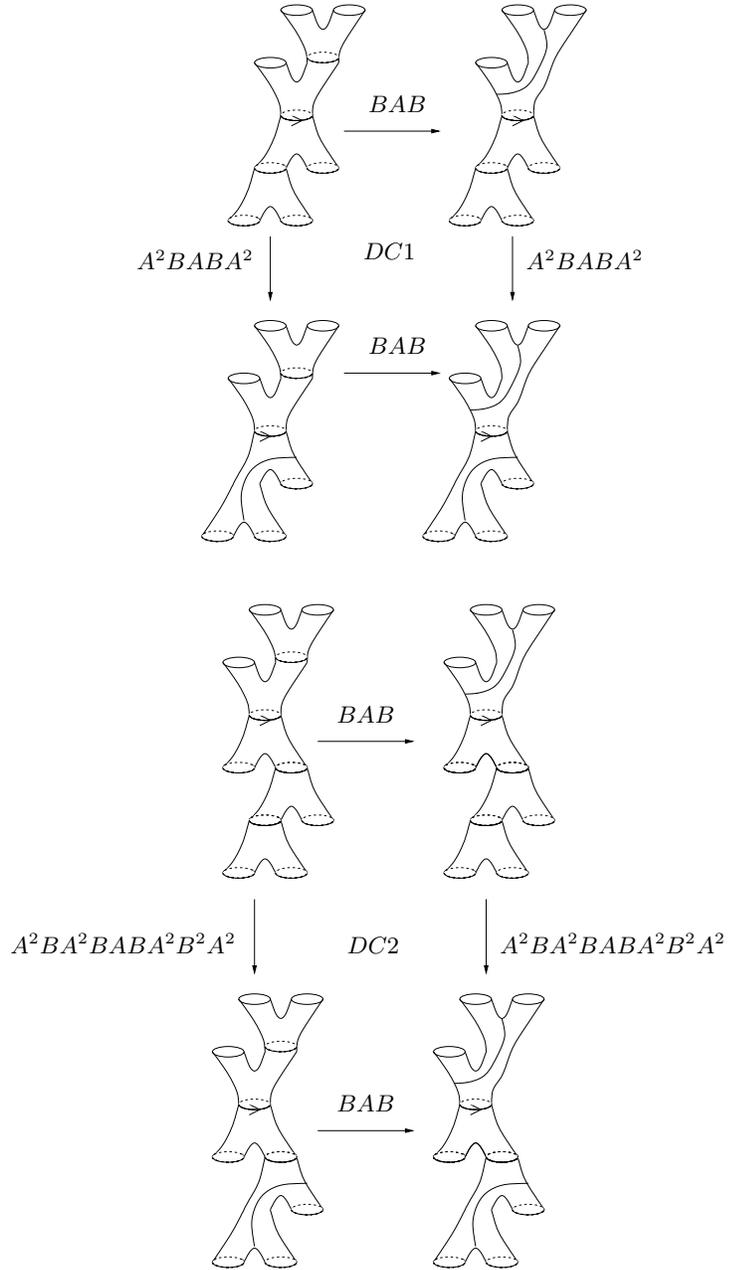}%
\end{picture}%
\setlength{\unitlength}{2901sp}%
\begingroup\makeatletter\ifx\SetFigFont\undefined%
\gdef\SetFigFont#1#2#3#4#5{%
  \reset@font\fontsize{#1}{#2pt}%
  \fontfamily{#3}\fontseries{#4}\fontshape{#5}%
  \selectfont}%
\fi\endgroup%
\begin{picture}(4917,10816)(811,-10779)  
\put(3826,-2131){\makebox(0,0)[lb]{\smash{\SetFigFont{9}{10.8}{\rmdefault}{\mddefault}{\updefault}$DC1$}}}  
\put(5221,-2221){\makebox(0,0)[lb]{\smash{\SetFigFont{9}{10.8}{\rmdefault}{\mddefault}{\updefault}$A^2BABA^2$}}}  
\put(3691,-8071){\makebox(0,0)[lb]{\smash{\SetFigFont{9}{10.8}{\rmdefault}{\mddefault}{\updefault}$DC2$}}}  
\put(1891,-2221){\makebox(0,0)[lb]{\smash{\SetFigFont{9}{10.8}{\rmdefault}{\mddefault}{\updefault}$A^2BABA^2$}}}  
\put(811,-8071){\makebox(0,0)[lb]{\smash{\SetFigFont{9}{10.8}{\rmdefault}{\mddefault}{\updefault}$A^2BA^2BABA^2B^2A^2$}}}  
\put(3601,-6091){\makebox(0,0)[lb]{\smash{\SetFigFont{9}{10.8}{\rmdefault}{\mddefault}{\updefault}$BAB$}}}  
\put(3601,-9421){\makebox(0,0)[lb]{\smash{\SetFigFont{9}{10.8}{\rmdefault}{\mddefault}{\updefault}$BAB$}}}  
\put(4996,-8071){\makebox(0,0)[lb]{\smash{\SetFigFont{9}{10.8}{\rmdefault}{\mddefault}{\updefault}$A^2BA^2BABA^2B^2A^2$}}}  
\put(3871,-2941){\makebox(0,0)[lb]{\smash{\SetFigFont{9}{10.8}{\rmdefault}{\mddefault}{\updefault}$BAB$}}}  
\put(3871,-871){\makebox(0,0)[lb]{\smash{\SetFigFont{9}{10.8}{\rmdefault}{\mddefault}{\updefault}$BAB$}}}  
\end{picture}  
            
\caption{Cells (DC1) and (DC2)}\label{fi2}   
\end{center}           
\end{figure}

\begin{definition}{\bf (Reduced Hatcher-Thurston complex)}\\            
Let ${\cal HT}_{red}(\s)$ be the subcomplex of ${\cal  HT}(\s)$ whose   
vertices, edges, triangular, and pentagonal 2-cells are those of   
${\cal  HT}(\s)$, but whose square 2-cells are of two types:            
\begin{enumerate}            
\item Square $(DC1)$, associated with two commuting moves on two 4-holed            
  spheres adjacent along a common boundary circle (hence $d_{(DC1)}=0$), cf. Figure \ref{fi2};            
\item Square $(DC2)$, associated with two commuting moves on two 4-holed            
  spheres separated by a pants surface (hence $d_{(DC2)}=0$), cf. Figure \ref{fi2}.            
\end{enumerate}            
\end{definition}

\noindent The next proposition plays a key role in the proof of Theorem \ref{main}. It relies on the existence of a cellular map $\nu:\Gr({\cal P}t)\rightarrow {\cal            
  HT}(\s)$.

\begin{proposition}            
The reduced Hatcher-Thurston complex ${\cal HT}_{red}(\s)$ is            
connected and simply connected. The group $\B$ acts cellularly on it, with one orbit of 0-cells, one orbit of 1-cells, one   
orbit of triangular 2-cells, one orbit of pentagonal 2-cells, and two orbits   
of square 2-cells.    
\end{proposition}           
     
\begin{proof} It suffices to prove that each square cycle $(DC)$ corresponding to two            
commuting moves which are supported on arbitrarily far disjoint 4-holed spheres is    
a            
product of conjugates of square cycles of types $(DC1)$ and $(DC2)$ and of            
pentagonal cycles. We use the            
Ptolemy groupoid. Let $\Gr({\cal P}t)$ be the graph introduced in          
Definition \ref{PT}. There is an obvious            
forgetful cellular map    
$$\nu:\Gr({\cal P}t)\rightarrow {\cal   HT}(\s)$$      
defined on the set of vertices by forgetting the orientation of the   
distinguished circle and the fact that it is distinguished. The cycles of $\Gr({\cal P}t)$ corresponding to the first            
three relations of ${\cal P}t$ project by $\nu$ onto cycles of ${\cal   
  HT}(\s)$. The projection by $\nu$ of the cycle $C_{A^4}$ corresponding to $A^4$ is             
a degenerate cycle of the form $ee^{-1}ee^{-1}$, where $e$ is an    
oriented edge and $e^{-1}$ the edge with the opposite orientation;    
$\nu(C_{B^3})$ is a            
0-cell and  $\nu(C_{(AB)^5})$ is a pentagonal cycle. It is elementary    
to check that  $\nu(C_{[BAB,A^2BAB A^2]})$ is            
a square cycle of type $(DC1)$, and $\nu(C_{[BAB, A^2 BA^2BABA^2B^2A^2]})$   
a square cycle of type $(DC2)$.   
   
\vspace{0.1cm}    
\noindent          
Now let $(DC)$ be an            
arbitrary square cycle of ${\cal HT}(\s)$. We can find an   
asymptotically    
trivial prerigid structure $r$      
such that the four vertices of $(DC)$, that is, the four pants decompositions,      
are all compatible with $r$. At the price of replacing $(DC)$ by an equivalent   
(hence homeomorphic) cycle modulo $\B$, we may suppose that      
$r$ is the canonical prerigid structure. It follows that we can find a ``lift"            
$\widetilde{(DC)}$ of $(DC)$ in $\Gr({\cal P}t)$ as follows. We first lift the    
4 edges $e_1,\ldots,e_4$ of $(DC)$ to 4 edges    
$\widetilde{e_1},\ldots,\widetilde{e_4}$, corresponding to $A$ moves in ${\cal    
  P}t$. Since the terminal vertex $t(\widetilde{e_i})$ of $\widetilde{e_i}$    
and the origin $o(\widetilde{e_{i+1}})$ of    
$\widetilde{e_{i+1}}$ (for $i=1,\ldots,4$ mod 4) are pants decompositions    
which differ by the positions of the distinguished circle, we must (and can)    
find an edge-path $\gamma_i$ in $\Gr({\cal P}t)$ joining $t(\widetilde{e_i})$ to    
$o(\widetilde{e_{i+1}})$. By Remark \ref{PSL2}, $\gamma_i$ is a composite of    
$B$ moves and squares of $A$ moves. Let now $\widetilde{(DC)}$ be the cycle $\widetilde{e_1}\gamma_1    
\widetilde{e_2}\gamma_2\widetilde{e_3}\gamma_3\widetilde{e_4}\gamma_4$. Since    
$\nu(B)$ is a point and $\nu(A^2)$ is homotopic to a point, $\nu(\widetilde{(DC)})$ is    
homotopically equivalent to $(DC)$. But the cycle $\widetilde{(DC)}$ can be expressed as a product of            
conjugates of cycles of the presentation of ${\cal P}t$. We then project this product of cycles             
onto ${\cal HT}(\s)$, and obtain $\nu(\widetilde{(DC)})$, expressed            
as a product of conjugates of cycles of the subcomplex ${\cal HT}_{red}(\s)$    
(namely, pentagons,    
squares $(DC1)$ and $(DC2)$, and degenerate cycles of the form  
$\nu(A^2)$). Therefore, the cycle $(DC)$ is homotopically trivial in ${\cal HT}_{red}(\s)$.\\   
   
\noindent The last assertion of the proposition is a direct consequence of Proposition \ref{HT}.             
\end{proof}           
           
\subsection{An auxiliary pants decomposition complex}              
Our main object will be a ${\cal      
  B}$-complex ${\cal        DP}(\s)$, which is connected, simply      
  connected and finite modulo  ${\cal B}$. In order to prove its simple      
  connectivity, we introduce an auxiliary complex $\widetilde{\cal DP}(\s)$            
  containing ${\cal DP}(\s)$, and prove that  it is simply connected. By      
  studying the inclusion ${\cal  DP}(\s)\subset \widetilde{\cal      
  DP}(\s)$, we finally prove the same property for ${\cal DP}(\s)$.      
              
\begin{definition}[Complex $\widetilde{\cal DP}(\s)$]\label{k2}             
The complex $\widetilde{\cal DP}(\s)$ is a two-di\-mensional cellular            
complex whose vertices are triples $(p,S,r)$, where:   
\begin{itemize}    
\item $p$ is an asymptotically trivial pants decomposition of $\s$,         
\item $S$ is a surface of level $l\in \{3,\ldots,7\}$, {\it compatible with}   
  $p$ (or {\it p-compatible}), that is, a      
connected compact subsurface of $\s$ of level $l$ bounded by circles of $p$, and endowed      
with the pants decomposition induced by $p$ (but devoid of its seams), and   
\item  $r$ is a rigid structure on            
$\s\setminus S$, to which $p$ is subordinate      
outside $S$.   
\end{itemize}      
One says that the            
level of $(p,S,r)$ is $l$ (see Figure \ref{fi3}), and that $S$ is    
its {\it support}.                
\begin{figure}           
\begin{center}      
\includegraphics{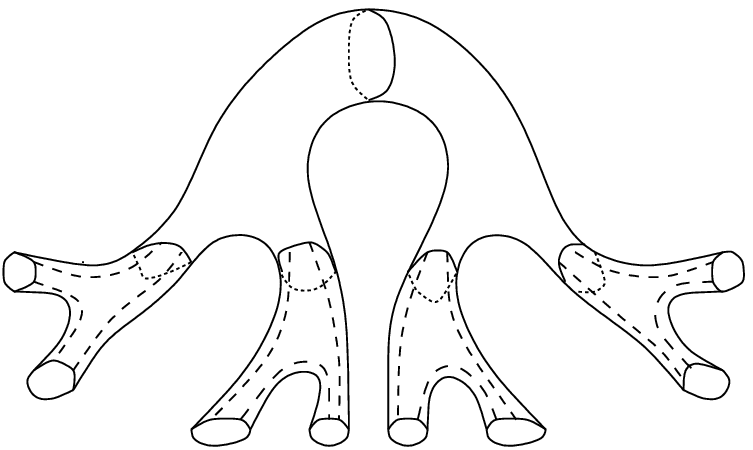}      
\caption{A vertex $(p,S,r)$ of level 4 of $\widetilde{\cal DP}(\s)$}\label{fi3}          
\end{center}            
\end{figure}       
   
\vspace{0.1cm}            
\noindent There are 1-cells associated with moves of three kinds:            
\begin{enumerate}            
\item $A$ moves: they are defined  as in the Hatcher-Thurston    
 complexes, but we restrict to those  which preserve the            
  support $S$ and the rigid structure $r$ of a vertex $(p,S,r)$,            
  and change the pants decomposition $p$ inside $S$.            
\item Propagation moves $P$: a $P$ move on $(p,S,r)$ consists of choosing            
  (finitely many) pairs of pants of $p$ inside $\s\setminus S$,            
  {\it all} adjacent to $S$,            
  and erasing their seams. This changes $(p,S,r)$ into $(p,S',r')$, where            
  $S'\supset S$, and $r'$ coincides with $r$ on $\s\setminus S'$. One also says that $(p,S,r)$ results from $(p,S',r')$ by a            
  $P^{-1}$ move.            
\item Braiding moves $Br$:  a braiding move on   
$(p,S,r)$, supported    
on a pair of pants ${\cal P}\subset \s\setminus S$ belonging to $p$, consists of changing    
only the rigid structure $r$ on ${\cal P}$. The terminology is justified by  
the fact that a braiding move may be realized by an element of ${\cal M}({\cal P})$.             
   
  
\end{enumerate}             
The 2-cells are introduced to fill in the following cycles of moves:            
\begin{enumerate}            
\item Triangles (see Figure \ref{fi7}), pentagons, squares of commutativity of            
  $A$ moves of type $(DC1)$ and $(DC2)$ (see Figure \ref{fi2}, forgetting          
  about the orientation of the circles, or Figures \ref{fi9} and   
  \ref{fi10}):          
  they all involve vertices with the same            
  support, of level $4,5,6$ or $7$, respectively.             
\item Triangles of $P$ moves: suppose that the action    
of a $P$ move on $(p,S,r)$ yields $(p,S',r')$ and consider   
another $P$ move on $(p, S',r')$ which is erasing seams            
of pairs of pants adjacent to $S$, and hence to $S'$.   
Then the composition of the two            
$P$ moves is again a $P$ move,   
and there results a triangle $PP=P$ (see Figure \ref{fi4}).            
\item Squares of commutativity of $P$ moves with $A$ moves (see Figure \ref{fi4}).             
\begin{figure}           
\begin{center}            
\includegraphics{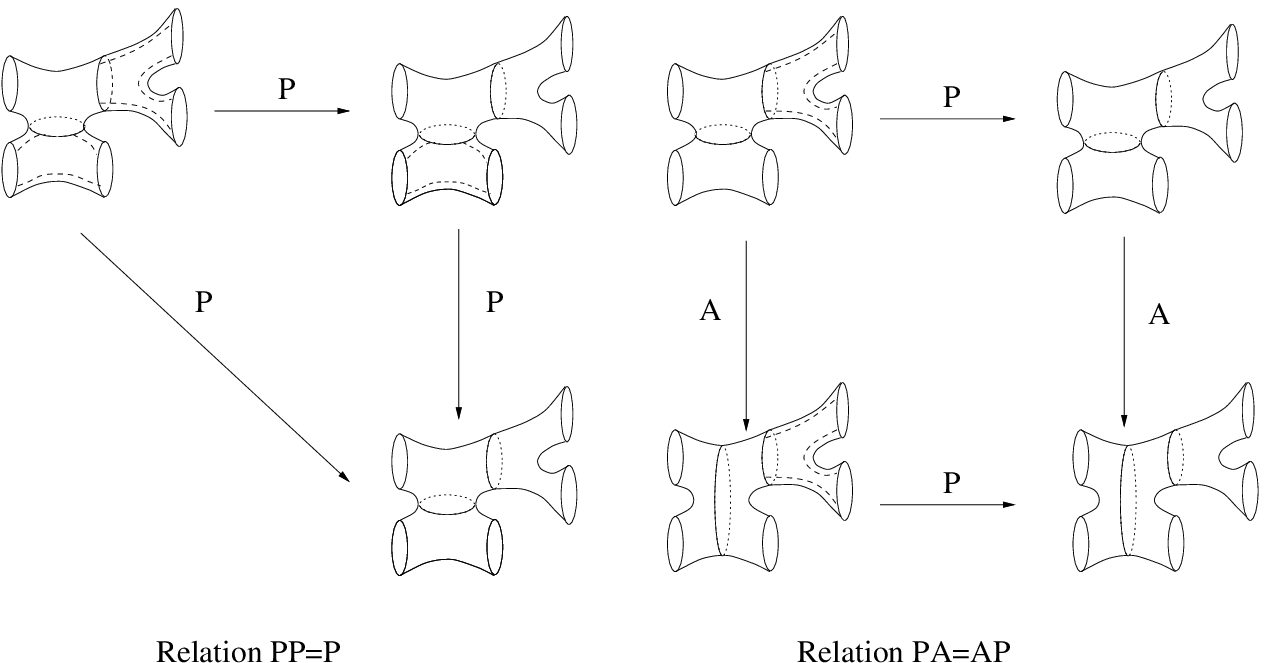}           
\caption{$2$-cells $PP=P$ and $PA=AP$}\label{fi4}            
\end{center}           
\end{figure}               
\item   
\begin{enumerate}  
\item Triangles of $Br$ moves: suppose that the action    
of a $Br$ move on $(p,S,r)$, supported on a pair of pants ${\cal  
  P}\subset \s\setminus S$ belonging to $p$, yields $(p,S',r')$, and consider   
another $Br$ move on $(p, S',r')$ supported on the same ${\cal P}$. Then the composition of the two            
$Br$ moves is again a $Br$ move,   
and there results a triangle $Br\,Br=Br$.  
  
\item Squares of commutativity of $Br$ moves supported on two distinct  
  pairs  of pants: let $(p,S,r)$ be a vertex, ${\cal P}$ and ${\cal Q}$ be two  
  pairs of pants of $\s\setminus S$, belonging to $p$. If $Br_{\cal P}$ and  
  $Br_{\cal Q}$ denote braiding moves supported on ${\cal P}$ and ${\cal Q}$,  
  respectively, then the actions of $Br_{\cal P}Br_{\cal Q}$ (i.e. $Br_{\cal  
    Q}$ followed by $Br_{\cal P}$)  
  and $Br_{\cal Q}Br_{\cal P}$ on  
  $(p,S,r)$ yield the same vertex.

\end{enumerate}      
\item Squares of            
  commutativity of $Br$ moves with $A$ moves, and squares of commutativity of            
  $Br$ moves with $P$ moves (see Figure \ref{fi5}).             
\item Triangles $P^{-1}P=Br$ (see Figure \ref{fi5}): a $P$ move erasing the rigid structure on a            
  single pair            
  of pants followed by a $P^{-1}$ move introducing a new rigid structure on            
  the same pair of pants has the same effect as a $Br$ move. (Notice that   
  these may be seen as degenerate squares of commutativity of $P$ moves with   
  $Br$ moves).                
\end{enumerate}             
\end{definition}           
           
\begin{figure}           
\includegraphics{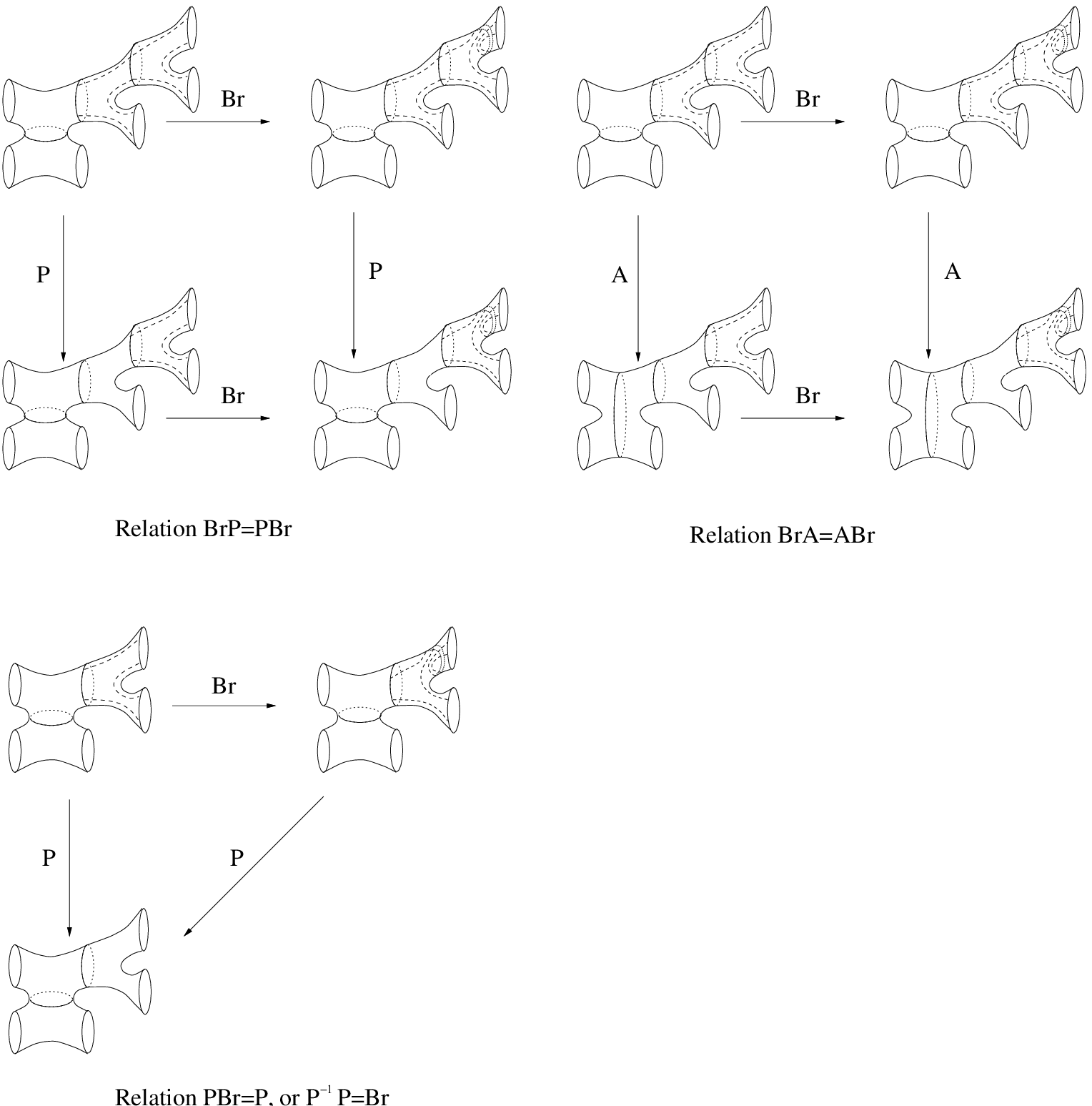}   
          
\caption{$2$-cells involving $Br$ moves}\label{fi5}       
          
\end{figure}            
                 
\begin{proposition}\label{k1}            
The complex $\widetilde{\cal DP}(\s)$ is connected and            
simply connected.            
\end{proposition}             
\noindent The proof will use twice the following lemma        
of algebraic topology from (\cite{BK2}, prop. 6.2, see also a variant   
of it in \cite{FG}):       
           
\begin{lemma}\label{sta}       
Let ${\cal M}$ and ${\cal C}$ be two $CW$-complexes of dimension            
  $2$, with oriented edges, and $f :{\cal M}^{(1)}\rightarrow {\cal C}^{(1)}$            
  be a cellular map between their 1-skeletons, surjective on $0$-cells and            
  $1$-cells. Suppose that:              
\begin{enumerate}            
\item ${\cal C}$ is connected and simply connected;            
\item For each vertex $c$ of ${\cal C}$,              
$f^{-1}(c)$ is connected and simply connected;            
\item Let $c_1\stackrel{e}{\longrightarrow } c_2$ be an oriented edge   
  of ${\cal C}$,            
  and let $m_1 '\stackrel{e'}{\longrightarrow } m_2 '$ and $m_1            
  "\stackrel{e"}{\longrightarrow } m_2"$ be two lifts in ${\cal M}$. Then            
  we can find two paths $m_1 '\stackrel{p_1}{\longrightarrow } m_1 ''$ in            
  $f^{-1}(c_1)$ and  $m_2 '\stackrel{p_2}{\longrightarrow } m_2 ''$ in            
  $f^{-1}(c_2)$ such that the loop          
            
\begin{center}        
\begin{picture}(0,0)%
\includegraphics{dia.pstex}%
\end{picture}%
\setlength{\unitlength}{2486sp}%
\begingroup\makeatletter\ifx\SetFigFont\undefined%
\gdef\SetFigFont#1#2#3#4#5{%
  \reset@font\fontsize{#1}{#2pt}%
  \fontfamily{#3}\fontseries{#4}\fontshape{#5}%
  \selectfont}%
\fi\endgroup%
\begin{picture}(1845,1650)(1396,-2041)       
\put(2926,-2041){\makebox(0,0)[lb]{\smash{\SetFigFont{8}{9.6}{\rmdefault}{\mddefault}{\updefault}$m_2''$}}}       
\put(1576,-2041){\makebox(0,0)[lb]{\smash{\SetFigFont{8}{9.6}{\rmdefault}{\mddefault}{\updefault}$m_1''$}}}       
\put(1396,-1366){\makebox(0,0)[lb]{\smash{\SetFigFont{8}{9.6}{\rmdefault}{\mddefault}{\updefault}$p_1$}}}       
\put(3241,-1366){\makebox(0,0)[lb]{\smash{\SetFigFont{8}{9.6}{\rmdefault}{\mddefault}{\updefault}$p_2$}}}       
\put(2431,-556){\makebox(0,0)[lb]{\smash{\SetFigFont{8}{9.6}{\rmdefault}{\mddefault}{\updefault}$e'$}}}       
\put(2431,-1861){\makebox(0,0)[lb]{\smash{\SetFigFont{8}{9.6}{\rmdefault}{\mddefault}{\updefault}$e''$}}}       
\put(2971,-736){\makebox(0,0)[lb]{\smash{\SetFigFont{8}{9.6}{\rmdefault}{\mddefault}{\updefault}$m_2'$}}}       
\put(1621,-736){\makebox(0,0)[lb]{\smash{\SetFigFont{8}{9.6}{\rmdefault}{\mddefault}{\updefault}$m_1'$}}}       
\end{picture}       
         
\end{center}              
is contractible in ${\cal M}$;            
\item For any $2$-cell $X$ of ${\cal C}$, its boundary $\partial X$ can be            
 lifted to a contractible loop of ${\cal M}$.            
\end{enumerate}            
Then ${\cal M}$ is connected and simply connected.           
\end{lemma}       
          
\noindent {\it Proof of Proposition \ref{k1}.} Let $F:\widetilde{\cal DP}(\s)\rightarrow            
{\cal HT}_{red}(\s)$ be the cellular map induced by the map            
$(p,S,r)\mapsto p$ on the set of vertices, which forgets about the rigid            
structure (and support $S$). The restriction of $F$ to the 1-skeletons   
is indeed    
surjective. We have proved also that ${\cal            
  HT}_{red}({\s})$ is connected and simply connected. The squares of            
commutativity of $A$ moves with $P$ moves and $Br$ moves    
insure condition  {\it 3}.    
Condition {\it 4} is satisfied by using 2-cells of type 1.    
It remains to check condition  {\it 2}.           
      
 \vspace{0.1cm}       
\noindent              
 Let $p$ be a $0$-cell of ${\cal HT}_{red}(\s)$, i.e. an            
asymptotically trivial pants decomposition of $\s$. Then            
$F^{-1}(p)$ is the subcomplex of $\widetilde{\cal DP}(\s)$ with             
vertices $(p,S,r)$  and edges corresponding to $P$ moves and   
$Br$ moves.      
        
\vspace{0.1cm}    
\noindent          
To study the connectivity of $F^{-1}(p)$, we introduce a new map    
which forgets about the pants decomposition and the rigid structure:           
$$G:{F^{-1}(p)}\rightarrow sk_2({\cal BS}_{3,7}(p))$$             
where $sk_2(X)$ denotes the 2-skeleton of the complex $X$ and    
${\cal BS}_{3,7}(p)$ is a subcomplex of a  simplicial complex            
${\cal BS}(p)$ that we now define:             
\begin{definition}[A variant of Brown-Stein's complex of bases,    
cf. \cite{br2}]            
Let $S$ and $S'$ be two surfaces compatible with $p$ (see Definition      
\ref{k2}). One says that $S'$ is a {\em simple expansion}            
of $S$ if $S'$ is the union of $S$ with a pair of pants of $p$ adjacent to $S$   
(so that $l(S')=l(S)+1$). Further $S'$ is an            
{\em expansion} of $S$ if it is obtained from $S$ by a sequence of    
simple  expansions.    
   
\vspace{0.1cm}            
\noindent    
One denotes by $S\leq S'$ when $S'$ is an expansion of $S$: this endows            
the set of $p$-compatible surfaces with a poset structure. One denotes            
by ${\cal B}(p)$ the associated simplicial complex.           
   
\vspace{0.1cm}            
\noindent  One says that $S'$            
is an elementary expansion of $S$ if $S'$ is the union of $S$ with (finitely   
many) pairs of pants of $p$ adjacent to $S$.    
   
\vspace{0.1cm}            
\noindent   
One defines the {\em elementary            
$n$-simplices} as the $(n+1)$-tuples $(U_0,\ldots,U_n)$ with            
$U_0\leq\ldots\leq U_n$ and $U_n$ an elementary expansion of $U_0$    
(so that  each $U_j$ is an elementary expansion of $U_i$,    
$0\leq i\leq j\leq n$). The set of            
all elementary simplices forms a simplicial    
subcomplex ${\cal BS}(p)$ of ${\cal B}(p)$.            
\end{definition}             
\begin{remark}                        
The complexes ${\cal BS}(p)\subset {\cal B}(p)$ are closely related to simplicial complexes introduced by K. Brown and M. Stein            
(\cite{br2}) in the study of Thompson's group $V$.    
Since the complex ${\cal B}(p)$ is associated to a direct poset, it      
is contractible.    
It happens that ${\cal BS}(p)$ is also contractible: this can            
be proved by repeating step by step M. Stein's proof of the            
contractibility of the analogous bases complex of elementary simplices            
(\cite{br2}).           
\end{remark}             
      
\noindent       
Denote by ${\cal BS}_{3,n}(p)$, $n\geq 3$, the full subcomplex of ${\cal            
  BS}(p)$ whose vertices             
are the $p$-compatible surfaces $S\subset p$ of level in $\{3,\ldots,n\}$.\\             
\begin{lemma}            
For all $n\geq 5$, ${\cal BS}_{3,n}(p)$ is contractible.            
\end{lemma}             
\begin{proof} Since ${\cal BS}(p)$ is the ascending union of the            
  subcomplexes ${\cal BS}_{3,n}(p)$, $n\geq 5$, it suffices to prove that each            
  inclusion ${\cal  BS}_{3,n}(p)\subset {\cal  BS}_{3,n+1}(p)$ is a homotopy            
  equivalence. We pass from ${\cal  BS}_{3,n}(p)$ to ${\cal  BS}_{3,n+1}(p)$   
  by considering each $p$-compatible surface    
$S$ of level $n+1$, and            
  attaching to ${\cal  BS}_{3,n}(p)$ a cone on $S$ over the link of $S$. But   
  the link, which is contained in ${\cal  BS}_{3,n}(p)$, is contractible, since            
  it can be contracted on the minimal vertex $S'$, obtained from $S$ by            
  removing all its outermost pairs of pants. Note that this is possible because           
  $n\geq 5$.            
\end{proof}             
      
\noindent       
In particular, the $2$-skeletons $sk_2({\cal BS}_{3,n}(p))$, $n\geq 5$, are            
connected and simply connected.           
       
\vspace{0.1cm}      
\noindent           
The forgetful map $G:F^{-1}(p)\rightarrow sk_2({\cal BS}_{3,7}(p))$ is induced      
on the set of vertices by the map $(p,S,r)\mapsto S$. It maps a $P$ move onto   
an edge of ${\cal BS}_{3,7}(p)$, and a 2-cell of type $PP=P$ onto a 2-simplex   
of ${\cal BS}_{3,7}(p)$. It follows from the above lemma that condition            
{\em 1} of Lemma \ref{sta} is satisfied.    
The preimage of a 0-cell is contractible in the full            
complex, since the combinatorial cycles in this preimage can be filled   
in by 2-cells of      
types {\em 4a} and {\em 4b} of Definition \ref{k2}.    
Condition {\em 4} is satisfied, since $G$ is surjective on            
the set of $2$-cells. Finally the squares of commutativity of    
$P$ moves with  $Br$ moves together with the    
triangles $P^{-1}P=Br$ imply condition {\em 3}. This completes            
the proof that $F^{-1}(p)$ is connected and simply connected, and   
thus completes the proof            
of Proposition \ref{k1}. $\;\;\;\square$

\subsection{The complex ${\cal DP}(\s)$ and the main theorem for $\B$}             
\begin{definition}\label{k3}            
Let ${\cal DP}(\s)$ be the subcomplex of $\widetilde{\cal            
  DP}(\s)$ obtained from the latter by eliminating all edges and    
$2$-cells  involving $Br$ moves. In            
  other words, ${\cal DP}(\s)$ is the cellular complex whose            
  vertices are triples $(p,S,r)$ as in Definition \ref{k2}, with edges corresponding to            
  $A$ moves and $P$ moves, and three distinct families of $2$-cells:            
\begin{enumerate}            
\item Hatcher-Thurston-type cells: triangles, pentagons, and    
squares $(DC1)$ and   $(DC2)$;            
\item Brown-Stein-type cells ($2$-simplices): $PP=P$;            
\item Squares of commutativity $PA=AP$.        
\end{enumerate}            
\end{definition}    
\noindent Notice that the cells of third kind are  connecting    
Hatcher-Thurston-type cells  to Brown-Stein-type cells.                   
\begin{theorem}            
The complex ${\cal DP}(\s)$ is a connected and simply connected            
${\cal B}$-complex, with finite quotient ${\cal DP}(\s)/{\cal            
  B}$. In particular, ${\cal B}$ is a finitely presented group.             
\end{theorem}             
\begin{proof} By adapting the arguments of the proof of Proposition \ref{HT},   
  one proves that ${\cal DP}(\s)$ has a finite number of cells modulo   
  $\B$. The last assertion is therefore a direct application of K. Brown's   
theorem on  presentations for groups acting on simply connected    
$CW$-complexes, in the case of finitely presented vertex stabilizers    
and finitely generated edge stabilizers   
(which is the case here, cf. the next subsection). To            
prove the first assertion, we show that $\widetilde{\cal            
  DP}(\s)$ is obtained from ${\cal  DP}(\s)$ by            
attaching countably many $2$-spheres at each $0$-cell of ${\cal            
  DP}(\s)$. Our strategy is to examine below cells of type {\it 4a} in  
step 1, {\it 4b} in step 2, and {\it 5-6} in step 3.

\begin{enumerate}  
\item   Consider first a $2$-cell $\omega_a$ bounded by a triangle of $Br$  
  moves, all supported on the same pair of pants ${\cal P}\subset \s\setminus S$, with vertices  
  $(p,S,r_i)$, $i=1,2$ and 3,  cf. Definition \ref{k2}, {4a)}. The rigid structures $r_i$ only differ from each other on the pants            
${\cal P}$. We may find a sequence of $P$ moves and            
$P^{-1}$ moves $P_1^{\pm 1}P_2^{\pm 1}\cdots P_n^{\pm 1}$    
translating $S$ to another $p$-compatible surface $\tilde{S}$    
containing ${\cal P}$, since $l(S)\geq            
l({\cal P})$. By using the squares of commutativity $P_i Br=Br P_i$    
and the  degenerate squares (or triangles) $P_i Br=P_i$, the    
initial triangle of  $Br$ moves is replaced,    
$P^{\pm 1}$ move by $P^{\pm 1}$ move, by a            
conjugate triangle. The edge-paths of $P^{\pm 1}$ moves starting from the            
3 different vertices $(p,S,r_i)$ end at the {\it same} vertex            
$(p,\tilde{S},\tilde{r})$ (since $\tilde{S}$ contains the pants    
${\cal P}$ on which            
the 3 rigid structures differ from each other). Thus, adding the    
$2$-cell $\omega_a$ to            
${\cal DP}(\s)$ is homotopically equivalent to attaching            
$2$-spheres at the vertex $(p,\tilde{S}, \tilde{r})$.        
      
\item  Consider now a $2$-cell $\omega_b$ bounded by a square cycle of            
   commutativity between two $Br$ moves supported on two distinct pairs of            
   pants ${\cal P}$ and ${\cal P}'$ (subordinate to the same $(p,S)$,            
   cf. Definition \ref{k2}, {4b)}). The four vertices of the square correspond to            
   pairs $(p,S,r_i)$, $i=1,\ldots,4$, the pants decomposition $p$ being            
   fixed as well as the support $S$, but the rigid structure $r_i$ differing            
   from each other only on ${\cal P}$ and/or ${\cal P}'$.  
  
\vspace{0.1cm}   
\noindent           
Let $\Sigma_{{\cal P P}'}$ be the minimal compact subsurface,   
   compatible with $p$, connecting  ${\cal P}$ to ${\cal P}'$ (see Figure   
   \ref{atmost5}). By changing the pants decomposition $p$ inside   
   $\Sigma_{{\cal P P}'}$, we may find a pants decomposition $\tilde{p}$ such that ${\cal P}$ and            
   ${\cal P}'$ (which are still compatible with $\tilde{p}$) belong to a same   
   surface $\tilde{S}$ {\it of level at most 5}, compatible with $\tilde{p}$ (see Figure \ref{atmost5}).   
   
\begin{figure}   
\begin{center}   
\begin{picture}(0,0)%
\includegraphics{atmost5.pstex}%
\end{picture}%
\setlength{\unitlength}{2693sp}%
\begingroup\makeatletter\ifx\SetFigFont\undefined%
\gdef\SetFigFont#1#2#3#4#5{%
  \reset@font\fontsize{#1}{#2pt}%
  \fontfamily{#3}\fontseries{#4}\fontshape{#5}%
  \selectfont}%
\fi\endgroup%
\begin{picture}(8326,4985)(338,-4170)   
\put(4051,-4170){\makebox(0,0)[lb]{\smash{\SetFigFont{8}{9.6}{\rmdefault}{\mddefault}{\updefault}$\tilde{S}$ (of level 5)}}}   
\put(1171,-390){\makebox(0,0)[lb]{\smash{\SetFigFont{8}{9.6}{\rmdefault}{\mddefault}{\updefault}${\cal P}$}}}   
\put(1216,-3225){\makebox(0,0)[lb]{\smash{\SetFigFont{8}{9.6}{\rmdefault}{\mddefault}{\updefault}${\cal P}$}}}   
\put(7651,-390){\makebox(0,0)[lb]{\smash{\SetFigFont{8}{9.6}{\rmdefault}{\mddefault}{\updefault}${\cal P}'$}}}   
\put(7651,-3225){\makebox(0,0)[lb]{\smash{\SetFigFont{8}{9.6}{\rmdefault}{\mddefault}{\updefault}${\cal P}'$}}}   
\put(4366,-1380){\makebox(0,0)[lb]{\smash{\SetFigFont{8}{9.6}{\rmdefault}{\mddefault}{\updefault}${\Sigma}_{{\cal PP}'}$}}}   
\put(3556,659){\makebox(0,0)[lb]{\smash{\SetFigFont{8}{9.6}{\rmdefault}{\mddefault}{\updefault}Pants decomposition $p$}}}   
\put(3511,-2221){\makebox(0,0)[lb]{\smash{\SetFigFont{8}{9.6}{\rmdefault}{\mddefault}{\updefault}Pants decomposition $\tilde{p}$}}}   
\end{picture}

\caption{From $p$ to $\tilde{p}$ }\label{atmost5}   
\end{center}    
\end{figure}   
            
\noindent Consider next the pants decompositions $p$ and $\tilde{p}$ in the Hatcher-Thurston            
   complex of $\s$. Since they only differ on ${\Sigma}_{{\cal PP}'}$, the connectedness of the Hatcher-Thurston complex of ${\Sigma}_{{\cal PP}'}$ enables us to find            
   a sequence of $A$ moves $A_1 A_2 \cdots A_n$ connecting $p$ to $\tilde{p}$,   
   and all supported in ${\Sigma}_{{\cal PP}'}$. Now an $A$ move in            
   ${\cal DP}(\s)$ is            
   possible only inside a support devoid of its seams. Therefore, we lift to ${\cal DP}(\s)$ the            
   sequence            
 of moves $A_1 A_2 \cdots A_n$ into four edge paths starting from            
   each $(p,S,r_i)$, by inserting  an adequate number of    
 $P^{\pm 1}$ moves between the $A$ moves            
of the sequence $A_1 A_2 \cdots A_n$. The four edge paths end at a same vertex            
   $(\tilde{p}, \tilde{S}, \tilde{r})$, and we conclude as in 1), using the squares            
   $PBr=BrP$ and $ABr=BrA$.   
       
\item  Notice finally that each square $PBr=BrP$, $ABr=BrA$, belongs to a            
   2-sphere obtained in 1) or 2), and thus can be eliminated with the cells            
   $\omega_a$ and $\omega_b$.            
   
 \end{enumerate}  
\end{proof}

\section{A finite presentation for the group ${\cal B}$}         
           
\subsection{Reduction of ${\cal DP}(\s)$ to the complex ${\cal DP}^+_5(\s)$}                 
Consider the full subcomplex ${\cal DP}_5(\s)$ whose levels of            
vertices belong to $\{3,4,5\}$. The boundaries of the $2$-cells of type $(DC1)$   
and $(DC2)$ in ${\cal DP}_5(\s)$ do not belong to ${\cal DP}_5(\s)$, but are   
homotopic to cycles ${\cal C}_1$ and ${\cal C}_2$ of ${\cal DP}_5(\s)$,   
respectively (cf. Figure \ref{fi9} and Figure \ref{fi10}). Indeed, the ``area"   
between the boundary of $(DC1)$ and the cycle ${\cal C}_1$ (resp. $(DC2)$ and   
the cycle ${\cal C}_2$) may be filled in by squares $PA=AP$ and simplices   
$PP=P$.  
  
\vspace{0.1cm}   
\noindent   
Note that the cycles ${\cal C}_1$ (resp. ${\cal C}_2$), which are of length   
$12$ (resp. of length $20$), are all            
equivalent modulo ${\cal B}$.

\begin{definition}            
Let ${\cal DP}^+_5(\s)$ be the complex obtained from ${\cal            
  DP}_5(\s)$ by adding $2$-cells to fill in the cycles of types             
${\cal C}_1$ and ${\cal C}_2$.            
\end{definition}             
\begin{proposition}            
The ${\cal B}$-complex ${\cal DP}^+_5(\s)$ is connected and simply            
connected, and finite modulo ${\cal B}$.            
\end{proposition}             
\begin{proof} From Figures \ref{fi9} and \ref{fi10}, each cell of type (DC1) or (DC2) gives rise to countably            
many ``big squares", which are therefore homeomorphic to $2$-cells, attached to ${\cal DP}^+_5(\s)$ along            
contractible cycles ${\cal C}_1$ or ${\cal C}_2$. They are made up of            
squares (DC1) or (DC2), squares $PA=AP$, and simplices $PP=P$.   
   
\vspace{0.1cm}   
\noindent Similarly, consider the  other $2$-cells of the closure of ${\cal DP}(\s)\setminus {\cal            
  DP}_5(\s)$: the Hatcher-Thurston triangles and            
pentagons (on supports of levels $6$ or $7$) are attached along similar            
triangular and pentagonal cycles of ${\cal DP}_5(\s)$, via            
squares $PA=AP$; the simplices $PP=P$ also are attached along contractible cycles            
via other simplices $PP=P$.      
      
 \vspace{0.1cm}       
\noindent            
All squares $PA=AP$ and simplices $PP=P$ have been used in proceeding to the            
above attachments. It follows that the closure of ${\cal DP}(\s)\setminus {\cal            
  DP}_5(\s)$ is homeomorphic to a union of closed $2$-cells,            
all attached to ${\cal DP}^+_5(\s)$ along contractible            
cycles.            
\end{proof}
              
\begin{figure}           
\begin{center}            
\includegraphics{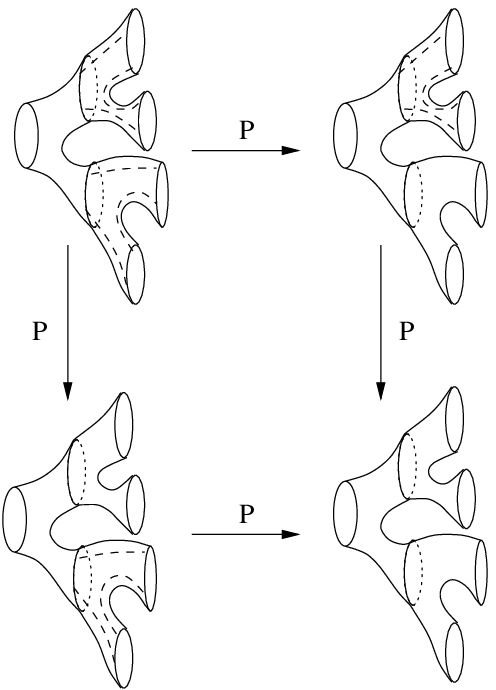}  
          
\caption{Square $PP=PP$}\label{fi6}           
\end{center}            
\end{figure}

\noindent{\bf Set ${\cal F}$ of essential $2$-cells:} Notice first that the triangles of            
${\cal DP}^+_5(\s)$ on supports of level $5$ are conjugate, via            
$P$-edges, to triangles on supports of level $4$. Further, replace the            
simplices $PP=P$ by squares $PP=PP$ involving two adjacent simplices            
(cf. Figure \ref{fi6}).   
             
This does not affect the homeomorphism type of ${\cal            
  DP}^+_5(\s)$, and by abuse of notation, we still denote the            
modified complex by ${\cal DP}^+_5(\s)$. It follows that there            
are 6 essential 2-cells in ${\cal DP}^+_5(\s)$ modulo ${\cal            
  B}$ (cf. Figures \ref{fi7}, \ref{fi8}, \ref{fi9} and \ref{fi10}).                  
\begin{figure}   
\begin{center}           
\begin{picture}(0,0)%
\includegraphics{penta.pstex}%
\end{picture}%
\setlength{\unitlength}{2693sp}%
\begingroup\makeatletter\ifx\SetFigFont\undefined%
\gdef\SetFigFont#1#2#3#4#5{%
  \reset@font\fontsize{#1}{#2pt}%
  \fontfamily{#3}\fontseries{#4}\fontshape{#5}%
  \selectfont}%
\fi\endgroup%
\begin{picture}(5109,11775)(1339,-10928)  
\put(1801,569){\makebox(0,0)[lb]{\smash{\SetFigFont{8}{9.6}{\familydefault}{\mddefault}{\updefault}$v_3$}}}  
\put(1891,-6721){\makebox(0,0)[lb]{\smash{\SetFigFont{8}{9.6}{\familydefault}{\mddefault}{\updefault}$v_2$}}}  
\put(5491,-6721){\makebox(0,0)[lb]{\smash{\SetFigFont{8}{9.6}{\familydefault}{\mddefault}{\updefault}$\alpha(v_2)$ }}}  
\put(2296,-10861){\makebox(0,0)[lb]{\smash{\SetFigFont{10}{12.0}{\rmdefault}{\bfdefault}{\updefault}Relation 2 (b):  $(\alpha\pi)^3=t_2^{-1}$}}}  
\put(2476,-5641){\makebox(0,0)[lb]{\smash{\SetFigFont{10}{12.0}{\rmdefault}{\bfdefault}{\updefault}Relation 2 (a): $(\beta\alpha)^5=1$}}}  
\put(2386,-1276){\makebox(0,0)[lb]{\smash{\SetFigFont{8}{9.6}{\rmdefault}{\mddefault}{\updefault}1}}}  
\put(5086,-4201){\makebox(0,0)[lb]{\smash{\SetFigFont{8}{9.6}{\rmdefault}{\mddefault}{\updefault}$g$}}}  
\put(3601,-646){\makebox(0,0)[lb]{\smash{\SetFigFont{8}{9.6}{\rmdefault}{\mddefault}{\updefault}5}}}  
\put(5356,-1366){\makebox(0,0)[lb]{\smash{\SetFigFont{8}{9.6}{\rmdefault}{\mddefault}{\updefault}4}}}  
\put(3601,-241){\makebox(0,0)[lb]{\smash{\SetFigFont{8}{9.6}{\rmdefault}{\mddefault}{\updefault}$g$}}}  
\put(5581,-1096){\makebox(0,0)[lb]{\smash{\SetFigFont{8}{9.6}{\rmdefault}{\mddefault}{\updefault}$g$}}}  
\put(2071,-1096){\makebox(0,0)[lb]{\smash{\SetFigFont{8}{9.6}{\rmdefault}{\mddefault}{\updefault}$g$}}}  
\put(2836,-3841){\makebox(0,0)[lb]{\smash{\SetFigFont{8}{9.6}{\rmdefault}{\mddefault}{\updefault}2}}}  
\put(4861,-3931){\makebox(0,0)[lb]{\smash{\SetFigFont{8}{9.6}{\rmdefault}{\mddefault}{\updefault}3}}}  
\put(2566,-4111){\makebox(0,0)[lb]{\smash{\SetFigFont{8}{9.6}{\rmdefault}{\mddefault}{\updefault}$g$}}}  
\put(4816,-8161){\makebox(0,0)[lb]{\smash{\SetFigFont{8}{9.6}{\rmdefault}{\mddefault}{\updefault}$\pi\alpha$}}}  
\put(3691,-6946){\makebox(0,0)[lb]{\smash{\SetFigFont{8}{9.6}{\rmdefault}{\mddefault}{\updefault}1}}}  
\put(3691,-6586){\makebox(0,0)[lb]{\smash{\SetFigFont{8}{9.6}{\rmdefault}{\mddefault}{\updefault}$\alpha$}}}  
\put(4411,-8161){\makebox(0,0)[lb]{\smash{\SetFigFont{8}{9.6}{\rmdefault}{\mddefault}{\updefault}2}}}  
\put(3061,-8161){\makebox(0,0)[lb]{\smash{\SetFigFont{8}{9.6}{\rmdefault}{\mddefault}{\updefault}3}}}  
\put(2386,-8161){\makebox(0,0)[lb]{\smash{\SetFigFont{8}{9.6}{\rmdefault}{\mddefault}{\updefault}$\pi\alpha$}}}  
\end{picture}

 \caption{Relations associated to the pentagon and the triangle}\label{fi7}           
\end{center}           
\end{figure}             
                 
\begin{figure}           
\begin{center}            
\begin{picture}(0,0)%
\includegraphics{rel1.pstex}%
\end{picture}%
\setlength{\unitlength}{3108sp}%
\begingroup\makeatletter\ifx\SetFigFont\undefined%
\gdef\SetFigFont#1#2#3#4#5{%
  \reset@font\fontsize{#1}{#2pt}%
  \fontfamily{#3}\fontseries{#4}\fontshape{#5}%
  \selectfont}%
\fi\endgroup%
\begin{picture}(7740,5241)(586,-4561)  
\put(1351,299){\makebox(0,0)[lb]{\smash{\SetFigFont{9}{10.8}{\rmdefault}{\mddefault}{\updefault}$v_3$}}}  
\put(856,-3931){\makebox(0,0)[lb]{\smash{\SetFigFont{9}{10.8}{\rmdefault}{\bfdefault}{\updefault}Relation $g^{-1}\alpha^{-1}\pi^{-1}= \beta\pi^{-1}$ in ${\cal B}_{v_3}$}}}  
\put(2476,524){\makebox(0,0)[lb]{\smash{\SetFigFont{9}{10.8}{\rmdefault}{\bfdefault}{\updefault}$g^{-1}$}}}  
\put(2341,-196){\makebox(0,0)[lb]{\smash{\SetFigFont{9}{10.8}{\rmdefault}{\mddefault}{\updefault}$A$}}}  
\put(5176,-2806){\makebox(0,0)[lb]{\smash{\SetFigFont{9}{10.8}{\rmdefault}{\mddefault}{\updefault}$v_2$}}}  
\put(5131,-3931){\makebox(0,0)[lb]{\smash{\SetFigFont{9}{10.8}{\rmdefault}{\bfdefault}{\updefault}Relation $\pi^{\beta}$ in ${\cal B}_{v_3}$ = $\beta^{-1}\pi\beta$ in ${\cal B}_{v_1}$}}}  
\put(5716,-4561){\makebox(0,0)[lb]{\smash{\SetFigFont{9}{10.8}{\rmdefault}{\mddefault}{\updefault}Commutation PP=PP}}}  
\put(1441,-4561){\makebox(0,0)[lb]{\smash{\SetFigFont{9}{10.8}{\rmdefault}{\mddefault}{\updefault}Commutation PA=AP}}}  
\put(2386,-511){\makebox(0,0)[lb]{\smash{\SetFigFont{9}{10.8}{\rmdefault}{\mddefault}{\updefault}1}}}  
\put(5941,-1321){\makebox(0,0)[lb]{\smash{\SetFigFont{9}{10.8}{\rmdefault}{\mddefault}{\updefault}1}}}  
\put(2296,-3481){\makebox(0,0)[lb]{\smash{\SetFigFont{9}{10.8}{\rmdefault}{\bfdefault}{\updefault}$\alpha^{-1}$}}}  
\put(6571,-2806){\makebox(0,0)[lb]{\smash{\SetFigFont{9}{10.8}{\rmdefault}{\mddefault}{\updefault}$P$}}}  
\put(8326,-1321){\makebox(0,0)[lb]{\smash{\SetFigFont{9}{10.8}{\rmdefault}{\mddefault}{\updefault}$\pi^{\beta}$}}}  
\put(4906,-1321){\makebox(0,0)[lb]{\smash{\SetFigFont{9}{10.8}{\rmdefault}{\mddefault}{\updefault}1}}}  
\put(2431,-2716){\makebox(0,0)[lb]{\smash{\SetFigFont{9}{10.8}{\rmdefault}{\mddefault}{\updefault}$A$}}}  
\put(2476,-2401){\makebox(0,0)[lb]{\smash{\SetFigFont{9}{10.8}{\rmdefault}{\mddefault}{\updefault}3}}}  
\put(2971,-1321){\makebox(0,0)[lb]{\smash{\SetFigFont{9}{10.8}{\rmdefault}{\mddefault}{\updefault}2}}}  
\put(586,-1321){\makebox(0,0)[lb]{\smash{\SetFigFont{9}{10.8}{\rmdefault}{\bfdefault}{\updefault}$\pi^{-1}$}}}  
\put(3286,-1321){\makebox(0,0)[lb]{\smash{\SetFigFont{9}{10.8}{\rmdefault}{\mddefault}{\updefault}$P^{-1}$}}}  
\put(4096,-1321){\makebox(0,0)[lb]{\smash{\SetFigFont{9}{10.8}{\rmdefault}{\bfdefault}{\updefault}$1$}}}  
\put(1441,-1321){\makebox(0,0)[lb]{\smash{\SetFigFont{9}{10.8}{\rmdefault}{\mddefault}{\updefault}$P$}}}  
\put(1801,-1321){\makebox(0,0)[lb]{\smash{\SetFigFont{9}{10.8}{\rmdefault}{\mddefault}{\updefault}4}}}  
\put(7291,-1321){\makebox(0,0)[lb]{\smash{\SetFigFont{9}{10.8}{\rmdefault}{\mddefault}{\updefault}3}}}  
\put(7606,-1321){\makebox(0,0)[lb]{\smash{\SetFigFont{9}{10.8}{\rmdefault}{\mddefault}{\updefault}$P^{-1}$}}}  
\put(5536,-1321){\makebox(0,0)[lb]{\smash{\SetFigFont{9}{10.8}{\rmdefault}{\mddefault}{\updefault}$P$}}}  
\put(6661,-2356){\makebox(0,0)[lb]{\smash{\SetFigFont{9}{10.8}{\rmdefault}{\mddefault}{\updefault}2}}}  
\put(6616,-646){\makebox(0,0)[lb]{\smash{\SetFigFont{9}{10.8}{\rmdefault}{\mddefault}{\updefault}4}}}  
\put(6526,-286){\makebox(0,0)[lb]{\smash{\SetFigFont{9}{10.8}{\rmdefault}{\mddefault}{\updefault}$P^{-1}$}}}  
\put(5221,-421){\makebox(0,0)[lb]{\smash{\SetFigFont{9}{10.8}{\rmdefault}{\mddefault}{\updefault}$v_1$}}}  
\put(6706,-3391){\makebox(0,0)[lb]{\smash{\SetFigFont{9}{10.8}{\rmdefault}{\mddefault}{\updefault}1}}}  
\end{picture}  
  
\caption{Relations associated to the cells $PA=AP$ and $PP=PP$}\label{fi8}            
\end{center}            
\end{figure}

\begin{figure}       
\begin{center}            
\begin{picture}(0,0)%
\includegraphics{C1.pstex}%
\end{picture}%
\setlength{\unitlength}{2486sp}%
\begingroup\makeatletter\ifx\SetFigFont\undefined%
\gdef\SetFigFont#1#2#3#4#5{%
  \reset@font\fontsize{#1}{#2pt}%
  \fontfamily{#3}\fontseries{#4}\fontshape{#5}%
  \selectfont}%
\fi\endgroup%
\begin{picture}(8460,11203)(856,-10703)   
\put(5086,-3256){\makebox(0,0)[lb]{\smash{\SetFigFont{8}{9.6}{\rmdefault}{\mddefault}{\updefault}$A$}}}   
\put(5086,-6271){\makebox(0,0)[lb]{\smash{\SetFigFont{8}{9.6}{\rmdefault}{\mddefault}{\updefault}$A$}}}   
\put(6436,-4696){\makebox(0,0)[lb]{\smash{\SetFigFont{8}{9.6}{\rmdefault}{\mddefault}{\updefault}$A$}}}   
\put(3781,-4741){\makebox(0,0)[lb]{\smash{\SetFigFont{8}{9.6}{\rmdefault}{\mddefault}{\updefault}$A$}}}   
\put(4951,-4741){\makebox(0,0)[lb]{\smash{\SetFigFont{8}{9.6}{\rmdefault}{\mddefault}{\updefault}$DC1$}}}   
\put(2071,-4786){\makebox(0,0)[lb]{\smash{\SetFigFont{8}{9.6}{\rmdefault}{\mddefault}{\updefault}$A$}}}   
\put(1711,-4786){\makebox(0,0)[lb]{\smash{\SetFigFont{8}{9.6}{\rmdefault}{\mddefault}{\updefault}$1$}}}   
\put(1666,-2176){\makebox(0,0)[lb]{\smash{\SetFigFont{8}{9.6}{\rmdefault}{\mddefault}{\updefault}$2$}}}   
\put(2071,-2176){\makebox(0,0)[lb]{\smash{\SetFigFont{8}{9.6}{\rmdefault}{\mddefault}{\updefault}$P^{-1}$}}}   
\put(7381,-871){\makebox(0,0)[lb]{\smash{\SetFigFont{8}{9.6}{\rmdefault}{\mddefault}{\updefault}$P^{-1}$}}}   
\put(7381,-466){\makebox(0,0)[lb]{\smash{\SetFigFont{8}{9.6}{\rmdefault}{\mddefault}{\updefault}$5$}}}   
\put(8056,-4741){\makebox(0,0)[lb]{\smash{\SetFigFont{8}{9.6}{\rmdefault}{\mddefault}{\updefault}$A$}}}   
\put(8596,-4741){\makebox(0,0)[lb]{\smash{\SetFigFont{8}{9.6}{\rmdefault}{\mddefault}{\updefault}$7$}}}   
\put(8596,-7396){\makebox(0,0)[lb]{\smash{\SetFigFont{8}{9.6}{\rmdefault}{\mddefault}{\updefault}$8$}}}   
\put(7111,-8656){\makebox(0,0)[lb]{\smash{\SetFigFont{8}{9.6}{\rmdefault}{\mddefault}{\updefault}$P$}}}   
\put(2926,-8656){\makebox(0,0)[lb]{\smash{\SetFigFont{8}{9.6}{\rmdefault}{\mddefault}{\updefault}$P^{-1}$}}}   
\put(2116,-7531){\makebox(0,0)[lb]{\smash{\SetFigFont{8}{9.6}{\rmdefault}{\mddefault}{\updefault}$P$}}}   
\put(5086,-9061){\makebox(0,0)[lb]{\smash{\SetFigFont{8}{9.6}{\rmdefault}{\mddefault}{\updefault}$10$}}}   
\put(5086,-8656){\makebox(0,0)[lb]{\smash{\SetFigFont{8}{9.6}{\rmdefault}{\mddefault}{\updefault}$A$}}}   
\put(7156,-9061){\makebox(0,0)[lb]{\smash{\SetFigFont{8}{9.6}{\rmdefault}{\mddefault}{\updefault}$9$}}}   
\put(9316,-4696){\makebox(0,0)[lb]{\smash{\SetFigFont{8}{9.6}{\rmdefault}{\bfdefault}{\updefault}$g^{-1}$}}}   
\put(9316,-7351){\makebox(0,0)[lb]{\smash{\SetFigFont{8}{9.6}{\rmdefault}{\bfdefault}{\updefault}$\pi^{\beta}$}}}   
\put(7021,-9781){\makebox(0,0)[lb]{\smash{\SetFigFont{8}{9.6}{\rmdefault}{\bfdefault}{\updefault}$\alpha^2$}}}   
\put(8146,-2131){\makebox(0,0)[lb]{\smash{\SetFigFont{8}{9.6}{\rmdefault}{\mddefault}{\updefault}$P$}}}   
\put(8641,-2131){\makebox(0,0)[lb]{\smash{\SetFigFont{8}{9.6}{\rmdefault}{\mddefault}{\updefault}$6$}}}   
\put(9316,-2131){\makebox(0,0)[lb]{\smash{\SetFigFont{8}{9.6}{\rmdefault}{\bfdefault}{\updefault}$\alpha^2$}}}   
\put(3286,299){\makebox(0,0)[lb]{\smash{\SetFigFont{8}{9.6}{\rmdefault}{\bfdefault}{\updefault}$\alpha^2$}}}   
\put(3286,-466){\makebox(0,0)[lb]{\smash{\SetFigFont{8}{9.6}{\rmdefault}{\mddefault}{\updefault}$3$}}}   
\put(3286,-871){\makebox(0,0)[lb]{\smash{\SetFigFont{8}{9.6}{\rmdefault}{\mddefault}{\updefault}$P$}}}   
\put(5446,344){\makebox(0,0)[lb]{\smash{\SetFigFont{8}{9.6}{\rmdefault}{\bfdefault}{\updefault}$g^{-1}$}}}   
\put(5446,-466){\makebox(0,0)[lb]{\smash{\SetFigFont{8}{9.6}{\rmdefault}{\mddefault}{\updefault}$4$}}}   
\put(5446,-871){\makebox(0,0)[lb]{\smash{\SetFigFont{8}{9.6}{\rmdefault}{\mddefault}{\updefault}$A$}}}   
\put(7381,344){\makebox(0,0)[lb]{\smash{\SetFigFont{8}{9.6}{\rmdefault}{\bfdefault}{\updefault}$\pi^{\beta}$}}}   
\put(7741,-7396){\makebox(0,0)[lb]{\smash{\SetFigFont{8}{9.6}{\rmdefault}{\mddefault}{\updefault}$P^{-1}$}}}   
\put(2926,-9061){\makebox(0,0)[lb]{\smash{\SetFigFont{8}{9.6}{\rmdefault}{\mddefault}{\updefault}$11$}}}   
\put(2926,-9826){\makebox(0,0)[lb]{\smash{\SetFigFont{8}{9.6}{\rmdefault}{\bfdefault}{\updefault}$1$}}}   
\put(1576,-7531){\makebox(0,0)[lb]{\smash{\SetFigFont{8}{9.6}{\rmdefault}{\mddefault}{\updefault}$12$}}}   
\put(2386,-10636){\makebox(0,0)[lb]{\smash{\SetFigFont{10}{12.0}{\rmdefault}{\bfdefault}{\updefault}Relation 6 (a): $(\alpha^2\pi^{\beta}\alpha^3\beta^{2})^2=(\pi^{\beta}\alpha^3\beta^2\alpha^2)^2$}}}   
\put(856,-2176){\makebox(0,0)[lb]{\smash{\SetFigFont{8}{9.6}{\rmdefault}{\bfdefault}{\updefault}$1$}}}   
\put(856,-4786){\makebox(0,0)[lb]{\smash{\SetFigFont{8}{9.6}{\rmdefault}{\bfdefault}{\updefault}$g$}}}   
\put(901,-7531){\makebox(0,0)[lb]{\smash{\SetFigFont{8}{9.6}{\rmdefault}{\bfdefault}{\updefault}$\alpha^2$}}}   
\put(4816,-9826){\makebox(0,0)[lb]{\smash{\SetFigFont{8}{9.6}{\rmdefault}{\bfdefault}{\updefault}$\pi^{\beta}g$}}}   
\end{picture}

\caption{Relation associated to the cell $(DC1)$}\label{fi9}           
\end{center}           
\end{figure}

\begin{figure}           
\begin{center}            
\begin{picture}(0,0)%
\includegraphics{C2.pstex}%
\end{picture}%
\setlength{\unitlength}{2486sp}%
\begingroup\makeatletter\ifx\SetFigFont\undefined%
\gdef\SetFigFont#1#2#3#4#5{%
  \reset@font\fontsize{#1}{#2pt}%
  \fontfamily{#3}\fontseries{#4}\fontshape{#5}%
  \selectfont}%
\fi\endgroup%
\begin{picture}(8370,13642)(541,-13673)   
\put(1981,-2986){\makebox(0,0)[lb]{\smash{\SetFigFont{8}{9.6}{\rmdefault}{\mddefault}{\updefault}$P$}}}   
\put(1981,-3391){\makebox(0,0)[lb]{\smash{\SetFigFont{8}{9.6}{\rmdefault}{\mddefault}{\updefault}$P$}}}   
\put(2971,-11086){\makebox(0,0)[lb]{\smash{\SetFigFont{8}{9.6}{\rmdefault}{\mddefault}{\updefault}$P$}}}   
\put(2521,-11086){\makebox(0,0)[lb]{\smash{\SetFigFont{8}{9.6}{\rmdefault}{\mddefault}{\updefault}$P$}}}   
\put(7516,-9511){\makebox(0,0)[lb]{\smash{\SetFigFont{8}{9.6}{\rmdefault}{\mddefault}{\updefault}$P$}}}   
\put(7516,-9871){\makebox(0,0)[lb]{\smash{\SetFigFont{8}{9.6}{\rmdefault}{\mddefault}{\updefault}$P$}}}   
\put(4726,-6541){\makebox(0,0)[lb]{\smash{\SetFigFont{8}{9.6}{\rmdefault}{\mddefault}{\updefault}$DC2$}}}   
\put(4771,-1771){\makebox(0,0)[lb]{\smash{\SetFigFont{8}{9.6}{\rmdefault}{\mddefault}{\updefault}$A$}}}   
\put(7381,-6406){\makebox(0,0)[lb]{\smash{\SetFigFont{8}{9.6}{\rmdefault}{\mddefault}{\updefault}$A$}}}   
\put(1981,-6451){\makebox(0,0)[lb]{\smash{\SetFigFont{8}{9.6}{\rmdefault}{\mddefault}{\updefault}$A$}}}   
\put(4726,-11086){\makebox(0,0)[lb]{\smash{\SetFigFont{8}{9.6}{\rmdefault}{\mddefault}{\updefault}$A$}}}   
\put(1261,-3211){\makebox(0,0)[lb]{\smash{\SetFigFont{8}{9.6}{\rmdefault}{\mddefault}{\updefault}$9/10$}}}   
\put(8911,-6406){\makebox(0,0)[lb]{\smash{\SetFigFont{8}{9.6}{\rmdefault}{\bfdefault}{\updefault}$g$}}}   
\put(8056,-6406){\makebox(0,0)[lb]{\smash{\SetFigFont{8}{9.6}{\rmdefault}{\mddefault}{\updefault}$1$}}}   
\put(7966,-3166){\makebox(0,0)[lb]{\smash{\SetFigFont{8}{9.6}{\rmdefault}{\mddefault}{\updefault}$2/3$}}}   
\put(7291,-196){\makebox(0,0)[lb]{\smash{\SetFigFont{8}{9.6}{\rmdefault}{\bfdefault}{\updefault}$\beta$}}}   
\put(586,-3031){\makebox(0,0)[lb]{\smash{\SetFigFont{8}{9.6}{\rmdefault}{\bfdefault}{\updefault}$\beta^{-1}$}}}   
\put(586,-3436){\makebox(0,0)[lb]{\smash{\SetFigFont{8}{9.6}{\rmdefault}{\bfdefault}{\updefault}$\alpha^2$}}}   
\put(1396,-6451){\makebox(0,0)[lb]{\smash{\SetFigFont{8}{9.6}{\rmdefault}{\mddefault}{\updefault}$11$}}}   
\put(541,-9421){\makebox(0,0)[lb]{\smash{\SetFigFont{8}{9.6}{\rmdefault}{\bfdefault}{\updefault}$(\pi^{\beta})^{-1}$}}}   
\put(541,-9961){\makebox(0,0)[lb]{\smash{\SetFigFont{8}{9.6}{\rmdefault}{\bfdefault}{\updefault}$\alpha^2$}}}   
\put(2251,-12706){\makebox(0,0)[lb]{\smash{\SetFigFont{8}{9.6}{\rmdefault}{\bfdefault}{\updefault}$\beta$}}}   
\put(3061,-12706){\makebox(0,0)[lb]{\smash{\SetFigFont{8}{9.6}{\rmdefault}{\bfdefault}{\updefault}$\alpha^2$}}}   
\put(4636,-12706){\makebox(0,0)[lb]{\smash{\SetFigFont{8}{9.6}{\rmdefault}{\bfdefault}{\updefault}$g^{-1}$}}}   
\put(4726,-11581){\makebox(0,0)[lb]{\smash{\SetFigFont{8}{9.6}{\rmdefault}{\mddefault}{\updefault}$16$}}}   
\put(8821,-9961){\makebox(0,0)[lb]{\smash{\SetFigFont{8}{9.6}{\rmdefault}{\bfdefault}{\updefault}$\beta^{-1}$}}}   
\put(8821,-9466){\makebox(0,0)[lb]{\smash{\SetFigFont{8}{9.6}{\rmdefault}{\bfdefault}{\updefault}$\alpha^2$}}}   
\put(3466,-196){\makebox(0,0)[lb]{\smash{\SetFigFont{8}{9.6}{\rmdefault}{\bfdefault}{\updefault}$1$}}}   
\put(2611,-196){\makebox(0,0)[lb]{\smash{\SetFigFont{8}{9.6}{\rmdefault}{\bfdefault}{\updefault}$\alpha^2$}}}   
\put(4861,-196){\makebox(0,0)[lb]{\smash{\SetFigFont{8}{9.6}{\rmdefault}{\bfdefault}{\updefault}$\pi^{\beta}g$}}}   
\put(6526,-196){\makebox(0,0)[lb]{\smash{\SetFigFont{8}{9.6}{\rmdefault}{\bfdefault}{\updefault}$\alpha^2\pi$}}}   
\put(8821,-3346){\makebox(0,0)[lb]{\smash{\SetFigFont{8}{9.6}{\rmdefault}{\bfdefault}{\updefault}$1$}}}   
\put(7066,-12706){\makebox(0,0)[lb]{\smash{\SetFigFont{8}{9.6}{\rmdefault}{\bfdefault}{\updefault}$\alpha^2$}}}   
\put(6301,-12706){\makebox(0,0)[lb]{\smash{\SetFigFont{8}{9.6}{\rmdefault}{\bfdefault}{\updefault}$\pi^{\beta}g$}}}   
\put(541,-6451){\makebox(0,0)[lb]{\smash{\SetFigFont{8}{9.6}{\rmdefault}{\bfdefault}{\updefault}$g^{-1}$}}}   
\put(2521,-1816){\makebox(0,0)[lb]{\smash{\SetFigFont{8}{9.6}{\rmdefault}{\mddefault}{\updefault}$P^{-1}P^{-1}$}}}   
\put(2836,-1321){\makebox(0,0)[lb]{\smash{\SetFigFont{8}{9.6}{\rmdefault}{\mddefault}{\updefault}$7/8$}}}   
\put(6706,-1321){\makebox(0,0)[lb]{\smash{\SetFigFont{8}{9.6}{\rmdefault}{\mddefault}{\updefault}$4/5$}}}   
\put(1936,-9466){\makebox(0,0)[lb]{\smash{\SetFigFont{8}{9.6}{\rmdefault}{\mddefault}{\updefault}$P^{-1}$}}}   
\put(1936,-9871){\makebox(0,0)[lb]{\smash{\SetFigFont{8}{9.6}{\rmdefault}{\mddefault}{\updefault}$P^{-1}$}}}   
\put(6841,-11086){\makebox(0,0)[lb]{\smash{\SetFigFont{8}{9.6}{\rmdefault}{\mddefault}{\updefault}$P^{-1}$}}}   
\put(6346,-11086){\makebox(0,0)[lb]{\smash{\SetFigFont{8}{9.6}{\rmdefault}{\mddefault}{\updefault}$P^{-1}$}}}   
\put(6526,-1771){\makebox(0,0)[lb]{\smash{\SetFigFont{8}{9.6}{\rmdefault}{\mddefault}{\updefault}$P$}}}   
\put(6931,-1771){\makebox(0,0)[lb]{\smash{\SetFigFont{8}{9.6}{\rmdefault}{\mddefault}{\updefault}$P$}}}   
\put(7336,-3031){\makebox(0,0)[lb]{\smash{\SetFigFont{8}{9.6}{\rmdefault}{\mddefault}{\updefault}$P^{-1}$}}}   
\put(7336,-3346){\makebox(0,0)[lb]{\smash{\SetFigFont{8}{9.6}{\rmdefault}{\mddefault}{\updefault}$P^{-1}$}}}   
\put(1576,-13606){\makebox(0,0)[lb]{\smash{\SetFigFont{10}{12.0}{\rmdefault}{\bfdefault}{\updefault}Relation 6 (b): $\alpha^2\beta\alpha^2\pi\pi^{\beta}\alpha^3\beta^2\alpha^2=(\alpha^3\beta^2\alpha^2\beta\alpha^2\pi^{\beta})^2$}}}   
\put(1216,-9691){\makebox(0,0)[lb]{\smash{\SetFigFont{8}{9.6}{\rmdefault}{\mddefault}{\updefault}$12/13$}}}   
\put(4816,-1321){\makebox(0,0)[lb]{\smash{\SetFigFont{8}{9.6}{\rmdefault}{\mddefault}{\updefault}$6$}}}   
\put(8821,-2851){\makebox(0,0)[lb]{\smash{\SetFigFont{8}{9.6}{\rmdefault}{\bfdefault}{\updefault}$\alpha^2$}}}   
\put(6436,-11536){\makebox(0,0)[lb]{\smash{\SetFigFont{8}{9.6}{\rmdefault}{\mddefault}{\updefault}$17/18$}}}   
\put(8101,-9736){\makebox(0,0)[lb]{\smash{\SetFigFont{8}{9.6}{\rmdefault}{\mddefault}{\updefault}$19/20$}}}   
\put(2566,-11581){\makebox(0,0)[lb]{\smash{\SetFigFont{8}{9.6}{\rmdefault}{\mddefault}{\updefault}$14/15$}}}   
\end{picture}

\caption{Relation associated to the cell $(DC2)$}\label{fi10}           
\end{center}           
\end{figure}              
\subsection{The method: K. Brown's theorem for    
groups acting on simply   connected complexes}             
\begin{figure}           
\begin{center}            
\begin{picture}(0,0)%
\includegraphics{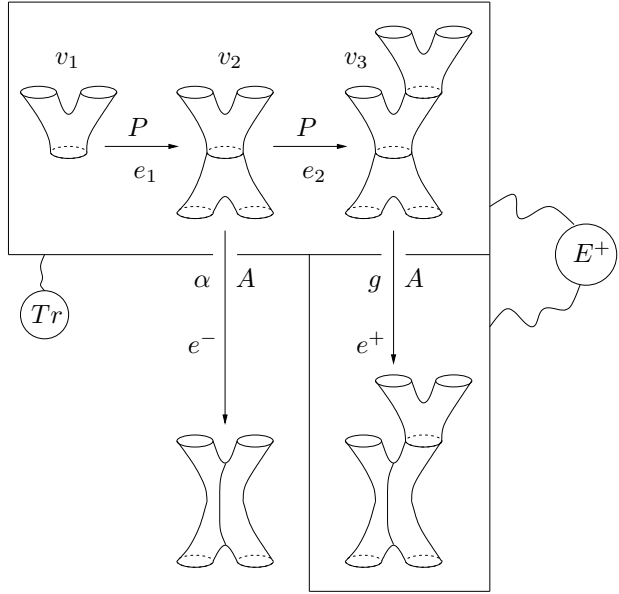}%
\end{picture}%
\setlength{\unitlength}{3315sp}%
\begingroup\makeatletter\ifx\SetFigFont\undefined%
\gdef\SetFigFont#1#2#3#4#5{%
  \reset@font\fontsize{#1}{#2pt}%
  \fontfamily{#3}\fontseries{#4}\fontshape{#5}%
  \selectfont}%
\fi\endgroup%
\begin{picture}(4569,4434)(169,-3673)  
\put(541,299){\makebox(0,0)[lb]{\smash{\SetFigFont{10}{12.0}{\rmdefault}{\mddefault}{\updefault}$v_1$}}}  
\put(2701,299){\makebox(0,0)[lb]{\smash{\SetFigFont{10}{12.0}{\rmdefault}{\mddefault}{\updefault}$v_3$}}}  
\put(1756,299){\makebox(0,0)[lb]{\smash{\SetFigFont{10}{12.0}{\rmdefault}{\mddefault}{\updefault}$v_2$}}}  
\put(1576,-1366){\makebox(0,0)[lb]{\smash{\SetFigFont{10}{12.0}{\rmdefault}{\mddefault}{\updefault}$\alpha$}}}  
\put(361,-1636){\makebox(0,0)[lb]{\smash{\SetFigFont{10}{12.0}{\rmdefault}{\mddefault}{\updefault}$Tr$}}}  
\put(1891,-1366){\makebox(0,0)[lb]{\smash{\SetFigFont{10}{12.0}{\rmdefault}{\mddefault}{\updefault}$A$}}}  
\put(2341,-241){\makebox(0,0)[lb]{\smash{\SetFigFont{10}{12.0}{\rmdefault}{\mddefault}{\updefault}$P$}}}  
\put(2881,-1366){\makebox(0,0)[lb]{\smash{\SetFigFont{10}{12.0}{\rmdefault}{\mddefault}{\updefault}$g$}}}  
\put(3151,-1366){\makebox(0,0)[lb]{\smash{\SetFigFont{10}{12.0}{\rmdefault}{\mddefault}{\updefault}$A$}}}  
\put(4411,-1186){\makebox(0,0)[lb]{\smash{\SetFigFont{10}{12.0}{\rmdefault}{\mddefault}{\updefault}$E^+$}}}  
\put(1531,-1861){\makebox(0,0)[lb]{\smash{\SetFigFont{10}{12.0}{\rmdefault}{\mddefault}{\updefault}$e^-$}}}  
\put(2791,-1861){\makebox(0,0)[lb]{\smash{\SetFigFont{10}{12.0}{\rmdefault}{\mddefault}{\updefault}$e^+$}}}  
\put(1081,-241){\makebox(0,0)[lb]{\smash{\SetFigFont{10}{12.0}{\rmdefault}{\mddefault}{\updefault}$P$}}}  
\put(1126,-556){\makebox(0,0)[lb]{\smash{\SetFigFont{10}{12.0}{\rmdefault}{\mddefault}{\updefault}$e_1$}}}  
\put(2386,-556){\makebox(0,0)[lb]{\smash{\SetFigFont{10}{12.0}{\rmdefault}{\mddefault}{\updefault}$e_2$}}}  
\end{picture}  
  
\caption{Tree $Tr$ and set of edges $E^+$}\label{fi11}           
\end{center}             
\end{figure}

\begin{figure}           
\begin{center}             
\begin{picture}(0,0)%
\includegraphics{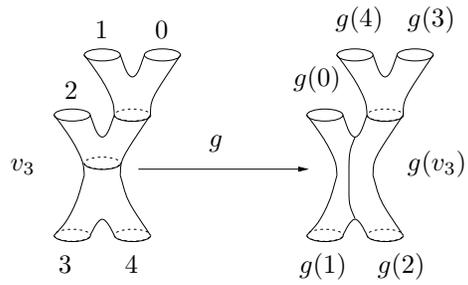}%
\end{picture}%
\setlength{\unitlength}{3315sp}%
\begingroup\makeatletter\ifx\SetFigFont\undefined%
\gdef\SetFigFont#1#2#3#4#5{%
  \reset@font\fontsize{#1}{#2pt}%
  \fontfamily{#3}\fontseries{#4}\fontshape{#5}%
  \selectfont}%
\fi\endgroup%
\begin{picture}(3162,2059)(766,-1334)  
\put(766,-511){\makebox(0,0)[lb]{\smash{\SetFigFont{10}{12.0}{\rmdefault}{\mddefault}{\updefault}$v_3$}}}  
\put(3736,-511){\makebox(0,0)[lb]{\smash{\SetFigFont{10}{12.0}{\rmdefault}{\mddefault}{\updefault}$g(v_3)$}}}  
\put(2251,-331){\makebox(0,0)[lb]{\smash{\SetFigFont{10}{12.0}{\rmdefault}{\mddefault}{\updefault}$g$}}}  
\put(1846,479){\makebox(0,0)[lb]{\smash{\SetFigFont{10}{12.0}{\rmdefault}{\mddefault}{\updefault}$0$}}}  
\put(1396,479){\makebox(0,0)[lb]{\smash{\SetFigFont{10}{12.0}{\rmdefault}{\mddefault}{\updefault}$1$}}}  
\put(1171, 29){\makebox(0,0)[lb]{\smash{\SetFigFont{10}{12.0}{\rmdefault}{\mddefault}{\updefault}$2$}}}  
\put(1621,-1276){\makebox(0,0)[lb]{\smash{\SetFigFont{10}{12.0}{\rmdefault}{\mddefault}{\updefault}$4$}}}  
\put(3511,-1276){\makebox(0,0)[lb]{\smash{\SetFigFont{10}{12.0}{\rmdefault}{\mddefault}{\updefault}$g(2)$}}}  
\put(3736,569){\makebox(0,0)[lb]{\smash{\SetFigFont{10}{12.0}{\rmdefault}{\mddefault}{\updefault}$g(3)$}}}  
\put(3196,569){\makebox(0,0)[lb]{\smash{\SetFigFont{10}{12.0}{\rmdefault}{\mddefault}{\updefault}$g(4)$}}}  
\put(1126,-1276){\makebox(0,0)[lb]{\smash{\SetFigFont{10}{12.0}{\rmdefault}{\mddefault}{\updefault}$3$}}}  
\put(2881,119){\makebox(0,0)[lb]{\smash{\SetFigFont{10}{12.0}{\rmdefault}{\mddefault}{\updefault}$g(0)$}}}  
\put(2926,-1276){\makebox(0,0)[lb]{\smash{\SetFigFont{10}{12.0}{\rmdefault}{\mddefault}{\updefault}$g(1)$}}}  
\end{picture}  
  
  \caption{Homeomorphism $g$}\label{fi12}   
\end{center}          
\end{figure}

Let us denote by $Tr$ a tree of            
representatives of ${\cal DP}_5^+(\s)$ modulo ${\cal B}$, with            
set of vertices $Vert=\{v_1,v_2,v_3\}$, linked by $P$ edges $e_1$ and $e_2$ as indicated in      
Figure \ref{fi11}. We assume that all the representatives are subordinate to   
the canonical rigid structure. We have only displayed the supports of the vertices, so that one has to imagine the canonical   
rigid  structure outside each support.             
Each edge of ${\cal DP}_5^+(\s)$ is equivalent modulo ${\cal            
  B}$ to an edge with origin in $Tr$. There are exactly two   
such edges (modulo      
  $\B$) which do not            
belong to $Tr$:  
\begin{itemize}  
\item  the edge $e^-$: $v_2\stackrel{A}{\longrightarrow} \alpha(v_2)$, and    
\item  the edge $e^+$: $v_3\stackrel{A}{\longrightarrow} g(v_3)$, where $g$ rotates the            
  boundary components of the support of $v_3$ as indicated    
in Figure \ref{fi12}.   
\end{itemize}             
             
\noindent The orientation of the first edge $e^-$ can be inverted by the            
 action of the group $\B$. We denote by $\sigma^-$ the cell            
underlying the oriented edge $e^-$.           
      
\vspace{0.1cm}       
\noindent             
Let $E^+$ be the set of edges of $Tr$ plus the  edge $e^+$, whose orientation cannot be  inverted by $\B$.   
   
\vspace{0.1cm}            
\noindent If $e$ is an edge, denote by $o(e)$ its initial vertex, and    
by $t(e)$ its terminal vertex. For each $e\in E^+$, denote by $w_e$    
the unique vertex of $Tr$            
equivalent to $t(e)$ modulo ${\cal B}$, and choose $g_e\in {\cal B}$    
such that  $t(e)=g_e\cdot w_e$: for instance $g_e=1$ if $e$ is an edge of $Tr$,    
and $g_e=g$ if $e$ is the            
edge $v_3\stackrel{A}{\longrightarrow} g(v_3)$.           
      
 \vspace{0.1cm}       
\noindent            
Theorem 1 of \cite{br} states that  ${\cal B}$ is generated by the    
stabilizer groups of the vertices ${\cal B}_{v_i}$, $i=1,2,3$, the            
stabilizer group ${\cal B}_{\sigma^-}$ of the cell $\sigma^-$, the elements $g_e$,            
$e\in E^+$ (thus, essentially the element $g$), subject to the           
following relations:           
\begin{itemize}           
\item {\it Pres. (i):} For each $e\in E^+$, $g_e^{-1} i_e(h) g_e=c_e(h)$ for      
  all $h\in {\cal B}_e$,            
where $i_e$ is the                         
inclusion ${\cal B}_e\hookrightarrow {\cal B}_{o(e)}$ and $c_e:{\cal B}_e\rightarrow {\cal B}_{w(e)}$            
is the conjugation morphism $h\mapsto  g_e^{-1} h  g_e$.           
           
\item {\it Pres. (ii):}  $i_{e^-}(h)=j_{e^-}(h)$ for $h\in {\cal B}_{e^-}$, where            
  $i_{e^-}:{\cal B}_{e^-}\hookrightarrow {\cal B}_{o(e^-)}$ and $j_{e^-}:{\cal B}_{e^-}\hookrightarrow {\cal B}_{\sigma^-}$            
  are inclusions.           
            
 \item {\it Pres. (iii):}  $r_{\tau}=1$ for each essential $2$-cell $\tau\in {\cal F}$,            
where $r_{\tau}$ is a word in the generators of ${\cal B}_{v_i}$, $\alpha$ and            
$g$, associated with the $2$-cell $\tau$ in the way described in \cite{br}   
(see \S 6.4).           
\end{itemize}   
           
\subsection{Presentations of the stabilizers}           
           
\begin{figure}           
\begin{center}           
\begin{picture}(0,0)%
\includegraphics{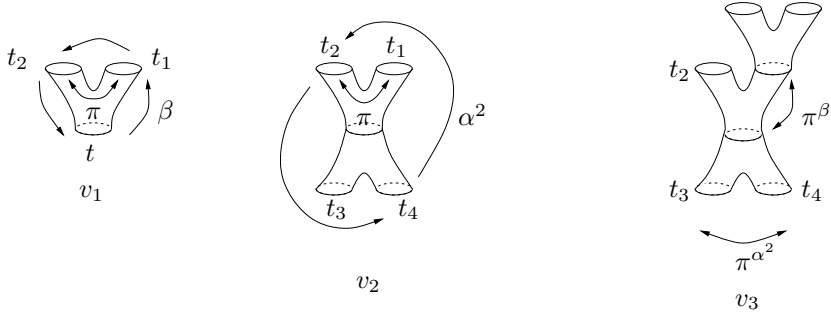}%
\end{picture}%
\setlength{\unitlength}{3315sp}%
\begingroup\makeatletter\ifx\SetFigFont\undefined%
\gdef\SetFigFont#1#2#3#4#5{%
  \reset@font\fontsize{#1}{#2pt}%
  \fontfamily{#3}\fontseries{#4}\fontshape{#5}%
  \selectfont}%
\fi\endgroup%
\begin{picture}(6087,2258)(811,-2041)  
\put(3646,-151){\makebox(0,0)[lb]{\smash{\SetFigFont{10}{12.0}{\rmdefault}{\mddefault}{\updefault}$t_1$}}}  
\put(4186,-691){\makebox(0,0)[lb]{\smash{\SetFigFont{10}{12.0}{\rmdefault}{\mddefault}{\updefault}$\alpha^2$}}}  
\put(811,-241){\makebox(0,0)[lb]{\smash{\SetFigFont{10}{12.0}{\rmdefault}{\mddefault}{\updefault}$t_2$}}}  
\put(6751,-691){\makebox(0,0)[lb]{\smash{\SetFigFont{10}{12.0}{\rmdefault}{\mddefault}{\updefault}$\pi^{\beta}$}}}  
\put(6256,-1771){\makebox(0,0)[lb]{\smash{\SetFigFont{10}{12.0}{\rmdefault}{\mddefault}{\updefault}$\pi^{\alpha^2}$}}}  
\put(6751,-1231){\makebox(0,0)[lb]{\smash{\SetFigFont{10}{12.0}{\rmdefault}{\mddefault}{\updefault}$t_4$}}}  
\put(5761,-331){\makebox(0,0)[lb]{\smash{\SetFigFont{10}{12.0}{\rmdefault}{\mddefault}{\updefault}$t_2$}}}  
\put(6256,-2041){\makebox(0,0)[lb]{\smash{\SetFigFont{10}{12.0}{\rmdefault}{\mddefault}{\updefault}$v_3$}}}  
\put(3151,-151){\makebox(0,0)[lb]{\smash{\SetFigFont{10}{12.0}{\rmdefault}{\mddefault}{\updefault}$t_2$}}}  
\put(3196,-1366){\makebox(0,0)[lb]{\smash{\SetFigFont{10}{12.0}{\rmdefault}{\mddefault}{\updefault}$t_3$}}}  
\put(3736,-1366){\makebox(0,0)[lb]{\smash{\SetFigFont{10}{12.0}{\rmdefault}{\mddefault}{\updefault}$t_4$}}}  
\put(3421,-1906){\makebox(0,0)[lb]{\smash{\SetFigFont{10}{12.0}{\rmdefault}{\mddefault}{\updefault}$v_2$}}}  
\put(1936,-646){\makebox(0,0)[lb]{\smash{\SetFigFont{10}{12.0}{\rmdefault}{\mddefault}{\updefault}$\beta$}}}  
\put(1891,-241){\makebox(0,0)[lb]{\smash{\SetFigFont{10}{12.0}{\rmdefault}{\mddefault}{\updefault}$t_1$}}}  
\put(1351,-1231){\makebox(0,0)[lb]{\smash{\SetFigFont{10}{12.0}{\rmdefault}{\mddefault}{\updefault}$v_1$}}}  
\put(1396,-961){\makebox(0,0)[lb]{\smash{\SetFigFont{10}{12.0}{\rmdefault}{\mddefault}{\updefault}$t$}}}  
\put(1396,-646){\makebox(0,0)[lb]{\smash{\SetFigFont{10}{12.0}{\rmdefault}{\mddefault}{\updefault}$\pi$}}}  
\put(3421,-691){\makebox(0,0)[lb]{\smash{\SetFigFont{10}{12.0}{\rmdefault}{\mddefault}{\updefault}$\pi$}}}  
\put(5761,-1231){\makebox(0,0)[lb]{\smash{\SetFigFont{10}{12.0}{\rmdefault}{\mddefault}{\updefault}$t_3$}}}  
\end{picture}

\caption{Stabilizers of the vertices}\label{fi13}           
\end{center}           
\end{figure}

\noindent{\bf The method:} Each stabilizer group $\B_{v_i}$ ($i=1,2$ or 3)   
permutes the $2i+1$ circles of the pants decomposition of the support of   
$v_i$. Thus, there is a morphism $\B_{v_i}\rightarrow {\cal S}_{2i+1}$, whose   
image is denoted ${\cal P}_i$, and whose kernel, which is generated   
by the Dehn   
twists around the circles of the pants decomposition, is therefore isomorphic to $\Z^{2i+1}$. It follows   
that $\B_{v_i}$ determines an extension   
$$1\rightarrow \Z^{2i+1}\longrightarrow \B_{v_i} \longrightarrow {\cal   
  P}_i\rightarrow 1.$$   
We first find a presentation of ${\cal P}_i$, and use Hall's Lemma   
(cf. \cite{ro}) to deduce a presentation of $\B_{v_i}$. Recall that Hall's Lemma states   
that if a group $N$ is normal in a group $G$, with given presentations $N=\langle   
x_1,\ldots,x_m|r_1=1,\ldots,r_k=1\rangle$ and $G/N=\langle   
\bar{y}_1,\ldots,\bar{y}_n|\bar{s}_1=1,\ldots,\bar{s}_l=1\rangle$, then $G$ is   
generated by $x_1,\ldots,x_m, y_1,\ldots,y_n$ with relations:   
\begin{enumerate}   
\item $r_i=1$, $i=1,\ldots,k$   
\item $s_i=t_i(x)$, $i=1,\ldots,l$, where $t_i(x)$ is a word in $x_1,\ldots,x_m$   
\item $y_j x_i y_j^{-1}=u_{ij}(x)$, $y_j ^{-1} x_i y_j=v_{ij}(x)$,   
  $i=1,\ldots,m$, $j=1,\ldots,n$, where $u_{ij}(x)$ and  $v_{ij}(x)$ are words   
  in $x_1,\ldots,x_m$    
\end{enumerate}

\begin{proposition}             
The stabilizer ${\cal B}_{v_1}$ is generated by the rotation $\beta$, the      
braiding $\pi$ and the Dehn twist $t$, subject to the following relations:            
\begin{enumerate}            
\item $[t,t_i]=[t,\pi]=1$, $i=1,2$, where $t_1=\beta t \beta^{-1}$,            
  $t_2=\beta^{-1} t \beta$ (cf. Rel. {\it 1. (a)})            
\item  $t_1\pi=\pi t_2$, $t_2\pi=\pi t_1$ (cf. Rel. {\it 1. (b)})            
\item $\pi^2=t t_1^{-1}t_2^{-1} $ (cf. Rel. {\it 1. (c)})            
\item $\beta^3=1$ (cf. Rel. {\it 1. (d)})          
\item $\beta=t\pi^{\beta}\pi$ (cf. Rel. {\it 1. (e)})           
\end{enumerate}             
Denote by $\alpha^2$ the rotation which interchanges the adjacent pairs of pants            
of the support of $v_2$, by $\pi,t_1,t_2$ respectively the braiding and the Dehn twists as above.              
The stabilizer ${\cal B}_{v_2}$ is generated by $\alpha^2$, $\pi$, $t_1$ and            
$t_2$, subject to the following relations:              
\begin{enumerate}            
\item $(\alpha^2)^2=1$ (cf. {\it Rel. 3. (a)})            
\item $[\pi,\alpha^2\pi\alpha^2]=1$ (cf. {\it Rel.  4. (a)})            
\item $t_1\pi=\pi t_2$, $t_2\pi=\pi t_1$, $[t_1,\alpha^2\pi\alpha^2]=1$,   
  $[t_2,\alpha^2\pi\alpha^2]=1$ (cf. {\it Rel.  1. (b)}, {\it Rel.  4. (b)})            
\item $[t_1,t_2]=1$ (redundant in the presentation of $\B$), and if $t_3=\alpha^2 t_1\alpha^2$, $t_4=\alpha^2  t_2\alpha^2$, then $[t_1,t_3]=1$ (cf. {\it Rel. 4. (d)}), $[t_1,t_4]=[t_2,t_3]=[t_2,t_4]=1$   
  (redundant in the presentation of $\B$)            
 \item $\pi^2t_1t_2= \alpha^2 \pi^2t_1t_2  \alpha^2   $ ({\it Rel. 5. (a)}, using {\it Rel. 1. (c)} )            
\end{enumerate}    
            
The stabilizer ${\cal B}_{v_3}$ is generated by $\pi^{\beta}$, which is the            
natural extension of the braiding $\beta\pi\beta^{-1}$ of ${\cal B}_{v_1}$, by            
$\pi^{\alpha^2}$ which corresponds to the element $\alpha^2\pi\alpha^2$ in            
${\cal B}_{v_2}$, and by the Dehn twists $t_1, t_2,t_3,t_4$ and $t$, subject to the   
following relations:   
   
\begin{enumerate}   
\item $(\pi^{\alpha^2})^2=tt_3^{-1}t_4^{-1}$ (cf. {\it Rel. 1. (c)}, after   
    conjugation by $\alpha^2$)    
\item $(\pi^{\beta})^2=t_2t^{-1}t_1^{-1}$ (cf. {\it Rel. 1. (c)}, after   
    conjugation by $\beta$)    
\item $[\pi^{\alpha^2},\pi^{\beta}\pi^{\alpha^2}(\pi^{\beta})^{-1}]=1$ (cf. {\it Rel. 4. (e)})   
\item $t_3\pi^{\alpha^2}=\pi^{\alpha^2}t_4$,   
  $t_4\pi^{\alpha^2}=\pi^{\alpha^2}t_3$,   
  $[\pi^{\alpha^2},t]=1$ (cf. {\it Rel. 1. (a), (b)}, after   
    conjugation by $\alpha^2$), $[\pi^{\alpha^2},t_1]=1$ (cf. {\it Rel. 4. (b)}, after   
    conjugation by $\alpha^2$),   
  $[\pi^{\alpha^2},(\pi^{\beta})^{-1}t_3\pi^{\beta}]=1$ (cf. {\it   
    Rel. 4. (d)}), $[\pi^{\alpha^2},(\pi^{\beta})^{-1}t_4\pi^{\beta}]=1$   
  (redundant in the presentation of $\B$)   
\item $\pi^{\beta}t=t_1\pi^{\beta}$, $\pi^{\beta}t_1=t\pi^{\beta}$,   
  $[\pi^{\beta}, t_2]=1$ (cf. {\it   
    Rel. 1. (a), (b)}, after conjugation by $\beta$)   
\item $[t_2,t_3]=[t_2,t_1]=1$ (redundant in the presentation of $\B$),   
  $ [t_3,t]=[t_3,t_1]=[t_3,\pi^{\beta}t_3(\pi^{\beta})^{-1}]=1$ (cf. {\it   
    Rel. 4. (d)}), $[t,t_1]=1$ (cf. {\it Rel. 1. (a)})   
\end{enumerate}

The stabilizer ${\cal B}_{\sigma^-}$ of the cell $\sigma^-$ is generated by      
$\alpha, t_1,t_2,t_3$ and $t_4$, subject to the relations:      
\begin{enumerate}      
\item $\alpha^4=1$ (cf. {\it Rel. 3. (a)})      
\item $\alpha t_i\alpha^{-1}=t_{i+1}$ ($\forall\; i\; {\rm mod}\; 4$) (cf. {\it    
    Rel. 4. (c)})      
\item $[t_i,t_j]=1$  ($\forall\; i,j\;{\rm mod}\; 4$) (cf. {\it Rel. 4. (d)})      
\end{enumerate}      
Equivalently, ${\cal B}_{\sigma^-}\cong \langle\alpha,t_1|\; \alpha^4=1,\; [t_1,\alpha t_1\alpha^{-1}]=1,\;[t_1,\alpha^2 t_1\alpha^2]=1\rangle.$       
           
\end{proposition}             
\begin{proof}   
   
\begin{enumerate}   
   
\item $\B_{v_1}$ is precisely ${\cal M}(0,3)$. Its presentation is   
  well known. It may be deduced from the presentation of ${\cal S}_3$ as   
  ${\cal S}_3=\langle  \beta,\pi|\beta^3=1,\; \pi^2=1,\;   
  \beta=\pi^{\beta}\pi\rangle$. We note that the relation $[t_1,t_2]=1$ has   
  been deliberately omitted, as it is a consequence of the others.   
       
\item The permutation group ${\cal P}_2$ is isomorphic to the semi-direct   
  product $(\Z/2\Z\oplus \Z/2\Z)\bowtie \Z/2\Z$. The generators of the three   
  factors are the images of   
  $\pi, \alpha^2\pi\alpha^2$ and $\alpha^2$, respectively. The presentation of ${\cal P}_2$   
  is therefore ${\cal P}_2=\langle \pi,\alpha^2| \pi^2=1,\;(\alpha^2)^2=1,\;   
  [\pi, \alpha^2\pi\alpha^2]=1\rangle$. One applies Hall's Lemma and obtains a   
  presentation of $\B_{v_2}$ with generators $\alpha^2$, $\pi$,   
  $t_1,t_2,t_3,t_4$ and $t$, and relations:    
   
\begin{itemize}   
          
\item $(\alpha^2)^2=1$    
\item ($\pi^2=tt_1^{-1}t_2^{-1}$), from which one can eliminate   
the generator $t$            
\item $[\pi,\alpha^2\pi\alpha^2]=1$             
\item $t_1\pi=\pi t_2$, $t_2\pi=\pi t_1$, ($[\pi,t]=1$),   
  $[t_1,\alpha^2\pi\alpha^2]=1$ (equivalent to $[t_3,\pi]=1$),   
  $[t_2,\alpha^2\pi\alpha^2]=1$ (equivalent to $[t_4,\pi]=1$)            
\item $t_3=\alpha^2 t_1\alpha^2$, $t_4=\alpha^2 t_2\alpha^2$, $t=\alpha^2   
  t\alpha^2$   
\item The 10 commutation relations between the 5 Dehn twists are:\\   
$[t_1,t_2]=1$, ($[t,t_1]=1$), ($[t,t_2]=1$), $[t_1,t_3]=1$,   
  $[t_2,t_3]=1$, $[t_1,t_4]=1$, $[t_2,t_4]=1$ (and the conjugates by   
  $\alpha^2$ of the first three:   
  $[t_3,t_4]=[t,t_3]=[t,t_4]=1$)            
\end{itemize}   
            
The parentheses indicate redundant relations, and omitting them provides the   
presentation of $\B_{v_2}$.

\item  The permutation group ${\cal P}_3$ is isomorphic to the semi-direct   
  product $(\Z/2\Z\oplus \Z/2\Z)\bowtie \Z/2\Z$. The generators of the three   
  factors are the images of   
  $\pi^{\alpha^2}$, $\pi^{\beta}\pi^{\alpha^2}\pi^{\beta}$ and $\pi^{\beta}$,   
  respectively. The presentation of ${\cal P}_3$   
  is therefore ${\cal P}_3=\langle \pi^{\alpha^2}, \pi^{\beta}|   
  (\pi^{\alpha^2})^2=1,\;(\pi^{\beta})^2=1,\;  [\pi^{\alpha^2},   
  \pi^{\beta}\pi^{\alpha^2}\pi^{\beta}]=1\rangle$. One applies Hall's Lemma and obtains a  presentation of $\B_{v_3}$ with generators $\pi^{\alpha^2}$, $\pi^{\beta}$,   
  $t_1,t_2,t_3,t_4$ and $t$, and relations:    
\begin{itemize}   
\item $(\pi^{\alpha^2})^2=tt_3^{-1}t_4^{-1}$   
\item $(\pi^{\beta})^2=t_2t^{-1}t_1^{-1}$   
\item $[\pi^{\alpha^2},\pi^{\beta}\pi^{\alpha^2}(\pi^{\beta})^{-1}]=1$   
\item $t_3\pi^{\alpha^2}=\pi^{\alpha^2}t_4$,   
  $t_4\pi^{\alpha^2}=\pi^{\alpha^2}t_3$,   
  $[\pi^{\alpha^2},t]=[\pi^{\alpha^2},t_1]=$   
  $[\pi^{\alpha^2},(\pi^{\beta})^{-1}t_3\pi^{\beta}]=[\pi^{\alpha^2},(\pi^{\beta})^{-1}t_4\pi^{\beta}]$ $=1$   
\item $\pi^{\beta}t=t_1\pi^{\beta}$, $\pi^{\beta}t_1=t\pi^{\beta}$,   
  $[\pi^{\beta}, t_2]=1$   
\item The 21 commutation relations between the 7 Dehn twists are:\\   
  $[t_2,t_3]=([t_2,t_4])=[t_2,t_1]=([t_2,t])=([t_2,\pi^{\beta}t_3(\pi^{\beta})^{-1}])=([t_2,\pi^{\beta}t_4(\pi^{\beta})^{-1}])=1$\\   
$([t_3,t_4])=[t_3,t]=[t_3,t_1]=[t_3,\pi^{\beta}t_3(\pi^{\beta})^{-1}]=([t_3,\pi^{\beta}t_4(\pi^{\beta})^{-1}])=1$\\   
$([t_4,t])=([t_4,t_1])=([t_4,\pi^{\beta}t_3(\pi^{\beta})^{-1}])=([t_4,\pi^{\beta}t_4(\pi^{\beta})^{-1}])=1$\\   
$[t,t_1]=([t,\pi^{\beta}t_3(\pi^{\beta})^{-1}])=([t,\pi^{\beta}t_4(\pi^{\beta})^{-1}])$=1\\   
$([t_1,\pi^{\beta}t_3(\pi^{\beta})^{-1}])=([t_1,\pi^{\beta}t_4(\pi^{\beta})^{-1}])=1$\\   
  $([\pi^{\beta}t_3(\pi^{\beta})^{-1},\pi^{\beta}t_4(\pi^{\beta})^{-1}])=1$   
\end{itemize}   
   
The commutations relations in parentheses may be deduced from the others by   
conjugation by $\pi^{\beta}$ or $\pi^{\alpha^2}$.

\item $\B_{\sigma^-}$ is an extension of $\Z/4\Z$ (generated by $\alpha$) by   
  $\Z^4$ (generated by $t_1,t_2,t_3$ and $t_4$).                   
\end{enumerate}

\end{proof}

\subsection{The relations}            
\begin{enumerate}            
\item  By the relations of type  {\it Pres.  (i)} applied to the edges $e$ inside the            
tree $Tr$, we immediately identify the generator $\pi$ of ${\cal  B}_{v_1}$ with the generator of ${\cal  B}_{v_2}$ denoted the same way,            
and the generators $t_1,t_2, t_3,t_4$ and $t$ denoted the same way in the presentations of            
the stabilizers of the vertices.           
\item By the relation of type {\it Pres.  (ii)}, the generator denoted $\alpha^2$ in ${\cal  B}_{v_2}$            
is identified with the square of the generator $\alpha$ of $\B_{\sigma^{-}}$, and the      
generators $t_1,t_2,t_3,t_4$ of $\B_{\sigma^{-}}$ are identified with the      
generators denoted the same way in the stabilizers of the vertices.           
\item  By the relation of type  {\it Pres.  (i)} applied to the edge            
$v_2\longrightarrow g(v_2)$, we obtain the relation            
$g^{-1}\pi^{\beta}\pi^{\alpha^2}(\pi^{\beta})^{-1} g= \pi^{\alpha^2}$ and the      
relations $g^{-1} t_2 g= \pi^{\beta}t_3 (\pi^{\beta})^{-1}$,      
  $g^{-1} t_3 g=\pi^{\beta}t_4 (\pi^{\beta})^{-1}$, $g^{-1} t_4 g=t_2$. They      
  give the relations {\it Rel. 4. (f), (g)}, after $g$ is identified with $(\beta\alpha)^{-1}$.           
\item  For each of the 6 essential cells, we compute an associated relation by   
  the procedure described in \cite{br}. Following closely the exposition of   
  \cite{br}, we recall it for the convenience of the reader:\\     
   
Each edge of the complex starting in $Tr$ has one of the following forms:   
   
\begin{enumerate}   
\item $v_1\stackrel{h(e_1)}{\longlongrightarrow} h(v_2)$, $h\in\B_{v_1},\;\;$   
  $v_2\stackrel{h(e_2)}{\longlongrightarrow} h(v_3)$, $h\in\B_{v_2}$   
\item $v_2\stackrel{h(e^-)}{\longlongrightarrow} h(\alpha(v_2))$, $h\in\B_{v_2}$   
\item $v_3\stackrel{h(e^+)}{\longlongrightarrow} h(g(v_3))$, $h\in\B_{v_3}$   
\item $v_3\stackrel{hg^{-1}(\overline{e^+})}{\longlongrightarrow} h(g^{-1}(v_3))$,   
  $h\in\B_{v_3}$, where $\overline{e^+}$ is the edge obtained from $e^+$ by   
  inverting its orientation.   
\end{enumerate}   
   
To such an edge $e$ we associate an element of $\gamma\in \B$ such that $e$   
ends in $\gamma(Tr)$: $\gamma=h$ in case (a), $\gamma=h\alpha$ in case (b),   
$\gamma=hg$ in case (c), and $\gamma=hg^{-1}$ in case (d).  
  
\vspace{0.1cm}   
\noindent   
Let $\tau$ be one of the 6 cells. One chooses an orientation and a cyclic   
labeling of the boundary edges, such that the labeled 1 edge $a_1$  starts from   
a vertex $v$ of the tree $Tr$.  
  
\vspace{0.1cm}   
\noindent  
Let $\gamma_1$   
be associated to $a_1$ as above. It ends in $\gamma_1(Tr)$, so the second edge   
is of the form $\gamma_1(a_2)$ for some edge $a_2$ starting in $Tr$. Let   
$\gamma_2$ be associated to $a_2$. The second edge ends in   
$\gamma_1\gamma_2(Tr)$. If $n$ is the length of the cycle bounding $\tau$, one obtains this way a sequence $\gamma_1,\ldots,\gamma_n$ such that   
$\gamma_1\cdots\gamma_n(v)=v$.  
  
\vspace{0.1cm}   
\noindent  
Note that for each of the 6 cycles, we have   
indicated the corresponding $\gamma_i$ above the $i^{th}$ edge.  
  
\vspace{0.1cm}   
\noindent  
Let $\gamma$ be the element of the stabilizer $\B_v$ computed as   
$\gamma_1\cdots\gamma_n$ when each element $\gamma_i$ is viewed in $\B$. Then   
the relation associated to $\tau$ is   
$$ \gamma_1\cdots\gamma_n=\gamma$$   
where the left-hand side is viewed as a word in $g,\alpha$, and elements of   
$\B_{v_i}$ ($i=1,2,3$).  
  
\vspace{0.1cm}   
\noindent  
We finally give the corresponding 6 relations:    
   
\begin{enumerate}            
\item  Cell $PP=PP$: it gives the relation $\pi^{\beta}=\beta^{-1}\pi\beta$ in            
${\cal B}$.            
\item Cell $PA=AP$: its associated relation identifies the word            
$g^{-1}\alpha^{-1}\pi^{-1}$ on the generators            
$\alpha$, $\pi$ and $g$, with the element of $\beta\pi^{-1}$ of the stabilizer            
${\cal B}_{v_3}$. As can be guessed, this will give the relation            
$g=(\beta\alpha)^{-1}$, but this is more tricky than it seems:  
$\pi\beta^{-1}$ (in ${\cal B}_{v_3}$) = $\pi\alpha g$ $\iff  
\pi^{\beta}\pi\beta^{-1}$  (in ${\cal B}_{v_3}$)= $\pi^{\beta}\pi\alpha g$. Now  
$\pi^{\beta}\pi\beta^{-1}$  (in ${\cal B}_{v_3}$) = $t^{-1}$ (in ${\cal  
  B}_{v_3}$), and by  {\it Pres.  (i)}, $t^{-1}$ (in ${\cal  B}_{v_3}$)= $t^{-1}$  
(in ${\cal  B}_{v_1}$). In ${\cal  B}_{v_1}$, $t^{-1}=\beta^{-1}  
\pi{\beta}\pi\beta^{-1}$, and using the previous relation  
$\pi^{\beta}=\beta^{-1}\pi\beta$, one obtains $\beta^{-1}\pi{\beta}\pi\beta^{-1}$ (in  
${\cal  B}_{v_1}$)=$\beta^{-1}\pi{\beta}\pi\alpha g$, hence  $g=(\beta\alpha)^{-1}$.           
\item  The remaining cells (Hatcher-Thurston type cells) give relations {\it            
  2. (a), (b)} and {\it 6. (a), (b)}.                    
            
\end{enumerate}          
\end{enumerate}

\section{The group $\B$ and the braided Thompson group of M. Brin}      
      
\begin{definition}\label{root}      
\begin{enumerate}      
\item Let $c$ denote the boundary circle of the support $S_0$ labeled 3 on Figure      
\ref{t}. The rooted dyadic surface $\s^r$ is the closure of the connected      
component of $\s\setminus c$ that contains $S_0$.         
\item A rooted admissible surface $\Sigma_{0,n+1}^r$ of level $n+1$ is an      
  admissible surface of $\s$ containing $S_0$ and contained in $\s^r$. It is      
  endowed with a cyclic labeling of its boundary circles by $1,\ldots,n+1$,      
  in such a way that $c$ corresponds to the label $n+1$.      
\end{enumerate}      
\end{definition}      
      
\begin{lemma}\label{tresse}      
Let $\Sigma^r_{0,n+1}$ be an admissible rooted surface of level $n+1$. If      
$n\geq 2$, there      
is a canonical embedding $\iota_{\Sigma^r_{0,n+1}}:B_n\longrightarrow {\cal      
  M}(\Sigma^r_{0,n+1})$, where $B_n$ is the Artin braid group.      
\end{lemma}      
      
\begin{proof}      
Let $\sigma_1,\ldots,\sigma_{n-1}$ denote the standard generators of      
$B_n$. For $i=1,\ldots,n$, denote by $c_i$ the $i^{th}$ boundary circle of      
$\Sigma^r_{0,n+1}$. For $i=1,\ldots,n-1$, there exists a unique (up to      
isotopy) pair of pants      
$P_i$ containing $c_i$ and      
$c_{i+1}$ in its boundary, which is homeomorphic to $S_0$ by a homeomorphism      
representing a class of $T$. Let $g_i\in T$ such that      
$g_i(S_0)=P_i$, which maps the circle labeled 1 (resp. labeled 2) of $S_0$ (see      
Figure \ref{t} or \ref{pi}) on the circle $c_i$ (resp. $c_{i+1}$) of $P_i$. One      
defines    
$$\iota_{\Sigma^r_{0,n+1}}(\sigma_i)=g_i \pi g_i^{-1}$$      
The result does not depend on the choice of $g_i$, since another choice $g_i'$      
should coincide with $g_i$ on $S_0$.   
   
\vspace{0.1cm}     
\noindent    
In fact, the embedding we have constructed $\iota_{\Sigma^r_{0,n+1}}:B_n\longrightarrow {\cal      
  M}(\Sigma^r_{0,n+1})$ is quite standard, and we have only expressed it using      
mapping classes in $\B$. By abuse of notation, we will denote by      
$\sigma_1,\ldots,\sigma_{n-1}$ the images of the Artin generators in ${\cal      
  M}(\Sigma^r_{0,n+1})$.      
\end{proof}

\begin{figure}    
\begin{center}      
\begin{picture}(0,0)%
\includegraphics{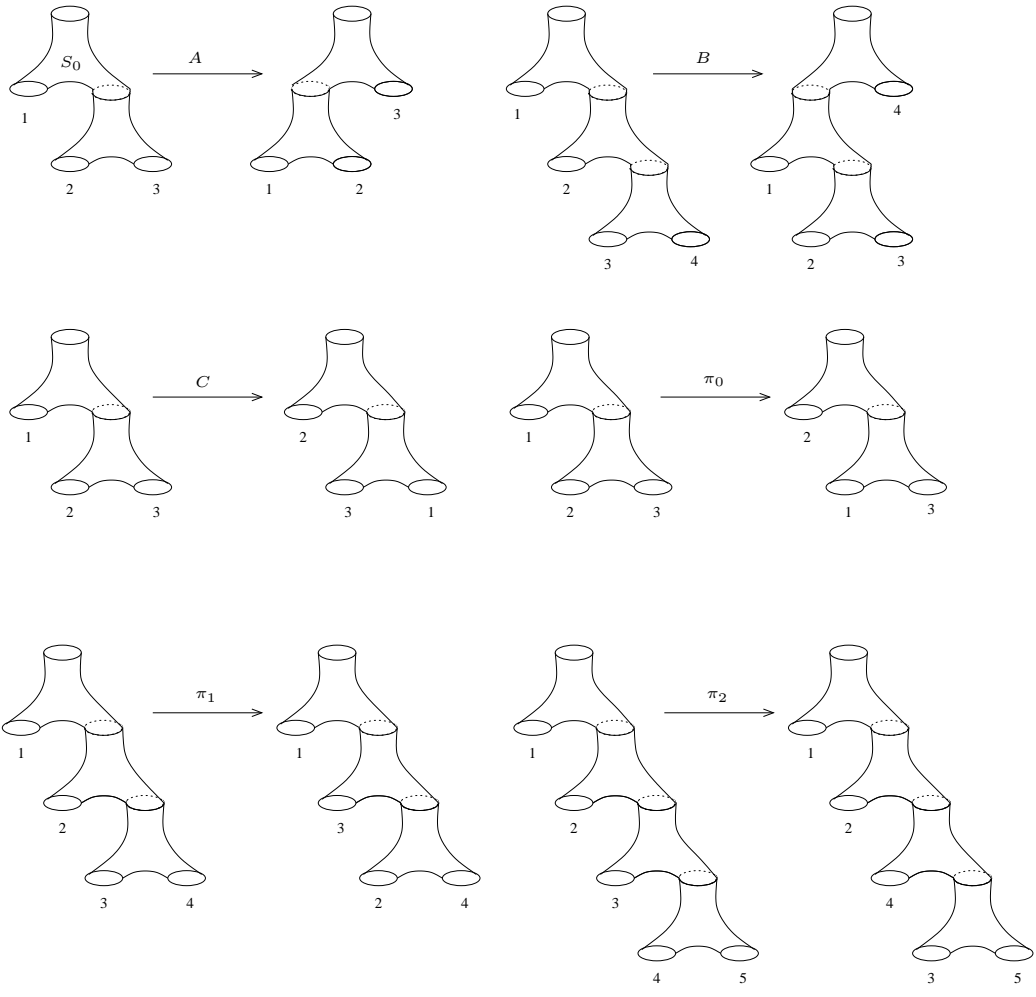}%
\end{picture}%
\setlength{\unitlength}{2072sp}%
\begingroup\makeatletter\ifx\SetFigFont\undefined%
\gdef\SetFigFont#1#2#3#4#5{%
  \reset@font\fontsize{#1}{#2pt}%
  \fontfamily{#3}\fontseries{#4}\fontshape{#5}%
  \selectfont}%
\fi\endgroup%
\begin{picture}(12346,11707)(443,-11311)  
\put(2791,-4156){\makebox(0,0)[lb]{\smash{\SetFigFont{6}{7.2}{\rmdefault}{\mddefault}{\updefault}$C$}}}  
\put(2701,-286){\makebox(0,0)[lb]{\smash{\SetFigFont{6}{7.2}{\rmdefault}{\mddefault}{\updefault}$A$}}}  
\put(1171,-331){\makebox(0,0)[lb]{\smash{\SetFigFont{6}{7.2}{\rmdefault}{\mddefault}{\updefault}$S_0$}}}  
\put(8776,-286){\makebox(0,0)[lb]{\smash{\SetFigFont{6}{7.2}{\rmdefault}{\mddefault}{\updefault}$B$}}}  
\put(8911,-7891){\makebox(0,0)[lb]{\smash{\SetFigFont{6}{7.2}{\rmdefault}{\mddefault}{\updefault}$\pi_2$}}}  
\put(2791,-7891){\makebox(0,0)[lb]{\smash{\SetFigFont{6}{7.2}{\rmdefault}{\mddefault}{\updefault}$\pi_1$}}}  
\put(8866,-4111){\makebox(0,0)[lb]{\smash{\SetFigFont{6}{7.2}{\rmdefault}{\mddefault}{\updefault}$\pi_0$}}}  
\end{picture}

\caption{Elements of the braided Thompson group}\label{brinmeier}      
\end{center}      
\end{figure}

\noindent We introduce a subgroup of $\B$ which has been recently   
introduced and studied by    
M. Brin in \cite{bri}:        
      
\begin{definition}      
The braided Thompson group $B{\bf V}$ is the subgroup of $\B$ generated by      
the four mapping classes $A,B,C$ and $\pi_0$ (cf. Figure \ref{brinmeier})   
 represented by homeomorphisms      
supported on $\s^r$ (fixing pointwise the root circle $c$):   
\begin{itemize}    
\item  $A$ and $B$ are the elements of $T\subset\B$ represented      
on Figure \ref{brinmeier} by symbols of the form $\Sigma_l\rightarrow      
\Sigma_r$, where $\Sigma_l$ and $\Sigma_r$ are admissible rooted   
surfaces. The      
left surface $\Sigma_l$ is canonically labeled (see Definition \ref{root}),      
while a circle labeled $i$ of the boundary of $\Sigma_r$ is the image of the      
$i^{th}$ boundary circle of $\Sigma_l$. Contrary to the conventions of      
Figures \ref{t} to \ref{F}, we have not represented on $\Sigma_l$ the change      
of rigid structure induced by the mapping class.   
\item  The mapping class $C$ is represented by a symbol      
$\Sigma_C\rightarrow\Sigma_C$ (with the conventions as above), where $\Sigma_C$      
is a rooted admissible surface of level 4. Via the embeddings      
$$B_3\stackrel{\iota_{\Sigma_C}}{\hookrightarrow}{\cal      
  M}(\Sigma_C)\hookrightarrow \B$$      
$C$ corresponds to $\sigma_2\sigma_1$ (see Lemma \ref{tresse}).     
\item  The mapping class $\pi_0$ corresponds to the braiding $\sigma_1$ occurring      
in the above definition of $C$.    
\end{itemize}        
\end{definition}

\begin{proposition}      
\begin{enumerate}      
\item Denote the natural images of $A,B,C$ and $\pi_0$ in $V$ by ${\bf A},      
  {\bf B}, {\bf C}$ and ${\bf \pi_0}$, respectively. They generate a subgroup      
  ${\bf V}$ of $V$ isomorphic to $V$.      
\item For an admissible rooted surface $\Sigma$, let $K(\Sigma)$ denote the image      
  by $\iota_{\Sigma}$ in ${\cal M}(\Sigma)$ of the Artin group of pure      
  braids, and $K_{\infty}$ be the      
  inductive limit $\displaystyle{\lim_{\stackrel{\rightarrow}{\Sigma}}      
  K(\Sigma)}$ induced by the inclusions $\Sigma\subset \Sigma'$. There is a short exact sequence      
$$1\rightarrow K_{\infty}\longrightarrow B{\bf V}\longrightarrow {\bf      
  V}\rightarrow 1$$      
which splits over the subgroup ${\bf F}$ of ${\bf V}$ generated by ${\bf A}$      
and ${\bf B}$.      
\end{enumerate}      
\end{proposition}          
      
\begin{proof}    
\begin{enumerate}     
\item The standard way to present the Thompson group acting on the Cantor set is      
precisely as the group ${\bf V}$ generated by ${\bf A},      
  {\bf B}, {\bf C}$ and ${\bf \pi_0}$ (see \cite{CFP}).     
     
\item  Since $A$ and $B$ belong to $T$, the extension splits over ${\bf F}$. It      
remains to prove that the kernel of the projection $B{\bf V}\longrightarrow      
{\bf  V}$ is $K_{\infty}$. Clearly, the kernel must be contained in      
$K_{\infty}$, and it suffices to check that $B{\bf V}$ contains $K_{\infty}$. For an integer $n\geq      
1$, set $C_n=A^{-n+1}CB^{n-1}$. Define      
$\pi_1=C_2^{-1}\pi_0C_2$, and $\pi_n=A^{-n+1}\pi_1A^{n-1}$ for $n\geq 2$. By      
an easy computation, $\pi_1=\iota_{\Sigma_{\pi_1}}(\sigma_2)$, where      
$\Sigma_{\pi_1}$ is the support of $\pi_1$ represented on Figure      
\ref{brinmeier}. Conjugating by $A^{-n+1}$, one immediately obtains that      
$\pi_n=\iota_{\Sigma_{\pi_n}}(\sigma_{n+1})$, where $\Sigma_{\pi_n}$, the      
support of $\pi_n$, has level $n+4$, and $\iota_{\Sigma_{\pi_n}}$ embeds      
$B_{n+3}$ into ${\cal M}(\Sigma_{\pi_n})$. It follows that $B{\bf V}$      
contains sufficiently many standard Artin generators (namely $\pi_0,   
C\pi_1^{-1},\pi_1,\ldots,\pi_n,\ldots$ and all their conjugates by    
$F=\langle A,B\rangle$) to contain the images of the Artin braid      
groups in ${\cal M}(\Sigma)$ for any rooted admissible surface $\Sigma$. In      
particular, $B{\bf V}$ contains $K_{\infty}$.            
\end{enumerate}   
\end{proof}

\bibliographystyle{plain}

\end{document}